\documentclass[aop,MSNbibl,number,citesort,seceqn,dvips]{arximspdf}
\usepackage{accents}
\usepackage{graphicx}

\doi{10.1214/11-AOP707}
\volume{41}
\issue{1}
\pubyear{2013}
\firstpage{50}
\lastpage{108}

\makeatletter
\newtheorem{theo}{Theorem}[section]
\newtheorem{coro}{Corollary}[section]
\newtheorem{prop}{Proposition}[section]
\newtheorem{lemm}{Lemma}[section]
\newproclaim{defn}{Definition}[section]
\newproclaim{rema}{Remark}
\newproclaim{remark}{Remark}

\newcommand{\cal}{\mathcal}
\def\N{\mathbb{N}}
\def\E{\mathbb{E}}
\def\S{{\cal S}}
\def\0{{\bf0}}
\def\Z{\mathbb{Z}}
\def\R{\mathbb{R}}
\def\ext{\operatorname{ext}}
\def\B{{ B}}
\renewcommand{\E}{{\mathbb E}}
\newcommand{\C}{{\cal C}}
\newcommand{\tod}{\stackrel{{\cal D}}{\longrightarrow}}
\newcommand{\eqd}{\stackrel{{\cal D}}{=}}

\newcommand{\X}{{\cal X}}
\def\la{{\lambda}}
\renewcommand{\P}{{{\cal P}}}
\newcommand{\Cov}{\operatorname{Cov}}
\newcommand{\Var}{\operatorname{Var}}
\newcommand{\conv}{\operatorname{conv}}
\newcommand{\Faces}{{\cal F}}
\newcommand{\Span}{\operatorname{aff}}
\newcommand{\phull}{\operatorname{p\mbox{-}hull}}
\newcommand{\uphull}{\operatorname{up\mbox{-}hull}}
\newcommand{\Vertices}{\operatorname{Vertices}}
\newcommand{\closure}{\operatorname{cl}}
\newcommand{\Top}{\operatorname{Top}}
\newcommand{\cone}{\operatorname{cone}}
\newcommand{\vcone}{\mbox{v-cone}}
\newcommand{\extdir}{\mbox{ext-dir}}
\newcommand{\Vol}{\operatorname{Vol}}
\newcommand{\F}{{\cal F}}
\newcommand{\1}{{\mathbf1}}
\def\R{\mathbb{R}}
\def\B^2{\mathbb{D}}
\def\S{\mathbb{S}}
\def\B{\mathbb{B}}
\def\la{{\lambda}}
\newcommand{\eqref}[1]{(\ref{#1})}
\newproclaim{remarks}{Remarks}
\newcommand{\cdotsm}{\cdots}
\renewcommand{\epsilon}{\varepsilon}
\newcommand{\fraca}[2]{{#1}/{#2}}
\newcommand{\fracd}[2]{({#1}/{#2})}
\makeatother

\begin{document}
\begin{frontmatter}

\title{Brownian limits, local limits and variance asymptotics for convex hulls in the ball}
\runtitle{Convex hulls in the ball}

\begin{aug}
\author[A]{\fnms{Pierre} \snm{Calka}\ead[label=e1]{pierre.calka@univ-rouen.fr}},
\author[B]{\fnms{Tomasz} \snm{Schreiber}\thanksref{t1,au2}}
\and
\author[C]{\fnms{J. E.} \snm{Yukich}\corref{}\thanksref{au3}\ead[label=e3]{joseph.yukich@lehigh.edu}}
\runauthor{P. Calka, T. Schreiber and J. E. Yukich}
\affiliation{Universit\'e de Rouen, Nicholas Copernicus University and Lehigh~University}
\dedicated{Dedicated to the memory of Tomasz Schreiber}
\address[A]{P. Calka\\
Laboratoire de Math\'ematiques Rapha\"el Salem\\
Universit\'e de Rouen\\
Avenue de l'Universit\'e, BP.12\\ \quad
Technop\^ole du Madrillet\\
F76801 Saint-Etienne-du-Rouvray\\
France\\
\printead{e1}} 
\address[B]{T. Schreiber\\
Faculty of Mathematics\\ \quad  and Computer Science\\
Nicholas Copernicus University\\ Toru\'n\\
 Poland}
\address[C]{J. E. Yukich\\ Department of Mathematics\\ Lehigh
University\\
Bethlehem, Pennsylvania 18015\\USA\\
\printead{e3}}
\end{aug}
\thankstext{t1}{Born June 25, 1975; died on December 1, 2010.}
\thankstext{au2}{Supported in part by the Polish
Minister of Science and Higher Education
Grant N N201 385234 (2008--2010).}
\thankstext{au3}{Supported in part by NSF Grant
DMS-08-05570.}

\received{\smonth{12} \syear{2009}}
\revised{\smonth{4} \syear{2011}}

%
\begin{abstract}
Schreiber and Yukich [\textit{Ann. Probab.} \textbf{36} (2008) 363--396]
establish an asymptotic representation for random convex
polytope geometry in the unit ball $\B^d,  d \geq2$, in terms of the
general theory
of stabilizing functionals of Poisson point processes as well as in
terms of \textit{generalized paraboloid growth processes}.
This paper further exploits this connection, introducing also a dual object
termed the \textit{paraboloid hull process}. Via these growth processes we
establish local functional
limit theorems for the properly scaled radius-vector and support functions
of convex polytopes
generated by high-density Poisson samples. We show that direct methods
lead to explicit
asymptotic expressions for the fidis of the properly scaled
radius-vector and support functions. Generalized paraboloid growth
processes, coupled with general techniques of stabilization theory, yield
Brownian sheet limits for the defect volume and mean width
functionals. Finally we provide explicit variance asymptotics and
central limit theorems for the $k$-face and intrinsic volume
functionals.
\end{abstract}

%

\end{frontmatter}

\section{Introduction}\label{sec1}\label{INTRO}

Let $K$ be a smooth convex set in $\R^d$ of unit volume. Letting
$\P_\la$ be a Poisson point process in $\R^d$ of intensity $\la$, we
let $K_\la$ be the convex hull of $K \cap\P_\la$. The random
polytope $K_\la$, together with the analogous polytope $K_n$,
obtained by considering $n$ i.i.d. uniformly distributed points in
$K$, are well-studied objects in stochastic geometry.

The study of the asymptotic behavior of the polytopes $K_\la$ and
$K_n$, as $\la\to\infty$ and $n \to\infty$, respectively, has a
long history originating with the work of R\'enyi and Sulanke
\cite{RS}. Letting $\S^{d-1}$ denote the unit sphere, the following
functionals of $K_\la$ have featured prominently:
\begin{itemize}
\item the volume $\Vol(K_\la)$ of $K_{\la},$ abbreviated as
$V(K_{\la})$;
%
\item the number of $k$-dimensional faces of $K_\la$, denoted
$f_k(K_\la), k \in\{0,1,\ldots,\allowbreak d-1\}$; in particular $f_0(K_{\la})$
is the number of vertices of $K_{\la};$
%
\item the mean width $W(K_{\la})$ of $K_\la$;
\item the distance between $\partial K_\la$ and $\partial K$ in the
direction $u \in
\S^{d-1}$, here denoted $r_\la(u), u \neq0$;
\item the distance between the boundary of the Voronoi flower,
defined by $\P_\la$ and $\partial K$, in the direction $u \in
\S^{d-1}$, here denoted $s_\la(u)$;
\item the $k$th intrinsic volumes of $K_\la,$ here denoted $V_k(K_\la
), k \in\{1,\ldots,d-1\}$.
\end{itemize}

The mean values of these functionals on general convex polytopes, as
well as their counterparts for $K_n$, have been widely studied, and
for a complete account we refer to the surveys of Affentranger
\cite{Af}, Buchta \cite{Bu1}, Gruber \cite{Gr}, Reitzner
\cite{ReBook}, Schneider \cite{Sc1,Sc2} and Weil and Wieacker
\cite{WW}, together with Chapter~8.2 in Schneider and Weil
\cite{SW}. There has been recent progress in establishing higher
order and asymptotic normality results for these functionals, for
various choices of~$K$. We signal the important breakthroughs by
Reitzner \cite{Re}, B\'ar\'any and Reitzner \cite{BR2}, B\'ar\'any
et al. \cite{BFV}, Pardon \cite{Pa} and Vu \cite{VV05,VV}.
These results, together with those of Schreiber and Yukich
\cite{SY}, are difficult and technical, with proofs relying upon
tools from convex geometry and probability, including martingales,
concentration inequalities and Stein's method. When $K$ is the unit
radius $d$-dimensional ball $\B^d$
centered at the origin, Schreiber and
Yukich \cite{SY} establish variance asymptotics for $f_0(K_\la)$ as
$\la\to\infty$,
but up to now little is known regarding explicit variance
asymptotics for other functionals of $K_\la$. 

This paper has the following goals. We first study two processes in
\textit{formal space--time} $\R^{d-1}\times\R_+$, one termed the \textit{paraboloid growth process} and denoted by~$\Psi$, and a second termed
the \textit{paraboloid hull process},
denoted by $\Phi$. While the first process was introduced in
\cite{SY}, the second has apparently not been considered before.
When $K= \B^d$, an embedding of convex sets into the space of
continuous functions on $\S^{d-1}$, together with a re-scaling, show
that these processes are naturally suited to the study of $K_\la$.
Their spatial localization can be exploited to describe first and
second order asymptotics of functionals of $K_\la$.
Many of our main results, described as follows, are obtained via
geometric properties of the processes $\Psi$ and $\Phi$. Our goals
are as follows:

$\bullet$ Show that the distance between $K_\la$ and $\partial\B^d$,
upon re-scaling in a local regime, converges in law as $\la\to
\infty$, to a continuous path stochastic process defined in terms of
$\Phi$, adding to Molchanov \cite{Mo}; similarly, we show that the
distance between $\partial\B^d$ and the Voronoi flower defined by
$\P_\la$ converges in law to a continuous path
stochastic process defined in terms 
of $\Psi$. 
In the two-dimensional case the fidis
(finite-dimensional distributions) of these distances, when
re-scaled, are shown to converge to the fidis of $\Psi$ and $\Phi$,
whose description is given explicitly, adding to work of Hsing
\cite{hsing}.

$\bullet$ Show, upon re-scaling in a global regime,
that the suitably integrated local defect
width and defect volume functionals, when considered as processes
indexed by points in $\R^{d-1}$ mapped on $\partial\B^d$ via the
exponential map, satisfy a functional central limit theorem, that
is, converge in the space of continuous functions on $\R^{d-1}$ to a
Brownian sheet on the injectivity region of the exponential map,
whose respective variance coefficients $\sigma^2_W$ and $\sigma^2_V$
are expressed in closed form in terms of $\Psi$ and $\Phi.$
To the best of our knowledge, this connection between the geometry
of random polytopes and Brownian sheets is new.
In particular we show
%
\begin{equation}\label{widthvar} \lim_{\la\to\infty} \la
^{(d +
3)/(d + 1)} \Var[W(K_{\la})] = \sigma^2_W
\end{equation}
and
%
\begin{equation}
\label{volvar} \lim_{\la\to\infty} \la^{(d + 3)/(d + 1)} \Var
[V(K_{\la})] = \sigma^2_V.
\end{equation}
%
This adds to Reitzner's central limit theorem (Theorem 1 of
\cite{Re}) and his variance approximation $\Var[V(K_\la)] \approx
\la^{-(d + 3)/(d + 1)}$ (Theorem 3 and Lemma~1 of \cite{Re}), both
valid when $K$ is an arbitrary smooth convex set. It also adds to
Hsing \cite{hsing}, which is confined to the case $K = \B^2$.

$\bullet$ Establish central limit theorems 
and variance asymptotics for the number of
$k$-dimensional
faces of $K_\la$, showing for all $k \in\{0,1,\ldots,d - 1\}$,
%
\begin{equation}\label{introvar}
\lim_{\la\to\infty} \la^{-(d-1)/(d + 1)} \Var[f_{k}(K_\la)] =
\sigma^2_{f_k},
\end{equation}
where $\sigma^2_{f_k}$ is
described in terms of the processes $\Psi$ and
$\Phi$. 
This improves
upon Reitzner (Lemma 2 of \cite{Re}), whose breakthrough paper
showed $\Var[f_{k}(K_\la)] \approx\la^{(d-1)/(d + 1)}$, and builds
upon \cite{SY}, which establishes (\ref{introvar}) when $k = 0$.

$\bullet$ Establish central limit theorems and variance asymptotics for the
intrinsic volumes
$V_k(K_\la),$ establishing for all $k \in\{1,\ldots,d-1\}$ that
%
\begin{equation}
\label{intrinsicvar} \lim_{\la\to\infty} \la^{(d + 3)/(d + 1)}
\Var[V_k(K_\la)] = \sigma^2_{V_k},
\end{equation}
where again
$\sigma^2_{V_k}$ is described in terms of the processes $\Psi$ and
$\Phi$. This adds to B\'ar\'any et al. (Theorem 1 of \cite{BFV}),
which shows $\Var[V_k(K_n)] \approx n^{-(d+3)/(d + 1)}$.
%

Limits (\ref{widthvar})--(\ref{intrinsicvar}) resolve the issue
of finding variance asymptotics for face functionals and intrinsic
volumes, a long-standing problem put forth this way in the 1993
survey of Weil and Wieacker (page 1431 of \cite{WW}): ``We finally
emphasize that the results described so far give mean values hence
first-order information on random sets and point processes\ldots\ There are also some less geometric methods to obtain higher-order
informations or distributions, but generally the determination of
the variance, for example, is a major open problem.''

These goals are stated in relatively simple terms, and yet they and the
methods 
behind them suggest further objectives involving additional
explanation.
One of our chief objectives is to carefully define the growth
processes $\Psi$ and $\Phi$ and exhibit their geometric
properties making them relevant to $K_\la$, including 
their
localization in space, known as \textit{stabilization}. The latter property is central to establishing
variance asymptotics and the limit theory of functionals of
$K_\la$.
A second objective is to describe two natural scaling regimes, one
suited for locally defined functionals of $K_\la$, and the other
suited for the integrated characteristics of $K_\la$, namely the
width and volume functionals. A third objective is to extend the
afore-mentioned results to ones holding on the level of measures. In
other words, functionals considered here are naturally associated
with random measures, and 
we shall show variance asymptotics for such measures and also
convergence of their fidis to those of a Gaussian process under
suitable global scaling.
We originally intended to restrict attention to convex hulls
generated from Poisson points with intensity density $\la$, but
realized that the methods easily extend to treat intensity densities
decaying as a power of the distance to the boundary of the unit ball
as given by (\ref{DELTAdef}) below, and so we shall include this
more general case without further complication.
These major objectives are discussed
further in the next section.

The extension of the variance
asymptotics \eqref{volvar} and \eqref{introvar} to smooth
compact convex sets with a $C^3$ boundary of positive Gaussian
curvature is nontrivial and is addressed in \cite{CY}. We expect
that much of the limit theory described here can be
``de-Poissonized,'' that is to say, extends to functionals of the
polytope $K_n$. This extension involves challenging technical
questions which we do not address here.


\section{Basic functionals and their scaled versions}\label{sec2}\label{INTRO-a}
Given a locally finite subset $\X$ of $\R^d$, we denote by $\conv(\X
)$ the
\textit{convex hull} generated by $\X.$
For a given compact convex set $K \subset\R^d$ containing the origin,
we let
$h_K\dvtx \S^{d-1} \to\R$ be the support function of $K,$
that is to say, for all $u \in\S^{d-1}$, we let
$h_K(u) := \sup\{\langle x, u\rangle,  x \in K\}$. 
It is easily seen for $\X\subset\R^d$ and $u \in\S^{d-1}$ that
%
\[
h_{\conv(\X)}(u) = \sup\bigl\{h_{\{ x \}}(u),  x \in\X\bigr\}= \sup\{
\langle x, u\rangle,  x \in\X\}.
\]
%
For $u \in\S^{d-1},$ the \textit{radius-vector function}
of $K$ in the direction of $u$ is given by
\[
r_K(u) := \sup\{ \varrho> 0,   \varrho u \in K \}.
\]
%
For $\la> 0$ and $\delta> 0$ we abuse notation and \textit{henceforth
denote by $\P_{\la}:= \P_\la(\delta)$ the Poisson point process
in $\B^d$ of intensity}
%
\begin{equation}\label{DELTAdef}
\la(1-|x|)^{\delta}\,dx, \qquad   x \in\B^d.
\end{equation}
The parameter $\delta$ shall remain fixed throughout, and therefore we
suppress mention of it.
Further, abusing notation we put
\[
K_{\la} := \conv(\P_{\la}).
\]
The principal characteristics of $K_\la$ studied here are the
following functionals, the first two of which represent $K_\la$ in
terms of continuous functions on $\S^{d-1}$:
$\bullet$ The \textit{defect support function.} For all $u \in\S^{d-1}$,
we define
%
\begin{equation}\label{FunctionS}
s_{\la}(u) := s(u,\P_{\la}),
\end{equation}
where for $\X\subseteq\B^d$ we define $ s(u,\X) := 1 - h_{\conv(\X)}(u).
$
In other words, $s_{\la}(u)$ is the defect support function of $K_\la$
in the direction $u.$ 
It is easily verified that $s(u,\X)$ is the distance in
the direction $u$
between the sphere $\S^{d-1}$ and the \textit{Voronoi flower}
%
\begin{equation}\label{flower}
F(\X) := \bigcup_{x \in\X} B_d\biggl(\frac{x}{2}, \frac{|x|}{2}\biggr),
\end{equation}
where for $x \in\R^d$ and $r > 0$ we let $B_d(x,r)$ denote the
$d$-dimensional radius $r$ ball centered at $x$.
%

$\bullet$ The \textit{defect radius-vector function.} For all $u \in\S^{d-1}$,
we define
%
\begin{equation}\label{FunctionR}
r_{\la}(u) := r(u,\P_{\la}),
\end{equation}
where for $\X\subseteq\B^d$ and $u \in\S^{d-1}$ we put $ r(u,\X
):= 1 -
r_{\conv(\X)}(u).$
Thus, $r_{\la}(u)$ is the distance in the direction $u$ between
$\S^{d-1}$ and $K_\la.$
The convex hull $K_\la$ contains the origin, except on a set
of exponentially small probability as $\la\to\infty$, and thus for
asymptotic purposes \textit{we assume without loss of generality
that $K_\la$ always contains the origin},
and therefore the radius vector function $r_{\la}(\cdot)$ is
well defined.
%

$\bullet$ The \textit{numbers of $k$-faces}. Let
$f_{k;\la} := f_k(K_{\la})$, $k \in\{0,1,,\ldots,d-1\}$, be the number
of $k$-dimensional faces of $K_\la.$ In particular,
$f_{0;\la}$ and $f_{1;\la}$ are the number of vertices and edges,
respectively.
The spatial distribution of $k$-faces is captured by
the $k$-face empirical measure (point process) $\mu^{f_k}_{\la}$ on
$\B^d$ given by
%
\begin{equation}\label{FaceEmpM}
\mu^{f_k}_{\la} := \sum_{f \in\Faces_{k}(K_\la)} \delta_{\Top(f)}.
\end{equation}
%
Here \textit{$\Faces_k(K_\la)$ is the collection of all $k$-faces
of $K_\la$} and $\Top(f)$, $f \in\Faces_k(K_\la)$,
is the point of $f$
which is closest to $\S^{d-1},$
with ties ignored as they occur with probability zero (there
are other conceivable choices for $\Top(f)$, but we find
this one to be as good as any).
The total mass $\mu^{f_k}_{\la}(\B^d)$
coincides with~$f_{k;\la}.$

$\bullet$ \textit{Projection avoidance functionals}.
Representing intrinsic volumes of $K_{\la}$ as the total masses of the
corresponding curvature measures,
while suitable in the local scaling regime, turns out to be less useful
in the global
scaling regime, as it leads to an asymptotically vanishing add-one cost
for related stabilizing
functionals, thus precluding normal use of stabilization theory. To
overcome this
problem, we shall use the following consequence of Crofton's general formula,
usually going under the name of Kubota's formula;
see (5.8) and (6.11) in \cite{SW}. We write
%
\begin{equation}\label{KUBOTA}
V_k(K_{\la}) = 
\frac{d! \kappa_d}{k! \kappa_k (d-k)! \kappa_{d-k}}
\int_{G(d,k)} \Vol_k(K_{\la} | L)\,d\nu_k(L),
\end{equation}
%
where $G(d,k)$ is the $k$th Grassmannian of $\R^d$, $\nu_k$ is the
normalized Haar measure
on $G(d,k)$ and $K_{\la} | L$ is the orthogonal projection of $K_{\la
}$ onto the
$k$-dimensional linear space $L \in G(d,k).$ We shall only focus on the
case $k \geq1$
because for $k=0$, we have $V_0(K_{\la}) \equiv1$ for all nonempty,
compact convex
$K_{\la}$;
see page 601 in \cite{SW}. Write
\[
\int_{G(d,k)} \Vol_k(K_{\la} | L)\,d\nu_k(L) = \int_{G(d,k)} \int
_L [1-\vartheta_L(x,\P_\la)]\,dx \, d\nu_k(L),
\]
where $\vartheta_L(x, \X) := {\mathbf 1}(\{x \notin\conv(\X)|L\}).$
Putting $x = ru, u \in\S^{d-1},  r \in[0,1],$ this yields
\begin{eqnarray*}
&&\int_{G(d,k)} \Vol_k(K_{\la} | L)\,d\nu_k(L)\\
&& \qquad  = \int_{G(d,k)} \int
_{\S^{d-1} \cap L} \int_0^1 [1-\vartheta_L(ru,\P_\la)] r^{k-1}\,dr \,d\sigma_{k-1}(u)\, d\nu_k(L)
\\
&& \qquad = \int_{G(d,k)} \int_{\S^{d-1} \cap L} \int_0^1 \frac{1}{r^{d-k}}
[1-\vartheta_L(ru,\P_\la)] r^{d-1}\,dr \,d\sigma_{k-1}(u)\, d\nu_k(L).
\end{eqnarray*}
Noting that $dx = r^{d-1}\,dr\,d\sigma_{d-1}(u)$ and interchanging the
order of integration, we conclude,
in view of the discussion on pages 590--591 of \cite{SW},
that the considered expression equals
%
\[
\frac{k \kappa_k}{d \kappa_d}
\int_{\B^d} \frac{1}{|x|^{d-k}} \int_{G(\operatorname{lin}[x],k)}
[1-\vartheta_L(ru,\P_\la)]\,d\nu^{\operatorname{lin}[x]}_k(L)\,dx,
\]
%
where 
$\operatorname{lin}[x]$ is the $1$-dimensional linear space spanned by $x$,
$G(\operatorname{lin}[x],k)$ is
the set of $k$-dimensional linear subspaces of $\R^d$ containing
$\operatorname{lin}[x]$ and
$\nu^{\operatorname{lin}[x]}_k$ is the corresponding normalized Haar measure;
see \cite{SW}.
Thus, putting
%
\begin{equation}\label{PAVOF}
\vartheta_k^{\X}(x) := \int_{G(\operatorname{lin}[x],k)} \vartheta
_L(x,\X)\,d\nu^{\operatorname{lin}[x]}_k(L),  \qquad   x \in
\B^d,  x \neq0,
\end{equation}
and using (\ref{KUBOTA}), we are led to
%
\begin{eqnarray}\label{KUBOTAREPR}
V_k(\B^d) - V_k(K_{\la}) &=&
\frac{ (^{d-1}_{k-1} )}{\kappa_{d-k}}
\int_{\B^d} \frac{1}{|x|^{d-k}} \vartheta_k^{\P_\la}(x)\,dx\nonumber
\\[-8pt]
\\[-8pt] &=&
\frac{ (^{d-1}_{k-1} )}{\kappa_{d-k}}
\int_{\B^d \setminus K_{\la}} \frac{1}{|x|^{d-k}} \vartheta_k^{\P
_\la}(x)\,dx.
\nonumber
\end{eqnarray}
%
We will refer to $\vartheta_k^{\P_\la}$ as the \textit{projection
avoidance function} for $K_{\la}$.

%

The large $\la$ asymptotics of the above characteristics of $K_\la$ are
studied in two natural scaling regimes, the \textit{local} and the \textit{global} one,
as discussed below.


\textit{Local scaling regime and locally re-scaled functionals.}
The first scaling we consider is referred to as the \textit{local scaling}
in the
sequel.
It stems from the following observation, which, while considered before
in \cite{BR2}, shall be
discussed here in the context of stabilization of growth processes. 
If we consider the local behavior of functionals of $K_\la$ in
the vicinity of two fixed boundary points $u, u' \in\S^{d-1},$
with $\la\to\infty$, then these behaviors become asymptotically
independent. Moreover,
if $u':=u'(\la)$ approaches $u$ slowly enough as $\la\to\infty,$
the asymptotic independence is preserved.
On the other hand, if the distance between $u$ and
$u':= u'(\la)$ decays rapidly enough, then
both behaviors coincide for large $\la$, and the resulting picture is rather
uninteresting. As in \cite{SY}, it is therefore natural to ask for the
frontier of these two asymptotic
regimes and to expect that this corresponds to the natural \textit{characteristic
scale} between the observation directions $u$ and $u'$ where the crucial
features of the local behavior of $K_\la$ are revealed. 

To render the characteristic scale as transparent as possible, we
start with some simple yet important observations, which shall
eventually lead to asymptotic independence of local convex hull
geometries and which shall also suggest the proper scaling limits of
convex hull statistics.
For arbitrary points $x_1,\ldots,x_k \in\B^d$, the support function
of the convex hull $\operatorname{conv}(x_1,\ldots,x_k)$ satisfies for all $u
\in\S^{d-1}$, the relation
\[
h_{\operatorname{conv}(x_1,\ldots,x_k)}(u) =
\max_{1 \leq i \leq k} h_{x_i}(u).
\]
\textit{We make the fundamental observation that the epigraph of $s(u,
\{x_i\}_{i=1}^k):=1 - h_{\operatorname{conv}(x_1,\ldots,x_k)}(u)$ is thus the
union of epigraphs which, locally near the apices, are of parabolic
structure.} Any scaling transformation for $K_\la$ on the
characteristic scale must preserve this structure, as should the
scaling limit for $K_\la$.

To determine the proper local scaling for
our model, we consider the following intuitive argument. To obtain a
nontrivial limit behavior we should re-scale $K_\la$ in a neighborhood
of $\S^{d-1}$, both in the $d-1$ surfacial (tangential) directions
with factor $\lambda^{\beta}$ and radial direction with factor
$\lambda^{\gamma}$
with suitable scaling exponents $\beta$ and $\gamma$ so that:
%
%

$\bullet$ The re-scaling compensates the intensity of $\P_{\la}$ with
growth factor
$\lambda$.
In other words, a subset of $\B^d$ in the vicinity
of $\S^{d-1}$, having a unit volume scaling image, should host
on average $\Theta(1)$ points of 
the point process $\P_{\la}.$ Since the integral of the intensity
density (\ref{DELTAdef})
scales as $\lambda^{\beta(d-1)}$, with respect to the $d-1$ tangential
directions, and since it scales as
$\lambda^{\gamma(1+\delta)}$ with respect to the radial direction,
where we take into account the integration
over the radial coordinate, we are led to
$\lambda^{\beta(d-1) + \gamma(1+\delta)} = \lambda$ and thus
%
\begin{equation}\label{EQBG1}
\beta(d-1) + \gamma(1+\delta) = 1.
\end{equation}

$\bullet$ The local behavior of the convex hull close to the boundary of
$\S^{d-1}$, as described by the locally parabolic structure of
$s_\la$, should preserve parabolic epigraphs, implying for
$u \in\S^{d-1}$ that
$(\la^{\beta} |u|)^2 = \la^{\gamma} |u|^2$, and thus
%
\begin{equation}\label{EQBG2}
\gamma= 2\beta.
\end{equation}
Solving the system (\ref{EQBG1}), (\ref{EQBG2}) we end up with the following
\textit{scaling exponents:}
%
\begin{equation}\label{BETAGAMMA}
\beta= \frac{1}{d+1+2\delta},  \qquad   \gamma= 2 \beta.
\end{equation}

We next describe scaling transformations for $K_\la.$
Fix $u_0 \in\S^{d-1}$,
and let $T_{u_0}:= T_{u_0}\S^{d-1}$ denote the tangent space to
$\S^{d-1}$ at $u_0$.
The \textit{exponential
map} $\exp_{u_0} \dvtx  T_{u_0} \S^{d-1}\to\S^{d-1}$ maps a vector $v$
of the tangent space
$T_{u_0}$ to the point $u \in\S^{d-1}$, such that $u$ lies at the end
of the geodesic of length $|v|$ starting at $u_0$ and having direction $v.$
Note that $\S^{d-1}$ is geodesically complete in that the exponential
map $\exp_{u_0}$ is well defined on the whole tangent space
$\R^{d-1} \simeq T_{u_0}\S^{d-1},$ although it is injective only
on $\{ v \in T_{u_0} \S^{d-1},  |v| < \pi\}.$ Instead of $\exp
_{u_0}$, we shall write
$\exp_{d-1}$ or simply $\exp$,
and we make the default choice
$u_0 := (0,0,\ldots,1).$ 
We use the isomorphism $T_{u_0} \S^{d-1} \simeq\R^{d-1}$
without further mention, and \textit{we shall denote the closure of
the injectivity region $\{ v \in T_{u_0} \S^{d-1},  |v| < \pi\}$ of
the exponential map simply by $\B_{d-1}(\pi).$} Thus we have
$\exp(\B_{d-1}(\pi)) = \S^{d-1}.$

Further, consider the following scaling transform $T^{\la}$ mapping
$\B^d$ into
$\R^{d-1} \times\R_+$
%
\begin{equation}\label{SCTRANSF}
T^{\la}(x) :=  \biggl(\la^{\beta} \exp^{-1}_{d-1}\biggl( \frac{x}{|x|}\biggr),
\la^{\gamma}
(1-|x|) \biggr),  \qquad  x \in\B^d\setminus\{\0\}.
\end{equation}
%
Here $\exp^{-1}(\cdot)$ is the inverse exponential map, which is well defined
on $\S^{d-1} \setminus\{ - u_0 \}$ and which takes values in the
injectivity region
$\B_{d-1}(\pi)$. For formal completeness, on the ``missing''
point $-u_0$, we let $\exp^{-1}$ admit an arbitrary value, say
$(0,0,\ldots,\pi),$ and likewise we put
$T^{\la}({\bf0}) := (\0,\la^{\gamma}),$
where $\0$ denotes either the origin of $\R^{d-1}$ or $\R^d$,
according to the context.
It is easily seen that $T^{\la}$ is a.e. (with respect to Lebesgue
measure on $\B^d$) a bijection from
$\B^d$ onto the
$d$-dimensional solid cylinders
%
\begin{equation}\label{RLA}
{\cal R}_{\la} := \la^{\beta} \B_{d-1}(\pi) \times[0,\la^{\gamma}).
\end{equation}
%
Throughout points in $\B^d$ are
written as $x: = (r, u)$, and we represent generic points in
$\R^{d-1} \times\R_+$ by $(v, h)$, whereas we write $(v',h')$ to
represent
points in the scaled region ${\cal R}_{\la}$. 
We assert that the transformation $T^\la$, defined at
(\ref{SCTRANSF}), maps the Poisson point process $\P_\la$ to $\P
^{(\la)}$,
where $\P^{(\la)}$ is the dilated Poisson point process in the region
${\cal R}_{\la}$ having intensity
%
\begin{equation}\label{DENS1}
  (v',h') \mapsto\frac{\sin^{d-2}(\la^{-\beta}|v'|)}{|\lambda
^{-\beta}v'|^{d-2}} (1-\la^{-\gamma}h')^{d-1}
{h'}^{\delta}\,dv' \,dh'
\end{equation}
at $(v',h') \in{\cal R}_{\la}.$
Indeed, this intensity measure is the image by the transformation
$T^{\la}$ of the measure on $\B^d$ given by
%
\begin{equation}\label{measure0}
\la(1-|x|)^{\delta}\,dx=\la(1-r)^{\delta}r^{d-1}\,dr\,d\sigma_{d-1}(u)
\end{equation}
introduced in (\ref{DELTAdef}), where we put $x = (r, u)$. To
obtain (\ref{DENS1}), we first make a change of variables,
\[
h':=\la^{\gamma}(1-r)  \quad \mbox{and} \quad  v':= \la^{\beta}\exp
_{d-1}^{-1} (u ) := \la^{\beta} v.
\]
Next, notice that the exponential map $\exp_{d-1} \dvtx  T_{u_0}
\S^{d-1}\to\S^{d-1}$ has the following expression:
%
\begin{equation}
\label{expomap}
\exp_{d-1}(v')= \cos(|v'|)(0,\ldots,0,1) +
\sin(|v'|)\biggl(\frac{v'}{|v'|},0\biggr),
\end{equation}
with $v'\in\R^{d-1}\setminus\{\0\}$.
Therefore, since $v:= \exp_{d-1}^{-1}(u)$, we have
\[
d\sigma_{d-1}(u)= \sin^{d-2}(|v|) d(|v|)\,d\sigma_{d-2}\biggl(
\frac{v} {|v|}\biggr) = \frac{\sin^{d-2}(|v|)\,dv} {|v|^{d-2}}.
\]
Since $v' = \la^{\beta} v$, this gives
%
\begin{equation}
d\sigma_{d-1}(u)= \frac{\sin^{d-2}(\la^{-\beta}|v'|)}{|\la
^{-\beta}v'|^{d-2}}\la^{-\beta(d-1)}\,dv'.
\label{CV2}
\end{equation}
We also have that
%
\begin{equation}\label{CV1}
(1-r)^{\delta}r^{d-1}\,dr=\la^{-\gamma\delta}{h'}^{\delta}(1-\la
^{-\gamma}{h'})^{d-1}\la^{-\gamma}d{h'}.
\end{equation}
Inserting (\ref{CV2}) and (\ref{CV1}) in (\ref{measure0}) and
using (\ref{EQBG1}) to obtain $\la\la^{-\beta(d-1)}
\la^{-\gamma(1+\delta)} = 1$, we obtain (\ref{DENS1}).


In Section~\ref{LocScal}, following \cite{SY}, we shall embed
$T^\la(K_\la)$ into a space of paraboloid growth processes on $\cal
R_\la$. One such process, denoted by $\Psi^{(\la)}$ and defined at
(\ref{gengrowth}), is a \textit{generalized growth process with
overlap} whereas the second, a dual process denoted by
$\Phi^{(\la)}$ and defined at (\ref{phullprocess}) 
is termed the \textit{paraboloid hull process}.
Infinite volume counterparts to $\Psi^{(\la)}$ and $\Phi^{(\la)}$,
described fully in Section~\ref{g+hullprocess} and denoted by $\Psi$
and $\Phi$, respectively, play a natural role in describing the
asymptotic behavior of our basic functionals of interest, re-scaled
as follows: 

$\bullet$ The \textit{re-scaled versions of the defect support function
(\ref{FunctionS}) and the radius support
function  (\ref{FunctionR})}, defined, respectively, by
%
\begin{eqnarray}\label{SRESC}
\hat{s}_{\la}(v) &:=& \la^{\gamma} s_{\la}(\exp_{d-1}(\la^{-\beta
}v)),  \qquad  v \in\R^{d-1},
\\\label{RRESC}
\hat{r}_{\la}(v) &:=& \la^{\gamma} r_{\la}(\exp_{d-1}(\la^{-\beta
}v)),  \qquad  v \in\R^{d-1}.
\end{eqnarray}



$\bullet$ The \textit{re-scaled version of the projection avoidance function
(\ref{PAVOF})} defined by
%
\begin{equation}\label{PAVOFRESC}
\hat\vartheta_k^{\P_\la}(x) := \vartheta_k^{\P_\la}([T^{\la
}]^{-1}(x)),  \qquad  x \in{\cal R}_{\la}.
\end{equation}
%

\textit{Global scaling regime and globally re-scaled functionals.}
The asymptotic independence of local convex hull geometries at
distinct points of $\S^{d-1},$ as discussed above, suggests that
the global behavior of both $s_{\la}$ and $r_{\la}$ is, in large $\la
$ asymptotics,
that of the \textit{white noise}.
Therefore it is natural to consider the corresponding
integral characteristics of $K_\la$ and to ask whether, under proper scaling,
they converge in law to a Brownian sheet. Define the processes
%
\begin{equation}\label{SGR}
W_{\la}(v) := \int_{\exp([\0,v])} s_{\la}(u)\,d\sigma_{d-1}(u),  \qquad  v
\in\R^{d-1},
\end{equation}
and
%
\begin{equation}\label{RGR}
V_{\la}(v) := \int_{\exp([\0,v])} r_{\la}(u)\,d\sigma_{d-1}(u), \qquad   v
\in\R^{d-1},
\end{equation}
%
where the ``segment'' $[\0,v]$ for $v \in\R^{d-1}$ is the rectangular solid
in $\R^{d-1}$ with vertices $\0$ and $v,$ that is to say, $[\0,v] :=
\prod_{i=1}^{d-1} [\min(0,v^{(i)}),\max(0,v^{(i)})]$, with $v^{(i)}$ standing
for the $i$th coordinate of $v$. 
We shall also consider the cumulative values
%
\begin{eqnarray}\label{GRINF}
W_{\la} &:=& W_{\la}(\infty) := \int_{\S^{d-1}} s_{\la}(u)\,d\sigma
_{d-1}(u);
\nonumber
  \\[-8pt]
  \\[-8pt]  V_{\la} &:=& V_{\la}(\infty) := \int_{\S^{d-1}} r_{\la}(u)\,d\sigma_{d-1}(u).
\nonumber
\end{eqnarray}
%

Notice that the radius-vector
function of the Voronoi flower $F(\P_\la)$ coincides with the
support function\vadjust{\goodbreak} of $K_\la$. In particular, the volume outside
$F(\P_\la)$ is equal to
%
\begin{equation}\label{ajout}
  \int_{\S^{d-1}}\biggl [\int_{1-s_{\la}(u)}^{1}\rho^{d-1}\,d\rho
 \biggr]\,d\sigma_{d-1}(u)=\int_{\S^{d-1}}\frac{1-(1-s_{\la
}(u))^d}{d}\,d\sigma_{d-1}(u).\hspace*{-35pt}
\end{equation}
Since
$s_{\la}$ goes to $0$ uniformly, the volume outside
$F(\P_\la)$ is asymptotically equivalent to the integral of the
defect support function, which in turn is proportional to the
defect mean width by Cauchy's formula.
Moreover, in two dimensions the mean width is the
ratio of the perimeter to $\pi$ (see page 210 of \cite{Schn}), and so
\textit{$W_{\la}(\infty)/\pi$ coincides with $2$ minus the mean width
of $K_\la$, and consequently $W_{\la}(\infty)$ itself equals $2\pi$
minus the perimeter of $K_\la$ for \mbox{$d=2.$}} On the other hand,
$V_{\la}(\infty)$ is asymptotic
to the volume of $\B^d \setminus K_\la,$ \textit{whence the notation
$W$ for (asymptotic)  width and $V$ for  (asymptotic)  volume}.

To get the desired
convergence to a Brownian sheet, we put 
%
\begin{equation}\label{ZETA}
\zeta:= \beta(d-1) + 2 \gamma= \frac{d+3}{d+1+2\delta};
\end{equation}
we show in Section~\ref{sec8} that it is natural to re-scale the
processes $(W_{\la}(v) - \E W_{\la}(v))$ and $(V_{\la}(v) - \E
V_{\la}(v))$ by $\la^{\zeta/2}$ and that the resulting re-scaled
processes
%
\begin{eqnarray}\label{SGRESC}
\hat{W}_{\la}(v) &:=& \la^{\zeta/2} \bigl(W_{\la}(v) - \E W_{\la}(v)\bigr)
  \quad \mbox{and}\nonumber
  \\[-8pt]
  \\[-8pt]
\hat{V}_{\la}(v) &:= &\la^{\zeta/2}\bigl (V_{\la}(v) - \E V_{\la}(v)\bigr), \qquad  v \in\R^{d-1},
\nonumber
\end{eqnarray}
converge in law to a Brownian sheet with an explicit variance
coefficient.\vspace*{6pt}

\textit{Putting the picture together.} The remainder of this
paper is organized as follows.

\textit{Section~\ref{g+hullprocess}}. Though the formulation of our
results might suggest otherwise, there are crucial connections
between the local and global scaling regimes. These regimes are
linked by stabilization 
and the objective method, which together show that the behavior of
locally defined processes on the finite volume rectangular solids
${\cal R}_{\la}$, defined at (\ref{RLA}), can be well approximated by
the local behavior of a related ``candidate object,'' either \textit{a
generalized growth process} $\Psi$ or a \textit{dual paraboloid hull
process} $\Phi$, on an \textit{infinite volume half-space}. While
generalized growth processes were developed in \cite{SY} in a larger
context, our limit theory depends heavily on a new object, the dual
paraboloid hull process. The purpose of Section~\ref{g+hullprocess}
is to carefully define these processes and to establish properties
relevant to their asymptotic analysis.

\textit{Section~\ref{LocScal}}. We show that as $\la\to\infty,$
both $\hat{s}_{\la}$ and $\hat{r}_{\la}$, defined, respectively, at
(\ref{SRESC}) and
(\ref{RRESC}),
converge in law to
continuous path stochastic processes explicitly constructed in
terms of the paraboloid generalized growth process $\Psi$ and the paraboloid
hull process $\Phi$, respectively.
This adds to Molchanov~\cite{Mo}, who considers the ``epiconvergence'' in the space $\S^{d-1}
\times\R$
of the random process, arising as the binomial counterpart of $\la
r_\la$. Molchanov's results \cite{Mo} are not framed in terms of the rescaled
function $\hat{r}_{\la}$, and thus they do not
involve the paraboloid growth processes described in this paper.

\textit{Section~\ref{PC1}.} When $d = 2$, after re-scaling in space by
a factor of $\la^{1/3}$ and in time (height coordinate)
by $\la^{2/3}$, we use nonasymptotic direct considerations to provide
explicit asymptotic expressions for the
fidis of $\hat{s}_\la$ and $\hat{r}_\la$
as $\la\to\infty$. These distributions coincide with the fidis 
of the parabolic growth process $\Psi$ and the
parabolic
hull process $\Phi$, respectively.

\textit{Section~\ref{StabRepr}}. Both the paraboloid growth
process $\Psi$ and its dual paraboloid hull process $\Phi$ are shown
to enjoy a localization property, which expresses, in geometric terms,
a type of spatial mixing. This provides a direct route toward
establishing first and second order asymptotics for the convex hull
functionals of interest.


\textit{Section~\ref{GAUSS}}.
This section establishes closed form variance asymptotics for the
total number of $k$-faces as well as the intrinsic volumes for the
random polytope $K_\la$. We also establish variance asymptotics and
a central limit theorem for the properly scaled integrals of
continuous test functions against the empirical measures associated
with the functionals under proper scaling.

\textit{Section~\ref{InvPrinc}}. Using the stabilization properties
established in Section~\ref{StabRepr},
we establish a functional central limit theorem for $\hat{W}_{\la}$
and $\hat{V}_{\la},$ showing that these processes converge, as $\la
\to\infty$ in the space of continuous functions on $\R^{d-1}$, to
Brownian sheets with variance coefficients given in terms of the
processes $\Psi$ and $\Phi$, respectively.

\section{Paraboloid growth and hull processes}\label{sec3}\label{g+hullprocess}

In this section 
we introduce the
paraboloid \textit{growth} and \textit{hull} processes in the upper half-space
$\R^{d-1} \times\R_+$ often interpreted as formal \textit{space--time}
below,
with $\R^{d-1}$ standing for the spatial dimension and $\R_+$
standing for
the time dimension. Although this interpretation is purely formal in
the convex
hull set-up, it provides a link to a well-established theory of growth
processes studied by means of stabilization theory; see below for further
details. These processes turn out to be infinite volume counterparts
to finite volume paraboloid growth processes, which are defined
in the next section, and which are used to describe the behavior of
our basic re-scaled functionals and measures.

\textit{Poisson point process on half-spaces.} Fix $\delta> 0$,
and let
$\P(\delta)$ be a Poisson point
process in $\R^{d-1} \times\R_+$ with intensity
density
%
\begin{equation}\label{basicPPP} h^{\delta}\,dh \,dv   \qquad \mbox{at } (v,h) \in\R
^{d-1} \times
\R_+.
\end{equation}
In the sequel we shall show that the scaled Poisson point process
$\P^{(\la)}:= T^\la(\P_\la)$ with intensity defined at (\ref{DENS1})
converges to $\P(\delta)$ on compacts, but for now we use the process
$\P$
to define growth processes on half-spaces. As with $\P_{\la}$ and
$\P^{(\la)},$ we suppress $\delta$ and simply write $\P$ for
$\P(\delta)$.

\textit{Paraboloid growth processes on half-spaces.}
We introduce the \textit{paraboloid generalized growth process with
overlap} (paraboloid growth
process for short), specializing to our present set-up the
corresponding general
concept defined in Section~1.1 of \cite{SY} and designed to
constitute the asymptotic
counterpart of the Voronoi flower $F(K_{\la}).$
Let $\Pi^{\uparrow}$ be the epigraph of the standard paraboloid $v
\mapsto|v|^2/2,$ that is,
\[
\Pi^{\uparrow} := \biggl\{ (v,h) \in\R^{d-1} \times\R_+,  h \geq\frac
{|v|^2} {2} \biggr\}.
\]

We introduce one of the fundamental objects of this
paper.

\begin{defn} \label{fund-def} Given a locally finite point set $\X$
in $\R^{d-1} \times\R_+$, the paraboloid growth model $\Psi(\X)$
is defined as the Boolean model
with parabo\-loid grain $\Pi^{\uparrow}$ and with germ collection
$\X$, namely
%
\begin{equation}\label{PSIDEF}
\Psi(\X) := \X\oplus\Pi^{\uparrow} = \bigcup_{x \in\X} x
\oplus\Pi^{\uparrow},
\end{equation}
where $\oplus$ stands for Minkowski addition. In particular, we define
the parabo\-loid growth process
$\Psi:= \Psi(\P)$, where $\P$ is the Poisson point
process defined at (\ref{basicPPP}).
\end{defn}

\begin{figure}[b]

\includegraphics{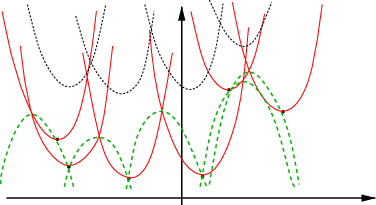}

\caption{Example of paraboloid and growth processes for $d=2$.}\label{parabolapicture}
\end{figure}

The model
$\Psi(\X)$ arises as the union of upwards paraboloids with apices at
the points of $\X$
(see Figure~\ref{parabolapicture}), in close analogy to the Voronoi 
flower $F(\X)$, where to each $x \in\X$ we attach a ball $B_d(x/2,|x|/2)$
(which asymptotically scales to an upwards paraboloid as we shall see
in the sequel)
and take the union thereof.

The name \textit{generalized growth process with overlap} comes from the original
interpretation of this construction \cite{SY}, where $\R^{d-1} \times
\R_+$
stands for \textit{space--time} with $\R^{d-1}$ corresponding to the \textit{spatial}
coordinates and the semi-axis $\R_+$ corresponding to the \textit{time}
(\textit{or height})
coordinate, and where the grain $\Pi^{\uparrow},$ possibly admitting
more general
shapes as well, arises as the graph of the growth of a germ born at the apex
of $\Pi^{\uparrow}$ and growing thereupon in time with properly
varying speed.
We say that the process \textit{admits overlaps} because the growth does
not stop
when two grains overlap, unlike in traditional growth schemes.
We shall often use this space--time interpretation
and refer to the respective coordinate axes as to the spatial and time (height)
axis.

The boundary $\partial\Psi$ of the random closed set $\Psi:=
\Psi(\P)$ constitutes a graph of a continuous function from
$\R^{d-1}$ (space) to $\R_+$ (time), also denoted by $\partial
\Psi$ in the sequel. In what follows we
interpret $\hat{s}_\la$, defined at (\ref{SRESC}),
as the boundary of a growth
process $\Psi^{(\la)}$, defined at (\ref{gengrowth}) below, on the
finite region ${\cal R}_{\la}$ at (\ref{RLA}); we shall see in
Section~\ref{sec4} that $\partial\Psi$ is the scaling limit for the
boundary of $\Psi^{(\la)}$.

A germ point $x \in\P$ is called \textit{extreme} in the
paraboloid growth process $\Psi$ if and only if its associated epigraph
$x \oplus\Pi^{\uparrow}$ is \textit{not} contained in the union of
the paraboloid epigraphs generated by other germ points in $\P,$
that is to say,
%
\begin{equation}\label{PSIEXTR}
x \oplus\Pi^{\uparrow} \not\subseteq\bigcup_{y \in\P, y \neq x} (y
\oplus\Pi^{\uparrow}).
\end{equation}
For $x$ to be extreme, it is sufficient, but not necessary, that $x$
fails to be contained in
paraboloid epigraphs of other germs. Write $\ext(\Psi)$ for the set
of all
extreme points.

%

%
\textit{Paraboloid hull process on half-spaces.}
The \textit{paraboloid hull process} $\Phi$ can be regarded
as the dual to the paraboloid growth process. At the same time, the
paraboloid hull process
is designed to exhibit geometric properties analogous to those of
convex polytopes with
\textit{paraboloids} playing the role of
\textit{hyperplanes}, with the
\textit{spatial coordinates} playing the role of \textit{spherical coordinates}
and with the \textit{height/time coordinate} playing the role of the \textit{radial
coordinate}. The motivation
of this construction is to mimic the convex geometry on
second order paraboloid structures in order to describe the local second
order geometry of convex polytopes, which dominates their limit
behavior in
smooth convex bodies. As we shall see, this intuition is
indeed correct and results in a detailed description of the limit behavior
of $K_\la.$

To proceed with our definitions, we let $\Pi^{\downarrow}$ be the downwards
space--time paraboloid hypograph
%
\begin{equation}\label{PDOWN}
\Pi^{\downarrow} := \biggl\{ (v,h) \in\R^{d-1} \times\R,  h \leq-
\frac{|v|^2} {2} \biggr\}.
\end{equation}
The idea behind our interpretation of the paraboloid process
is that the shifts of $\Pi^{\downarrow}$ correspond to half-spaces
not containing
${\bf0}$ in the Euclidean space $\R^{d}.$ 
We shall argue the \textit{paraboloid convex sets} have properties
strongly analogous to those
related to the usual concept of convexity. The corresponding proofs are not
difficult and will be presented in enough detail to make our presentation
self-contained, but it should be emphasized that alternatively
the entire argument of this paragraph could be re-written in terms of the
following \textit{trick}. Considering the transform $(v,h) \mapsto(v,h+|v|^2/2)$,
we see that it maps translates of $\Pi^{\downarrow}$ to half-spaces and
thus whenever we make a statement
below in terms of paraboloids and claim it is analogous to a standard statement
of convex geometry, we can alternatively apply the above auxiliary transform,
use the classical result and then transform back to our set-up. We do
not choose
this option here, finding it more aesthetic to work directly in the paraboloid
set-up, but we indicate at this point the availability of this alternative.

The next definitions are central
to the description of the paraboloid hull process. Recall that the
affine hull $\Span[v_1,\ldots,v_k]$ is the set of all affine
combinations
$\alpha_1 v_1 +\cdots + \alpha_k v_k,  \alpha_1 +\cdots + \alpha
_k = 1,
\alpha_i \in\R.$

\begin{defn} For any collection $x_1 := (v_1,h_1),\ldots
,x_k:=(v_k,h_k),  k \leq d,$
of points in $\R^{d-1} \times\R_+$ with affinely independent spatial
coordinates $v_i,$ we define
$\Pi^{\downarrow}[x_1,\ldots,x_k]$ to be the hypograph
in $\Span[v_1,\ldots,v_k] \times\R$ of the unique space--time paraboloid
in the affine space $\Span[v_1,\ldots,v_k] \times\R$ with quadratic
coefficient $-1/2$ and passing through $x_1,\ldots,x_k.$
\end{defn}

In other words $\Pi^{\downarrow}[x_1,\ldots,x_k]$ is the
intersection of
$\Span[v_1,\ldots,v_k] \times\R$ and a translate of $\Pi
^{\downarrow}$
having all $x_1,\ldots,x_k$ on its boundary; while such translates
are nonunique for $k < d,$ their intersections with
$\Span[v_1,\ldots,v_k]$ all coincide.

\begin{defn} For $x_1:= (v_1,h_1) \neq x_2:= (v_2,h_2) \in\R^{d-1}
\times\R_+,$
the \textit{parabolic segment} $\Pi^{[\cdot]}[x_1,x_2]$ 
is the unique parabolic segment with quadratic coefficient
$-1/2$ joining $x_1$ to $x_2$ in $\Span[v_1,v_2] \times\R.$ More generally,
for any collection $x_1:= (v_1,h_1),\ldots,x_k:= (v_k,h_k),  k \leq
d,$ of
points in $\R^{d-1} \times\R_+$ with affinely independent spatial
coordinates, we define the \textit{paraboloid face}
$\Pi^{[\cdot]}[x_1,\ldots,x_k]$ by
%
\begin{equation}\label{pface}\Pi^{[\cdot]}[x_1,\ldots,x_k]:= \partial
\Pi^{\downarrow}[x_1,\ldots,x_k] \cap
[\conv(v_1,\ldots,v_k) \times\R].
\end{equation}
\end{defn}

Clearly, $\Pi^{[\cdot]}[x_1,\ldots,x_k]$
is the smallest set containing $x_1,\ldots,x_k$ and with the \textit{paraboloid convexity}
property: For any two $y_1,y_2$ it contains, it also contains $\Pi
^{[\cdot]}[y_1,y_2].$
In these terms, $\Pi^{[\cdot]}[x_1,\ldots,x_k]$ is the \textit{paraboloid
convex hull}
$\phull(\{x_1,\ldots,x_k\}).$ In particular, we readily derive the property
%
\begin{eqnarray}\label{INTERSECT}
&&\Pi^{[\cdot]}[x_1,\ldots,x_i,\ldots,x_k] \cap\Pi^{[\cdot]}[x_i,\ldots
,x_k,\ldots,x_m] \nonumber
\\[-8pt]
\\[-8pt]
&& \qquad =
\Pi^{[\cdot]}[x_i,\ldots,x_k], \qquad 1 < i < k.\nonumber
\end{eqnarray}

Next, we say that $A \subseteq\R^{d-1} \times\R_+$ is \textit{upwards
paraboloid convex}
(up-convex for short) if and only if:
\begin{itemize}
\item for each two $x_1,x_2 \in A$ we have $\Pi^{[\cdot]}[x_1,x_2]
\subseteq A;$
\item and for each $x = (v,h) \in A$ we have $x^{\uparrow} := \{
(v,h'),  h' \geq h \} \subseteq A.$\vadjust{\goodbreak}
\end{itemize}
Whereas the first condition in the definition above is quite intuitive,
the second
will be seen to correspond to our requirement that ${\bf0} \in K_\la$
as ${\bf0}$ gets transformed to \textit{upper infinity} in the limit of
our re-scalings. Indeed, though $T^{\la}$ is not defined at $x=\mathbf{0}$, the last coordinate of $T^{\la}(x)$ goes to $\la^{\gamma}$ when
$x\to{\bf0}$, and $\la^{\gamma}$ goes to $\infty$ when $\la\to
\infty$.

With the notation introduced above, we now define the second
fundamental object of this paper.


\begin{defn}  Given $A \subseteq\R^{d-1} \times\R_+$, by the
paraboloid hull (up-hull for short)
of $A$, we mean the smallest up-convex set containing $A.$ Given
a locally finite point set $\X\in\R^{d-1} \times\R_+$, we
define the paraboloid hull $\Phi(\X)$ to be the up-hull of $\X$,
that is,
\[
\Phi(\X):= \uphull(\X).
\]
In particular, we define the
paraboloid hull process $\Phi$ in $\R^{d-1} \times\R_+$ as the
up-hull of $\P,$ that is to say,
%
\begin{equation}\label{PHIDEF}
\Phi:= \Phi(\P):= \uphull(\P).
\end{equation}
\end{defn}

For $A \subseteq\R^{d-1} \times\R_+$ we put $ A^{\uparrow} := \{
(v,h'),  (v,h) \in A \mbox{ for some } h \leq h' \}$
and observe that if $x'_1 \in x_1^{\uparrow},  x'_2 \in x_2^{\uparrow
},$ then
%
\begin{equation}\label{ODCZAW}
\Pi^{[\cdot]}[x'_1,x'_2] \subset\bigl[\Pi^{[\cdot]}[x_1,x_2]\bigr]^{\uparrow}
\end{equation}
and, more generally, by
definition of $\Pi^{[\cdot]}[x_1,\ldots,x_k]$ and by induction in $k,$
$\Pi^{[\cdot]}[x'_1,\ldots,x'_k] \subset[\Pi^{[\cdot]}[x_1,\ldots
,x_k]]^{\uparrow}.$
Consequently, we conclude that
%
\begin{equation}\label{PHIOKR}
\Phi= [\phull(\P)]^{\uparrow},
\end{equation}
which, in terms of our analogy between convex polytopes and paraboloid
hulls processes,
reduces to the trivial statement that a convex polytope containing
${\bf0}$ arises
as the union of radial segments joining ${\bf0}$ to convex
combinations of its
vertices. This statement is somewhat more interesting in the present
set-up where
${\bf0}$ \textit{disappears} at infinity, and we formulate it here for
further use.

\begin{lemm} \label{JYlem} With probability $1$ we have
%
\begin{equation}\label{PHIOPIS1}
\Phi= \bigcup_{\{x_1,\ldots,x_d\} \subset\P} \bigl[\Pi^{[\cdot]}[x_1,\ldots
,x_d]\bigr]^{\uparrow}.
\end{equation}
\end{lemm}

This statement corresponds to the property of $d$-dimensional polytopes
containing
${\bf0},$ stating that the convex hull of a collection of points
containing ${\bf0}$
is the union of all $d$-dimensional simplices with vertex sets running
over all
cardinality $(d+1)$ sub-collections of the generating collection which
contain ${\bf0}.$
Subsets $\{ x_1, \ldots, x_d \} \subset\P$ have their spatial coordinates
affinely independent with probability $1$ and thus the right-hand side
in (\ref{PHIOPIS1}) is a.s.
well defined; in the sequel we shall say that points of $\P$ are a.s.
\textit{in general
position}.

\begin{pf} Observe that, in view of (\ref{PHIOKR}) and the fact that
\[
\bigcup_{\{x_1,\ldots,x_d\} \subset\P} \Pi^{[\cdot]}[x_1,\allowbreak \ldots,x_d]
\subset\phull(\P),
\]
(\ref{PHIOPIS1}) will follow as soon as we show that
%
\begin{equation}\label{SPROW1}
\phull(\P) \subset\bigcup_{\{x_1,\ldots,x_d\} \subset\P} \bigl[\Pi
^{[\cdot]}[x_1,\ldots,x_d]\bigr]^{\uparrow}.
\end{equation}
To establish (\ref{SPROW1}) it suffices to show that
adding an extra point $x_{d+1}$ in general position to a set
$\bar{x} = \{x_1,\ldots,x_d\}$ results in having
%
\begin{equation}\label{SPROW2}
\phull (\bar{x}^+ := \bar{x} \cup\{ x_{d+1}\}  ) \subset
\bigcup_{i=1}^{d+1} \bigl[\Pi^{[\cdot]}[\bar{x}^+ \setminus\{ x_i \}
]\bigr]^{\uparrow},
\end{equation}
and inductive use of this fact readily yields the required relation
(\ref{SPROW1}).
To verify (\ref{SPROW2}) choose $y := (v,h) \in\phull(\bar{x}^+).$
Then there
exists $y' = (v',h') \in\Pi^{[\cdot]}[x_1,\ldots,x_d]$ such that $y \in
\Pi^{[\cdot]}[y',x_{d+1}].$
Consider the section of $\Pi^{[\cdot]}[x_1,\ldots,\allowbreak x_d]$ by the plane
$\Span[v',v_{d+1}] \times\R$
and $y''$ be its point with the lowest height coordinate. Clearly then
there exists
$x_i,  i \in\{ 1, \ldots, d \}$ such that $y'' \in\Pi^{[\cdot]}[\bar
{x} \setminus\{ x_i \}].$
On the other hand, by the choice of $y''$ and by (\ref{ODCZAW}), $y
\in\Pi^{[\cdot]}[y',x_{d+1}]
\subset[\Pi^{[\cdot]}[y',y'']]^{\uparrow} \cup[\Pi
^{[\cdot]}[y'',x_{d+1}]]^{\uparrow}.$
Consequently, $y \in[\Pi^{[\cdot]}[\bar{x}]]^{\uparrow}\cup[\Pi
^{[\cdot]}[\bar{x}^+ \setminus\{ x_i \}]]^{\uparrow}$, which completes
the proof of (\ref{SPROW2})
and thus also of (\ref{SPROW1}) and (\ref{PHIOPIS1}). This
completes the proof of Lemma~\ref{JYlem}.
\end{pf}

To formulate our next statement, we say that a collection $\{x_1,\ldots
,x_d\}$ is
\textit{extreme} in $\P$ if and only if $\Pi^{[\cdot]}[x_1,\ldots,x_d]
\subset\partial\Phi.$ Note
that, by (\ref{ODCZAW}) and Lemma~\ref{JYlem}, this is equivalent to having
%
\begin{equation}\label{PARAEXT}
\Phi\cap\Pi^{\downarrow}[x_1,\ldots,x_d] = \Pi^{[\cdot]}[x_1,\ldots,x_d].
\end{equation}
Each such $\Pi^{[\cdot]}[x_1,\ldots,x_d]$ is referred to as a \textit{paraboloid sub-face}.
Further, say that two extreme collections $\{x_1,\ldots,x_d\}$ and
$\{x'_1,\ldots,x'_d\}$ in $\P$ are co-paraboloid if and only if
$\Pi^{\downarrow}[x_1,\ldots,x_d] = \Pi^{\downarrow}[x'_1,\ldots,x'_d].$
By a \textit{$(d-1)$-dimensional paraboloid face} of $\Phi$, we shall understand
the union of each maximal collection of co-paraboloid sub-faces.
Clearly, these
correspond to $(d-1)$-dimensional faces of convex polytopes. It is not difficult
to check that $(d-1)$-dimensional paraboloid faces of $\Phi$ are
p-convex, and
their union is $\partial\Phi.$ In fact, since $\P$ is a Poisson
process, with
probability one all $(d-1)$-dimensional faces of $\Phi$ consist of precisely
one sub-face; in particular all $(d-1)$-dimensional faces of $\Phi$
are bounded.
By (\ref{PARAEXT}) we have for each $(d-1)$-dimensional face $f$,
%
\begin{equation}\label{PARAEXT2}
\Phi\cap\Pi^{\downarrow}[f] = f,
\end{equation}
which corresponds to the standard fact of the theory of convex
polytopes, stating
that the intersection of a $d$-dimensional polytope containing $\mathbf{0}$ with
a half-space determined by a $(d-1)$-dimensional face, and looking away from
${\bf0}$, is precisely the face itself. Further, pairs of adjacent
$(d-1)$-dimensional paraboloid faces intersect yielding $(d-2)$-dimensional
paraboloid manifolds, called $(d-2)$-dimensional paraboloid faces.
More
generally, $(d-k)$-dimensional paraboloid faces arise as $(d-k)$-dimensional
paraboloid manifolds obtained by intersecting suitable $k$-tuples of adjacent
$(d-1)$-dimensional faces. Finally, we end up with zero dimensional
faces, which
are the \textit{vertices} of $\Phi$, and which are easily seen to belong
to $\P.$
The set of vertices of $\Phi$ is denoted by $\Vertices(\Phi).$
In other words, we obtain a full analogy with the geometry of faces of
$d$-dimensional polytopes. Clearly, 
$\partial\Phi$ is the graph of a continuous piecewise paraboloid
function from $\R^{d-1}$ to $\R.$ 


As a consequence of the above description of the geometry of $\Phi$ in
terms of its faces, particularly (\ref{PARAEXT2}), we conclude that
%
\begin{equation}\label{INTREPR}
\Phi= \closure \biggl( \biggl [ \bigcup_{f \in\Faces_{d-1}(\Phi)}
\Pi^{\downarrow}[f]  \biggr]^c  \biggr) =
\bigcap_{f \in\Faces_{d-1}(\Phi)} \closure ( [\Pi
^{\downarrow}[f]]^c  ),
\end{equation}
with $\closure(\cdot)$ standing for the topological closure, and with
$(\cdot)^c$
denoting the complement in $\R^{d-1} \times\R_+.$ This is the
parabolic counterpart
to the standard fact that a convex polytope is the
intersection of closed half-spaces determined by its
$(d-1)$-dimensional faces and containing
${\bf0}.$ From (\ref{INTREPR}) it follows that for each point $x \notin\Phi$,
there exists a translate of $\Pi^{\downarrow}$ containing $x$, but
not intersecting
$\Phi$, hence in particular not intersecting $\P,$ which is the
paraboloid version of
the standard separation lemma of convex geometry. On the other hand, if
$x$ is contained
in a translate of $\Pi^{\downarrow}$ not hitting $\P$, then $x \notin\Phi.$ Consequently
%
\begin{eqnarray}\label{PHIDEF2}
\Phi&=&
 \biggl[ \bigcup_{x \in\R^{d-1} \times\R_+,
[x \oplus\Pi^{\downarrow}] \cap\P= \varnothing} x \oplus\Pi
^{\downarrow}  \biggr]^c \nonumber
\\[-8pt]
\\[-8pt]
&=&\bigcap_{x \in\R^{d-1} \times\R_+,  [x \oplus\Pi^{\downarrow}]
\cap\P= \varnothing}
 [x \oplus\Pi^{\downarrow}  ]^c.
\nonumber
\end{eqnarray}
Alternatively, $\Phi$ arises as the complement
of the morphological opening of $\R^{d-1} \times\R_+ \setminus\P$
with downwards
paraboloid structuring element $\Pi^{\downarrow},$ that is to say,
\[
\Phi^c = [\P^c \ominus\Pi^{\downarrow}] \oplus\Pi^{\downarrow}
\]
with $\ominus$ standing for Minkowski erosion. In intuitive terms this means
that the complement of $\Phi$ is obtained by trying to \textit{fill} $\R
^{d-1} \times\R_+$
with downwards paraboloids $\Pi^{\downarrow}$ forbidden to hit any of
the Poisson
points in $\P$---the random open set obtained as the union of such
paraboloids
is precisely $\Phi^c.$

\begin{figure}[b]

\includegraphics{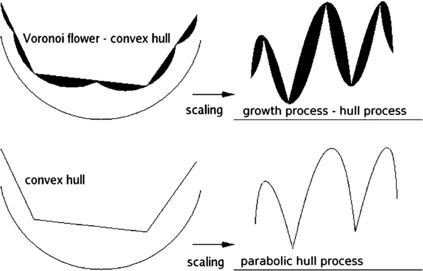}

  \caption{Convex hull, Voronoi flower and their scaling limits.}
\label{Fig1}
\end{figure}

To link the paraboloid hull and growth processes, note that
a point $x \in\P$ is a vertex of $\Phi$ if and only if $x \notin
\uphull(\P\setminus\{ x \}).$
By (\ref{PHIDEF2}) this means that $x \in\Vertices(\Phi)$ if and
only if there exists $y$
such that $[y \oplus\Pi^{\downarrow}] \cap\P= \{ x \}$ and, since
the set of
$y$ such that $y \oplus\Pi^{\downarrow} \ni x$ is simply $x \oplus
\Pi^{\uparrow},$
this condition is equivalent to having $x \oplus\Pi^{\uparrow}$ not entirely
contained in $[\P\setminus\{ x \}] \oplus\Pi^{\uparrow}$. In view of
(\ref{PSIEXTR}) this means that
%
\begin{equation}\label{EXTEQ}
\ext(\Psi) = \Vertices(\Phi).
\end{equation}

The theory developed in this section admits a particularly simple form
when $d=2.$ 
To see it, say that two points $x,y \in
\ext(\Psi)$ are neighbors in $\Psi,$ with notation $x \sim_{\Psi}
y$ or
simply $x \sim y,$ if and only if there is no point in $\ext(\Psi)$
with its spatial
coordinate between those of $x$ and $y.$ Then $\Vertices(\Phi) = \ext
(\Psi)$ as
in the general case, and $\Faces_1(\Phi) = \{ \Pi^{[\cdot]}[x,y],  x
\sim y \in\ext(\Psi) \}.$
In this context it is also particularly easy to display the
relationships between
the parabolic growth process $\Psi$ and the parabolic hull process
$\Phi$ in
terms of the analogous relations between the convex hull $K_\la$
and the Voronoi flower $F(\P_{\la})$ upon the transformation $T^\la$
at (\ref{SCTRANSF})
in large $\la$ asymptotics. To this end, see Figure~\ref{Fig1} and
note 
that, in large $\la$ asymptotics, we have:
%
%
%

$\bullet$ The extreme points in $\Psi,$ coinciding with $\Vertices(\Phi),$
correspond to the vertices of $K_\la.$

$\bullet$ Two points $x,y \in\ext(\Psi)$ are neighbors $x \sim y$ if and
only if
the corresponding vertices of $K_\la$ are adjacent, that is
to say, connected by an edge
of $\partial K_\la.$

$\bullet$ The circles $S^1(x/2,|x|/2)$ and $S^1(y/2,|y|/2)$ of two
adjacent vertices $x,y$ of $K_\la,$ whose pieces mark the boundary
of the Voronoi flower $F(\P_{\la}),$ are easily seen to
have their unique nonzero intersection point $z$ collinear
with $x$ and $y.$ Moreover, $z$ minimizes the distance to $0$
among the points on the $\overline{xy}$-line and $\overline{xy}
\,\bot\,\overline{0z}.$ For the parabolic processes this is
reflected by the fact that the intersection point of two upwards
parabolae with apices at two neighboring points $x$ and $y$
of $\Vertices(\Phi) = \ext(\Psi)$ coincides with the apex of the
downwards parabola $\Pi^{\downarrow}[x,y]$ as readily verified
by a direct calculation.

$\bullet$ Finally, relation (\ref{INTREPR}) becomes here
$\Phi= \bigcap_{x \sim y \in\ext(\Psi)} \closure ( [\Pi
^{\downarrow}[x,y]]^c  )
$ which is reflected by the fact that $K_\la$ coincides with the intersection
of all closed half-spaces containing ${\bf0}$ determined by
segments of the convex hull boundary $\partial K_\la.$

We conclude this paragraph by defining the paraboloid \textit{avoidance function}
$\hat\vartheta^{\infty}_k,  k \in\{1,2,\ldots,d\}.$
To this end, for each $x:= (v,h) \in\R^{d-1} \times\R_+$ let
$x^{\updownarrow} := \{ (v,h'),
h' \in\R\}$ be the infinite vertical ray (line) determined by $x$,
and let $A(x^{\updownarrow},k)$ be
the collection of all $k$-dimensional affine spaces in $\R^{d}$
containing $x^{\updownarrow},$
regarded as the asymptotic equivalent of the restricted Grassmannian
$G(\operatorname{lin}[x],k)$
considered in the definition (\ref{PAVOF}) of the nonrescaled
function $\vartheta_k^{\P_\la}.$
Next,
for $L \in A(x^{\updownarrow},k)$ we define the orthogonal paraboloid
surface $\Pi^{\bot}[x;L]$ to
$L$ at $x$ given by
%
\begin{eqnarray}\label{PARAORTHO}
&&\Pi^{\bot}[x;L]\nonumber\hspace*{-35pt}
\\[-8pt]
\\[-8pt] && \qquad := \biggl\{ x' = (v',h') \in\R^{d-1} \times\R,  (x-x')
\,\bot\, L,  h' = h - \frac{d(x,x')^2} {2} \biggr\}.\hspace*{-35pt}
\nonumber
\end{eqnarray}
Note that this is an analog of the usual orthogonal affine space
$L^{\bot} + x$ to $L$ at $x,$
with the second order parabolic correction typical in our asymptotic
setting---recall that
nonradial hyperplanes get asymptotically transformed onto downwards
paraboloids. Further,
for $L \in A(x^{\updownarrow},k)$, we put
%
\[
\vartheta^{\infty}_L(x) := {\bf1}(\{\Pi^{\bot}[x;L] \cap\Phi=
\varnothing\}).
\]
Observe that this is a direct analog of $\vartheta_L(x,\P_\la)$,
assuming the value
$1$ precisely when $x \notin K_{\la}|L \Leftrightarrow[L^{\bot} +
x] \cap K_{\la} = \varnothing.$
Finally, in full analogy to (\ref{PAVOF}) set
%
\begin{equation}\label{PARAVOF}
\vartheta^{\infty}_k(x) := \int_{A(x^{\updownarrow},k)} \vartheta
^{\infty}_L(x)\,d\mu^{x^{\updownarrow}}_k(L)
\end{equation}
with $\mu^{x^{\updownarrow}}_k$ standing for the normalized Haar
measure on $A(x^{\updownarrow},k);$
see page~591 in~\cite{SW}.\vspace*{6pt}

\textit{Duality relations between growth and hull processes.} %
As already signaled, there are close relationships between the
paraboloid growth and
hull processes, which we refer to as \textit{duality}. Here we discuss these
connections in more detail. The first observation is that
%
\begin{equation}\label{HULLTOGROWTH}
\Psi= \Phi\oplus\Pi^{\uparrow} = \Vertices(\Phi) \oplus\Pi
^{\uparrow}.
\end{equation}
This is verified either directly by the construction of $\Phi$ and
$\Psi,$
or, less directly but more instructively, by using the fact,
established in detail
in Section~\ref{LocScal} below, that $\Phi$ arises as the scaling
limit of $K_{\la}$,
whereas $\Psi$ is the scaling limit of the Voronoi flower %
\[
F(\P_{\la}) = \bigcup_{x \in\P_{\la}} B_d\biggl(\frac{x} {2}, \frac
{|x|} {2}\biggr) =
\bigcup_{x \in\Vertices(K_{\la})} B_d\biggl(\frac{x} {2}, \frac{|x|}
{2}\biggr),
\]
defined at (\ref{flower})
and then by noting that the balls
$B_d(x/2,|x|/2)$ asymptotically either scale into upward paraboloids or
they ``disappear at
infinity'';
see the proof of Theorem~\ref{LOCALT} below,
and recall that the support function of $K_{\la}$ coincides with the
radius-vector function
of $F(K_{\la})$ as soon as ${\bf0} \in K_{\la}$ (which, recall,
happens with overwhelming
probability).
Thus, it is straightforward to transform $\Phi$ into~$\Psi.$ To construct
the dual transform, say that $v \in\R^{d-1}$ is an \textit{extreme
direction} for $\Psi$ if
$\partial\Psi$ admits a local maximum at $v.$ Further, say that $x
\in\partial\Psi$ is
an extreme directional point for $\Psi,$ written $x \in\extdir(\Psi
),$ if and only if
$x = (v,\partial\Psi(v))$ for some extreme direction $v.$ Then we have
%
\begin{equation}\label{GROWTHTOHULL}
\Phi^c = \Psi^c \oplus\Pi^{\downarrow}   \quad  \mbox{and} \quad
\closure(\Phi^c) = \extdir(\Psi) \oplus\Pi^{\downarrow}.
\end{equation}
Again, this can be directly proved, yet it is more appealing to observe
that this statement
is simply an asymptotic counterpart of the usual procedure of restoring
the convex polytope
$K_{\la}$ given its support function. Indeed, the complement of the
polytope arises as the union
of all half-spaces of the form $H_x := \{ y \in\R^d,  \langle y-x, x
\rangle\geq0 \}$
(asymptotically transformed onto suitable translates of $\Pi
^{\downarrow}$ under the action
of $T^{\la},  \la\to\infty$) with $x$ ranging through $x = ru,  r
> h_{K_{\la}}(u),
r \in\R,  u \in\S^{d-1}$ which corresponds to taking $x$ in the
epigraph of
$h_{K_{\la}}$ (transformed onto $\Psi^c$ in our asymptotics). This
explains the first
equality in (\ref{GROWTHTOHULL}). The second one comes from the fact
that it
is enough in the above procedure to consider half-spaces $H_x$ for $x$
in extreme
directions only, corresponding to directions orthogonal to
$(d-1)$-dimensional faces
of $K_{\la}$ and marked by local minima of the support function
$h_{K_{\la}}$
(asymptotically mapped onto local maxima of $\partial\Psi$). It is
worth noting that all
extreme directional points of $\Psi$ arise as $d$-fold intersections
of boundaries of
upwards paraboloids $\partial[x \oplus\Pi^{\uparrow}],  x \in\ext
(\Psi)$, although
not all such intersections give rise to extreme directional points
[they do so
precisely when the apices of $d$ intersecting upwards paraboloids are
vertices of the same
$(d-1)$-dimensional face of $\Phi,$ which is not difficult to prove
but which is not needed
here].

\section{Local scaling limits}\label{sec4}\label{LocScal}

The re-scaled processes $\hat{s}_{\la}$ and $\hat{r}_{\la}$,
defined at (\ref{SRESC}) and (\ref{RRESC}), respectively, are
locally parabolic, and here we show that 
their graphs have scaling limits given by
the boundaries of the paraboloid growth processes $\Psi$ and $\Phi$,
respectively. Recall from Definition~\ref{fund-def} that both $\Psi$
and $\Phi$ are defined in terms of $\P$, the Poisson point
process in $\R^{d-1} \times\R_+$ with intensity
density $h^{\delta}\,dh\,dv.$ Recall
that $B_d(x, r)$ stands for the $d$-dimensional radius $r$ ball
centered at $x$.


\begin{theo}\label{LOCALT}
For any $R>0,$ the random functions $\hat{s}_{\la}$ and
$\hat{r}_{\la}$ converge in law as $\la\to\infty$ to $\partial
\Psi$ and
$\partial\Phi$, respectively, in the space ${\cal C}(B_{d-1}(\mathbf{0},R))$ of continuous
functions on $B_{d-1}({\bf0},R)$ endowed with the supremum norm.
\end{theo}

\begin{remark*} Theorem~\ref{LOCALT} adds to Molchanov \cite{Mo}, who
establishes convergence of the nonrescaled process $nr( \cdot ,
\{X_i\}_{i=1}^n)$ in $\S^{d-1} \times\R$, where $X_i$ are i.i.d.
uniform in $\B^d$.
It also adds to Eddy \cite{Eddy}, who considers convergence of the
properly scaled defect support function for i.i.d. random variables
with a circularly symmetric standard Gaussian distribution.
\end{remark*}

\begin{pf*}{Proof of Theorem~\ref{LOCALT}}
The convergence in law for $\hat{s}_{\la}$ may be shown to follow
from the more
general theory of generalized growth processes developed in \cite{SY},
but we
provide here an argument specialized to our present set-up. Recall that
we place ourselves on the event that $\0 \in K_\la$ which is exponentially
unlikely to fail as $\la\to\infty$ and thus, for our purposes, may be
assumed to hold without loss of generality. Further, the
support function $h_{\{x\}}\dvtx  \S^{d-1} \to\R$ of a point $x \in\B
^d$ is given
for all $u \in\S^{d-1}$ by
$ h_{\{x \}}(u) = |x| \cos(d_{\S^{d-1}}(u,x/|x|))$
with $d_{\S^{d-1}}$ standing for the geodesic distance in
$\S^{d-1}.$ 


Recall that $\P^{(\la)} := T^\la(\P_\la)$, where
$T^{\la}$ is defined at (\ref{SCTRANSF}) and where $\P^{(\la)}$ has
density given by (\ref{DENS1}). Write $x:= (v_x, h_x)$ for the
points in $\P^{(\la)}.$ Under $T^{\la}$ we may write
$\hat{s}_\la(v),  v \in\la^{\beta} \B_{d-1}(\pi),$ as
%
\begin{eqnarray}\label{sform}
\hat{s}_\la(v) &=&
\la^{\gamma}  \Bigl(1-\max_{x=(v_x,h_x) \in\P^{(\la)}}
[1-\la^{-\gamma}h_x]\nonumber\\
&&\hphantom{\la^{\gamma}  \Bigl(1-\max_{x=(v_x,h_x) \in\P^{(\la)}}}{}\times [\cos[ d_{\S^{d-1}}(\exp_{d-1}(\la^{-\beta}
v),\exp_{d-1}(\la^{-\beta} v_x))]] \Bigr)
\nonumber\\
&=& \la^{\gamma} \min_{x \in\P^{(\la)}} \bigl[ 1 - (1-\la^{-\gamma
}h_x)\nonumber
\\[-8pt]
\\[-8pt]
&&\hphantom{\la^{\gamma} \min_{x \in\P^{(\la)}} \bigl[ 1 {}-{}}{}\times\bigl(1-
\bigl(1-\cos[ d_{\S^{d-1}}(\exp_{d-1}(\la^{-\beta} v),\exp_{d-1}( \la
^{-\beta} v_x))]\bigr)\bigr)  \bigr]
\nonumber
\\
 &=& \min_{x \in\P^{(\la)}} \bigl [ h_x +
\la^{\gamma} \bigl(1-\cos(
d_{\S^{d-1}}(\exp_{d-1}(\la^{-\beta} v),\exp_{d-1}( \la^{-\beta}
v_x)))\bigr)
\nonumber\\
&&\hspace*{12pt}\hphantom{\min_{x \in\P^{(\la)}} \bigl [}{}-  h_x\bigl(1-\cos[ d_{\S^{d-1}}(\exp_{d-1}(\la^{-\beta} v),\exp
_{d-1}( \la^{-\beta} v_x))]\bigr)  \bigr].
\nonumber
\end{eqnarray}

Thus, by (\ref{FunctionS}) and (\ref{SRESC}), the graph of
$\hat{s}_{\la}$ coincides with the lower boundary of the following
\textit{generalized growth process}
%
\begin{equation}\label{gengrowth} \Psi^{(\la)}:= \bigcup_{x \in\P
^{(\la)}}
[\Pi^{\uparrow}]_x^{(\la)},
\end{equation}
where for $x:= (v_x, h_x)$ we have
%
\begin{eqnarray}\label{form4.2} [\Pi^{\uparrow}]_x^{(\la)} &=&\bigl \{
(v,h) \in
\R^{d-1} \times\R_+,
h \geq h_x \nonumber
\\[-8pt]
\\[-8pt]
&&\hphantom{\bigl \{}{}+ \la^{\gamma} \bigl(1-\cos[e_\la(v, v_x)]\bigr) - h_x\bigl(1- \cos
[e_\la(v,
v_x)]\bigr) \bigr\},
\nonumber
\end{eqnarray}
with
%
\begin{equation}\label{defe} e_\la(v, v_x):=
d_{\S^{d-1}}(\exp_{d-1}(\la^{-\beta} v),\exp_{d-1}(\la^{-\beta}
v_x)).
\end{equation}

We now show for fixed $R \in(0, \infty)$ that the lower boundary
of the process $\Psi^{(\la)}$ converges in law to $\partial\Psi$ in
the space ${\cal C}(B_{d-1}({\bf0},R))$. This goes as follows.
With $R$ fixed, for all $H \in\R^+$ and $\la\in\R^+$, let
$E_1(R,H, \la)$ be the event that the heights of the lower
boundaries of $\Psi$ and $\Psi^{(\la)}$ are at most $H$ over the
spatial region $B_{d-1}({\bf0},R)$. Interpreting the boundary
$\partial\Psi^{(\la)}$ as the graph of a function from $\R^{d-1}$
to $\R_+$, it follows from straightforward modifications of Lemma
3.2 in \cite{SY} that there is a $\la_0 \in(0, \infty)$ such that,
uniformly for $\la\geq\la_0$, we have
%
\begin{equation}\label{SuperexpBd}
  P\Bigl[\sup_{v \in B_{d-1}({\bf0},R)} \partial\Psi^{(\la)}(v) \geq H\Bigr]
\leq C(R) \exp\bigl(-c \bigl[H^{(d+1)/2} \wedge R^{d-1} H^{1+\delta}\bigr]\bigr)\hspace*{-35pt}
\end{equation}
with $c > 0$ and $C(R) < \infty$ (note that the extra term
$R^{d-1} H^{1+\delta}$ in the exponent corresponds to the
probability of having $B_{d-1}({\bf0},R) \times[0,H]$ devoid of
points of $\P$ and $\P^{(\la)}$). Lemma 3.2 in \cite{SY} likewise
gives a similar bound for $P[\sup_{v \in B_{d-1}({\bf0},R)}
\partial\Psi(v) \geq H]$. Thus $P[E_1(R, H, \la)^c]$ decays
exponentially fast in $H$, uniformly in $\la$ and it is enough to
show, conditional on $E_1(R, H, \la)$, that $\hat{s}_\la( \cdot)$
converges to $\partial\Psi$ in the space ${\cal C}(B_{d-1}(\mathbf{0},R))$.

Next, with $H$ fixed, observe that for each $R$ there exists a
constant $R':=R'(R,H)$ such that for all $\la$ large enough, the
behavior of $\Psi^{(\la)}$ and $\Psi$ restricted to $B_{d-1}(\mathbf{0},R) \times[0,H]$ only depends on the restriction to $B_{d-1}(\mathbf{0},R') \times[0,H]$ of the processes $\P^{(\la)}$ and $\P$,
respectively. For instance in the case of $\Psi$ it is enough that
the region $B_{d-1}({\bf0},R') \times[0,H]$ contain the apices of
all translates of $\Pi^{\uparrow}$ which hit $B_{d-1}({\bf0},R)
\times[0,H],$ that is to say, the choice $R' := R + \sqrt{2H}$ will
suffice.

We also assert for these fixed $R'$ and $H$ that $\P$ and
$\P^{(\la)}$ can be coupled on a common probability space so that on
a set $E_2(R', H, \la)$, with $P[E_2(R', H, \la)] \to1$ as $\la\to
\infty$, their restrictions agree on $B_{d-1}({\bf0},R') \times
[0,H]$. This assertion, referred to as ``total variation convergence
on compact sets,'' follows by combining Theorem 3.2.2 in \cite{REI},
which upper bounds total variation distance between Poisson measures
by a multiple of the $L^1$ norm of the difference of their
densities, with the observation that the intensity density of
$\P^{(\la)}$, as given by (\ref{DENS1}), converges in
$L^1(B_{d-1}({\bf0},R') \times[0,H])$ to the intensity density of
$\P$, as given by (\ref{basicPPP}).

Let $E(R, H, \la):= E_1(R, H, \la) \cap E_2(R', H, \la)$ and note
that $P[E(R, H, \la)] \to1$ as $\la\to\infty$. It is enough to
show, conditional on the event $E(R, H, \la)$, that $\hat{s}_\la(
\cdot)$ converges to $\partial\Psi$ in the space ${\cal
C}(B_{d-1}({\bf0},R))$.

Now we examine the lower boundary of $\P^{(\la)}$ given the event
$E(R, H, \la)$. On this event we have
\[
\Psi^{(\la)}:= \bigcup_{x \in\P} [\Pi^{\uparrow}]_x^{(\la)}
\]
with $[\Pi^{\uparrow}]_x^{(\la)}$ given by (\ref{form4.2}).
Recalling the definition of $e_\la(v, v_x)$ at (\ref{defe}) and
recalling $\gamma= 2 \beta$ from (\ref{EQBG2}) we have (using that
the ratio of the Euclidean norm and geodesic norm converges to $1$)
\[
\lambda^{\gamma} (e_\la(v, v_x))^2 = \lambda^{\gamma}  \biggl(
\frac{
e_\la(v, v_x) } { |\la^{-\beta}v - \la^{-\beta}v_x|}  \biggr)^2
|\la^{-\beta}v - \la^{-\beta}v_x|^2 \to| v - v_x|^2.
\]
Using
the Taylor expansion of the cosine function up to second order in
(\ref{form4.2}), it follows that on $E(R, H, \la)$ the graph of the
lower boundary of $[\Pi^{\uparrow}]_x^{(\la)}, x \in\P$, converges
with respect to the sup norm distance on $B_{d-1}({\bf0},R') \times
[0,H])$ to the graph of the lower boundary of the paraboloid $v
\mapsto h_x + |v-v_x|^2/2$, that is to say, the lower boundary of
$x \oplus\Pi^{\uparrow}.$ In the space ${\cal
C}(B_{d-1}({\bf0},R))$ the lower boundary of $\Psi^{(\la)}$ is with
probability one determined by a finite number of
$[\Pi^{\uparrow}]_x^{(\la)}$ and thus as $\la\to\infty$,
$\hat{s}_{\la}$ converges in law to $\partial\Psi$, as claimed.
This shows Theorem~\ref{LOCALT} for $\hat{s}_{\la}$.

To prove Theorem~\ref{LOCALT} for $\hat{r}_{\la}$, consider the
spherical cap
%
\begin{eqnarray}\label{scap}
\operatorname{cap}_{\la}[v^*,h^*] :=
\{ x \in\B^d,   \langle x, \exp_{d-1}(\la^{-\beta}v^*) \rangle
\geq1 - \la^{-\gamma} h^*
\}, \nonumber
\\[-8pt]
\\[-8pt]  \eqntext{(v^*,h^*) \in\R^{d-1} \times\R_+,}
\end{eqnarray}
and note that with $x:= (|x|,u) \in\B^d$, we equivalently have
\[
\operatorname{cap}_{\la}[v^*,h^*]: =  \biggl\{ x \in\B^d,  (1 -
|x|) \leq
\max\biggl(0,1- \frac{(1 - \la^{-\gamma} h^*)} { \cos\theta} \biggr)
 \biggr\},
\]
where
$\theta$ denotes the angle between $x$ and
$\exp_{d-1}(\la^{-\beta} v^*)$. Under the transformation
$T^{\la}$ the cap transforms into
\begin{eqnarray*}
&&\operatorname{cap}^{(\la)}[v^*,h^*] \\
&& \qquad :=
 \biggl\{ (v,h) \in{\cal R}_\la,  \\
 &&\hphantom{:=
 \biggl\{} \qquad  h \leq\la^{\gamma} \max \biggl(0,1-
\frac{1-\la^{-\gamma}h^*} {\cos(
d_{\S^{d-1}}(\exp_{d-1}(\la^{-\beta}v),\exp_{d-1}(\la^{-\beta
}v^*)) )}  \biggr)  \biggr\}
\\
&& \qquad =  \biggl\{ (v,h) \in{\cal R}_\la,   h \leq\la^{\gamma}
\max \biggl(0,1- \frac{1-\la^{-\gamma}h^*} {\cos(e_\la(v, v^*) )}
 \biggr)  \biggr\},
\end{eqnarray*}
where $e_\la(v, v^*)$ is as in (\ref{defe}).

Using that $\B^d \setminus K_\la$ is the union
of all spherical caps not hitting any of the points in $\P_{\la},$ we
conclude that under
the mapping $T^{\la}\dvtx  \P_\la\to\P^{(\la)},$ the complement of
$K_\la$ in $\B^d$
gets transformed into the union
%
\begin{equation}\label{TRANSUNION}
\bigcup \bigl\{ \operatorname{cap}^{(\la)}[v^*,h^*],
(v^*,h^*) \in{\cal R}_{\la},
\operatorname{cap}^{(\la)}[v^*,h^*] \cap\P^{(\la)} = \varnothing
 \bigr\}.
\end{equation}
Let the \textit{paraboloid hull process}
$\Phi^{(\la)}$ be the complement of this union in $\R^{d-1} \times
\R_+,$ that is,
%
\begin{equation}\label{phullprocess}
  \Phi^{(\la)} :=  \Bigl( \bigcup\bigl\{ \operatorname{cap}^{(\la
)}[v^*,h^*],
(v^*,h^*) \in{\cal R}_{\la},
\operatorname{cap}^{(\la)}[v^*,h^*] \cap\P^{(\la)} = \varnothing\bigr\}
 \Bigr)^c.\hspace*{-35pt}
\end{equation}

To prove the asserted convergence of $\hat{r}_\la$, we modify the
approach given for the convergence of $\hat{s}_\la$. Let $F_1(R,H,
\la)$ be the event that the heights of the lower boundaries of
$\Phi$ and $\Phi^{(\la)}$ are at most $H$ over the spatial region
$B_{d-1}({\bf0},R)$. As in (\ref{SuperexpBd}) we get that
$P[F_1(R,H, \la)^c]$ decays exponentially fast in $H$, uniformly in
$\la$, implying that it is enough to show, conditional on $F_1(R,H,
\la)$, that $\hat{r}_\la$ converges to $\partial\Phi$ in
$C(B_{d-1}({\bf0},R))$.

Both $\Phi$ and $\Phi^{(\la)}$ are locally determined in the sense
that for any $R,H,\epsilon> 0$ there exist $R'',H''
> 0$, such that, with probability at least $1-\epsilon,$ the restrictions
of $\Phi$ and $\Phi^{(\la)}$ to $B_{d-1}({\bf0},R) \times[0,H]$
are determined by the restrictions to $B_{d-1}({\bf0},R'') \times
[0,H'']$ of $\P^{(\la)}$ and $\P$, respectively. Indeed if the
geometry of $\Phi$ within $B_{d-1}({\bf0},R) \times[0,H]$ were
affected by the status of a point $x \in\R^{d-1} \times\R_+,$
there would exist a translate of $\Pi^{\downarrow}$ such that the
translate:
(i)~hits $B_{d-1}({\bf0},R) \times[0,H];$
(ii) contains $x$ on its boundary; (iii) is devoid of other points of
$\P.$
Thus the probability of such an influence being exerted by a faraway
point $x$ tends to $0$ with the distance of $x$ from $B_{d-1}(\mathbf{0},R) \times[0,H].$ The argument for $\Phi^{(\la)}$ and
$\P^{(\la)}$ is analogous. Statements of this kind, going under the
general name of stabilization, shall be discussed in more detail in
Lemma~\ref{LocalLem} below. 

As above, we may couple $\P$ and $\P^{(\la)}$ on a common
probability space so that their restrictions to $B_{d-1}(\mathbf{0},R'')
\times[0,H'']$ agree on a set $F_2(R'',H'', \la)$, with
$P[F_2^c(R'',H'', \la)] \to0$ as $\la\to\infty$. Put $F(R,H,
\la):= F_1(R,H, \la) \cap\break F_2(R'',H'', \la)$, and note that $P[F(R,H,
\la)] \to1$ as $\la\to\infty$. We now show on the event $F(R, H,
\la)$ that $\hat{r}_\la( \cdot)$ converges to $\partial\Psi$ as
$\la\to\infty$.

We Taylor-expand the cosine function up to second order to get that
\[
\operatorname{cap}^{(\la)}[v^*,h^*] = \biggl \{ (v,h) \in{\cal R}_\la
,  h
\leq
\max \biggl( 0, \la^\gamma- \frac{ \la^\gamma- h^*} { 1 - e_\la(v,
v^*)^2/2 +\cdots }  \biggr)  \biggr\}.
\]
Using the convergence
$\la^\gamma e_\la^2(v, v^*) \to|v - v^*|^2$
and the expansion $1/(1 - r) = 1 + r + r^2 +\cdots $ for $r$ small, we see that the upper boundary of $\operatorname
{cap}^{(\la)}[v^*,h^*]$ converges as $\la\to\infty$ with respect to
the sup norm distance on $B_{d-1}({\bf0},R) \times[0,H]$ to the
graph of the upper boundary of the paraboloid
\[
\biggl\{(v, h) \in\R^{d-1} \times\R_+,  h \leq h^* - \frac{|v -
v^*|^2}{2} \biggr\},
\]
that is, the graph of the upper boundary of $(v^*,h^*) \oplus
\Pi^{\downarrow}$. In the space ${\cal C}(B_{d-1}({\bf0},\allowbreak R))$ the
upper boundary of $\Phi^{(\la)}$ is with probability one determined
by a finite
number of $[\Pi^{\downarrow}]_x^{(\la)}$.

This observation, the definition of $\hat{r}_{\la}$, and the
relation (\ref{TRANSUNION}), show that $\hat{r}_{\la}$ converges in
law in the space ${\cal C}(B_{d-1}({\bf0},R))$ equipped with the
supremum norm to the continuous function determined by the upper
boundary of the process
\[
\bigcup_{x \in\R^{d-1} \times\R_+, [x \oplus\Pi^{\downarrow}] \cap\P= \varnothing} x \oplus
\Pi^{\downarrow},
\]
which coincides with $\partial\Phi$ in view of
(\ref{PHIDEF2}). This completes the proof of Theorem~\ref{LOCALT}.
\end{pf*}

\section{Exact distributional results for scaling limits}\label{sec5}
\label{PC1}

 This section is restricted to dimension $d = 2$ and to the
homogeneous Poisson point process in the unit-disk. Here we provide
explicit formulae for the fidis of the processes
$\hat{s}_{\lambda}$ and $\hat{r}_{\lambda}$ and give their explicit
asymptotics, confirming a posteriori the existence of the limiting
parabolic growth and hull processes of Section~\ref{g+hullprocess}.
\subsection{\texorpdfstring{The process $\hat{s}_{\lambda}$}{The process s lambda}}\label{sec5.1}

This subsection calculates the distribution of\break $s(\theta_0,{\cal
P}_{\lambda})$ and establishes the convergence of the fidis of
both the process and its re-scaled version. Throughout this
section we identify the unit sphere $\S^1$ with the segment
$[0,2\pi),$ whence the notation $s(\theta,\cdot),  \theta\in
[0,2\pi),$ and likewise for the radius-vector function
$r(\theta,\cdot).$ A first elementary result is the following:
%
\begin{lemm} 
For every $h>0$, $u\in\S^1$ and $\lambda>0$, we have
\[
P[s(u,{\mathcal P}_{\lambda})\ge h]=\exp\bigl \{-\lambda
\bigl(\arccos(1-h)- (1 - h)\sqrt{2h-h^2} \bigr) \bigr\}.
\]
\end{lemm}

\begin{pf} Notice that $(s(u,{\mathcal P}_{\lambda})\ge h)$
is equivalent to $\operatorname{cap}_1[u,h] \cap{\mathcal
P}_{\lambda} = \varnothing$, where $\operatorname{cap}_1[u,h]$ is
defined at (\ref{scap}).
Since the Lebesgue measure $\ell({\operatorname{cap}}_1[u,h])$ of
${\operatorname{cap}}_1[u,h]$ satisfies
%
\begin{equation}\label{eq:1}
\ell({\operatorname{cap}}_1[u,h])
=\arccos(1-h) -(1-h)\sqrt{2h-h^2},
\end{equation}
the lemma follows by the Poisson property of the process ${\cal
P}_{\lambda}$.
\end{pf}

  We focus on the asymptotic behavior of the process $s$
when $\lambda$ is large. When we scale in space, we obtain the
fidis of white noise and when we scale in both time and space to
get $\hat{s}$, we obtain the fidis of the parabolic growth process
$\Psi$ defined in Section~\ref{g+hullprocess}. Let $\N$
denote the positive integers.
In dimension two, by the representation \eqref{expomap}, we notice
that the exponential map obtained for the choice $u_0=(0,1)$ has
the following basic expression:
\[
\exp_1(\theta)=(\sin(\theta),\cos(\theta)),\qquad\theta\in
{\mathbb R}.
\]

%
\begin{prop} \label{fidis}
Let $n\in\N$, $0\le
\theta_1<\theta_2<\cdots<\theta_n<2\pi$ and $h_i \in(0,\infty)$
for all $i =1,\ldots,n$. 
Then
\begin{eqnarray*}
&&\lim_{\la\to\infty} P[\lambda^{2/3}s(\exp_1(\theta_1),{\mathcal
P}_{\lambda})\ge h_1;\ldots; \lambda^{2/3}
s(\exp_1({\theta_n}),{\mathcal P}_{\lambda})\ge h_n ] \\
&& \qquad = 
\prod_{k=1}^n \exp \biggl\{-\frac{4\sqrt{2}}{3}
h_k^{3/2} \biggr\}.
\end{eqnarray*}
Moreover, for every $v_1<v_2<\cdots<v_n\in\R$, we have
\begin{eqnarray*}
&& \lim_{\la\to\infty}
P[\lambda^{2/3}s(\exp_1(\lambda^{-1/3}v_1),{\mathcal P}_{\lambda
})\ge
h_1;\ldots; \lambda^{2/3}
s(\exp_1(\lambda^{-1/3}v_n),{\mathcal P}_{\lambda})\ge h_n]\\
&& \qquad
= \exp\biggl (-\int_{\inf_{1\le i\le
n}(v_i-\sqrt{2h_i})}^{\sup_{1\le i\le n}(v_i+\sqrt{2h_i})}\sup
_{1\le
i \le n}\biggl[\biggl(h_i-\frac{1}{2}(u-v_i)^2\biggr)\\
&& \qquad \hphantom{= \exp\biggl (-\int_{\inf_{1\le i\le
n}(v_i-\sqrt{2h_i})}^{\sup_{1\le i\le n}(v_i+\sqrt{2h_i})}\sup
_{1\le
i \le n}\biggl[}{}\times{\bf1} \bigl(|u-v_i|\le
\sqrt{2h_i}\bigr)\biggr]\,du \biggr).
\end{eqnarray*}
\end{prop}

\begin{pf}
The first assertion is obtained by noticing that the  events\break
$\{s(\exp_1(\theta_i), {\cal P}_{\lambda})\ge
\lambda^{-2/3}h_i\},   1 \leq i \leq n,$ are independent as soon
as $h_i \in\break(0, \frac{\lambda^{2/3}}{2}\min_{1\le k\le
n}(1-\cos(\theta_{k+1}-\theta_k)))$.
We then apply Lemma~\ref{PClem} to
estimate the probability of each of these events. Let us recall
beforehand that
$\arccos(1-x)$ is expanded as $\sqrt{2x}+\sqrt{2x^3}/12+\cdots$
when $x\to0$.
For every $1\le i\le n$, we have
\begin{eqnarray*}
&&-\log P[\lambda^{2/3}s(\exp_1(\theta_i),{\mathcal P}_{\lambda})\ge
h_i]\\
&& \qquad = \lambda\bigl [\arccos(1-\la^{-2/3}h_i)- (1 - \la
^{-2/3}h_i)\sqrt{2\la^{-2/3}h_i-\la^{-4/3}h_i^2} \bigr]\\
&& \qquad \mathop{=}\limits_{\la\to\infty}\la \biggl[\sqrt{2}\lambda
^{-1/3}\sqrt{h_i}+\frac{\sqrt{2}}{12}\lambda^{-1}h_i^{3/2} \\
&& \qquad   \hphantom{\mathop{=}\limits_{\la\to\infty}\la \biggl[} {}
-(1-\la^{-2/3}h_i)\sqrt{2}\la^{-1/3}\sqrt
{h_i} \biggl(1-\frac{1}{4}\la^{-2/3}h_i \biggr)+o(\la^{-1}) \biggr]\\
&& \qquad \mathop{=}\limits_{\la\to\infty} \biggl(\frac{1}{12}+1+\frac
{1}{4} \biggr)\sqrt{2}h_i^{3/2}+o(1)\\
&& \qquad \mathop{=}\limits_{\la\to\infty}\frac{4}{3}\sqrt{2}h_i^{3/2}+o(1).
\end{eqnarray*}
Here and elsewhere in this section, the terminology $f(\la)
\mathop{\sim}\limits_{\lambda\to\infty} g(\la)$ [resp., $f(\lambda)\mathop{=}
\limits_{\la\to\infty}o(g(\lambda))$] signifies that $\lim_{\la\to\infty}f(\la
)/ g(\la) = 1$ [resp., $\lim_{\la\to\infty}f(\la)/\allowbreak g(\la)=0$].
For the second assertion, it suffices to determine the area
$\ell({\mathcal D}_n)$ of the domain
\[
{\mathcal
D}_n:=\bigcup_{1\le i\le
n} {\operatorname{cap}}_{\lambda} [v_i,h_i].
\]
%
This set is contained in the angular sector between
$\alpha_n:=\inf_{1\le i\le
n}[\lambda^{-1/3}v_i-\arccos(1-\lambda^{-2/3}h_i)]$ and $\beta
_n:=\sup_{1\le i\le
n}[\lambda^{-1/3}v_i+\arccos(1-\lambda^{-2/3}h_i)]$.
Denote by $\rho_n(\cdot)$ the radial function which associates to
$\theta$ the distance between the origin and the point in ${\mathcal
D}_n$ closest to the origin lying on the half-line making angle
$\theta$ with the positive $x$-axis. Then
\begin{eqnarray*}
\ell\bigl({\mathcal
D}_n\bigr)&=&\int_{\alpha_n}^{\beta_n}\frac{1}{2}\bigl(1-\rho_n^2(\theta
)\bigr)\,d\theta\\
&=&
\lambda^{-1/3}\int_{\lambda^{1/3}\alpha_n}^{\lambda^{1/3}\beta
_n}\frac{1}{2}\bigl(1-\rho_n^2(\lambda^{-1/3}u)\bigr)\,du\\
&\mathop{\sim}\limits_{\lambda\to\infty}& \lambda^{-1/3}\int_{\inf
_{1\le
i\le
n}
(v_i-\sqrt{2h_i})}^{\sup_{1\le i\le n}(v_i+\sqrt{2h_i})}
\bigl(1-\rho_n(\lambda^{-1/3}u)\bigr)\,du.
\end{eqnarray*}
Each set
${\operatorname{cap}}_{\lambda}[v_i,h_i]$ is bounded by a line with the
polar equation
\[
\rho=\frac{1-\lambda^{-2/3}h_i}{\cos(\theta-\lambda^{-1/3}v_i)}.
\]
Consequently, the function $\rho_n(\cdot)$ satisfies, for every
$\theta\in
(0,2\pi)$,
\begin{eqnarray*}
1-\rho_n(\theta)&=&\sup_{1\le i\le
n} \biggl[
\frac{\cos(\theta-\lambda^{-1/3}v_i)-1+\lambda^{-2/3}h_i}{\cos
(\theta-\lambda^{-1/3}v_i)}\\
&&\hphantom{\sup_{1\le i\le
n} \biggl[}{}\times
{\bf1} \bigl(|\theta-\lambda^{-1/3}v_i|\le\arccos
(1-\lambda^{-2/3}h_i) \bigr)  \biggr].
\end{eqnarray*}
It remains to determine the
asymptotics of the above function. We obtain that
\[
1-\rho_n(\lambda^{-1/3}u)\mathop{\sim}\limits_{\lambda\to\infty}\lambda^{-2/3}\sup_{1\le
i \le n} \biggl[\biggl(h_i-\frac{1}{2}(u-v_i)^2\biggr){\bf1} \bigl(|u-v_i|\le
\sqrt{2h_i} \bigr) \biggr].
\]
\begin{figure}

\includegraphics{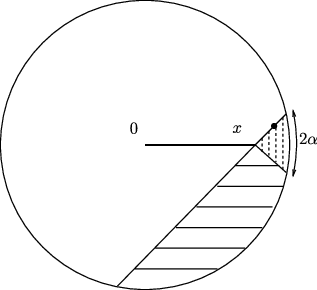}

\caption{When is a point included in the convex hull?}
\label{picture}
\end{figure}
\noindent \hspace*{-2pt}Considering that the required probability is equal to $\exp(-\lambda
\ell({\cal
D}_n))$, we complete the proof.
\end{pf}

\begin{rema}\label{rem1} Proposition~\ref{fidis} could have been obtained
through the use of the growth process $\Phi$. 
Indeed, we have $\partial\Phi(v_i)$
greater than $h_i$ for every $1\le i\le n$ if and only if none of the points
$(v_i,h_i)$ is covered by a parabola of $\Phi$. 
Equivalently, this means that there is no point of ${\cal P}$ in the
region arising as union of translated downward parabolae
$\Pi^{\downarrow}$ 
with peaks at $(v_i,h_i)$. Calculating the area of this region
yields Proposition~\ref{fidis}.
\end{rema}


\subsection{\texorpdfstring{The process $\hat{r}_{\lambda}$}{The process r lambda}}\label{sec5.2}

This subsection, devoted to distributional results for
$\hat{r}_{\lambda}$, follows the same lines as the previous one. The
problem of determining the distribution of $r(\cdot,{\cal
P}_{\lambda})$ seems to be a bit more tricky.\vadjust{\goodbreak} To proceed, we fix a
direction $u\in\S^{1}$ and a point $x=(1-h)u$
($h\in[0,1]$) inside the unit-disk $\B^2$.
Consider an angular sector centered at $x$ and opening
away from the origin. Open the sector until it first meets a point
of the Poisson point process at the angle ${\mathcal A}_{\lambda,h}$ (the
set with dashed lines must be empty in Figure~\ref{picture}).
Let ${\mathcal A}_{\lambda,h}$ be the minimal angle of opening from
$x=(1-h)u$ in order to meet a point of ${\mathcal
P}_{\lambda}$ in the opposite side of the origin. In particular,
when ${\mathcal A}_{\lambda,h}=\alpha$, there is no point of
${\mathcal P}_{\lambda}$ in
\[
{\cal S}_{\alpha,h}:=\{y\in{\B^2},\langle
y-x,u\rangle\ge\cos(\alpha)|y-x|\}.
\]
Consequently, we
have 
%
\begin{equation}
\label{eq:321}
P[{\mathcal A}_{\lambda,h}\ge
\pi/2]=P[s(u,{\cal
P}_{\lambda})\ge h].
\end{equation}
The next lemma provides the distribution of ${\mathcal A}_{\lambda,h}$.
%
\begin{lemm}\label{theta}
For every
$0\le\alpha\le\pi/2$ and $h\in[0,1]$, we have
%
\begin{equation}\label{eq:322}
P[{\mathcal A}_{\lambda,h}
\ge\alpha]=\exp \{-\lambda\ell({\cal S}_{\alpha,h}) \}
\end{equation}
with
%
\begin{eqnarray}
\label{eq:2}
\ell({\cal
S}_{\alpha,h})&=& \biggl(\alpha+\frac{(1-h)^2}{2}\sin(2\alpha
)-(1-h)\sin(\alpha)\sqrt{1-(1-h)^2\sin^2(\alpha)}\nonumber\hspace*{-30pt}
\\[-4pt]
\\[-12pt]
&&\hspace*{160.5pt}\hphantom{\biggl(}{}-\arcsin
\bigl((1-h)\sin(\alpha)\bigr) \biggr).
\nonumber\hspace*{-30pt}
\end{eqnarray}
When $\lambda$ goes to infinity, ${\mathcal A}_{\lambda,\lambda
^{-2/3}h}$ converges in
distribution to a measure with mass 0 on $[0,\pi/2)$ and mass
$(1-\exp\{-\frac{4\sqrt{2}}{3}h^{2/3}\})$ on $\{\pi/2\}$.
\end{lemm}
\begin{pf} A quick geometric consideration shows that the set ${\cal
S}_{\alpha,h}$ is seen from the origin with an
angle equal to
%
\begin{equation}\label{eq:323}
2\beta=2\bigl[\alpha-\arcsin\bigl((1-h)\sin(\alpha)\bigr)\bigr].
\end{equation}
To obtain (\ref{eq:2}), we first integrate in polar
coordinates, giving
\begin{eqnarray*}
\ell({\cal
S}_{\alpha,h})&=&
2\int_0^{\beta} \biggl[\int_{\fraca{\sin(\alpha-\gamma)}{\sin
(\alpha-\theta)}}^1\rho \,d\rho \biggr]\,d\theta\\
&=&
\int_0^{\beta} \biggl(1-\frac{(1-h)^2\sin^2(\alpha)}{\sin^2(\alpha
-\theta)} \biggr)\,d\theta\\
&=&\beta-(1-h)^2\sin^2(\alpha) \biggl(\frac{1}{\tan(\alpha-\beta
)}-\frac{1}{\tan(\alpha)} \biggr).
\end{eqnarray*}
We then use (\ref{eq:323}) to get (\ref{eq:2}).

Let us show now the last assertion of Lemma~\ref{theta}. Using
Proposition~\ref{fidis}
and (\ref{eq:321}), we get that
\[
\lim_{\lambda\to\infty}P[{\mathcal A}_{\lambda,\lambda
^{-2/3}h}\ge
\pi/2] = \exp \biggl( -\frac{ 4 \sqrt{2} } {3} h^{2/3}  \biggr).
\]
It
remains to remark that for every $\alpha<\pi/2$, $\lim_{\lambda\to
\infty} P[{\mathcal A}_{\lambda,\lambda^{-2/3}h}\ge\alpha]=1.$
Indeed, a
direct expansion in (\ref{eq:2}) shows that
\[
\ell({\cal
S}_{\alpha,\lambda^{-2/3}h})\mathop{\sim}\limits_{\lambda\to\infty} \biggl(\sin(\alpha)\cos(\alpha)+2\frac{\sin
^3(\alpha)}
{\cos(\alpha)}-\frac{\sin^3(\alpha)}{2\cos^3(\alpha)}
\biggr)\lambda^{-4/3}h^2.
\]
Inserting this estimation in (\ref{eq:322}) completes the proof.
\end{pf}

  The next lemma provides the explicit distribution of
$r(u,{\cal
P}_{\lambda})$ in terms of~${\mathcal A}_{\lambda,h}$.
%
\begin{lemm}\label{explicitr}
For all $h\in[0,1]$ and $u \in\S^1$,
%
\begin{eqnarray}\label{eq:3}
&&P[r(u,{\mathcal P}_{\lambda})\ge
h]\nonumber\\
&& \qquad =P[s(u,{\cal P}_{\lambda})\ge
h]
\\
&& \qquad  \quad {}+\lambda\int_0^{\pi/2}\frac{\partial\ell
({\cal S}_{\alpha,h})}{\partial\alpha} \exp\bigl\{-\lambda\ell\bigl({\emph{cap}}_1
\bigl[u,\bigl(1-(1-h)\sin(\alpha)\bigr)\bigr]\bigr)\bigr\}\,d\alpha,
\nonumber
\end{eqnarray}
where $\ell({\operatorname{cap}}_1[u,(1-(1-h)\sin({\alpha}))])$ and
$\ell({\cal S}_{\alpha,h})$ are defined at (\ref{eq:1}) and
(\ref{eq:2}), respectively.
\end{lemm}
\begin{pf} For fixed $h\in[0,1]$ and $\alpha\in[0,\pi/2)$, we
define the set (which is hatched in Figure~\ref{picture})
\[
{\mathcal F}_{h,\alpha}:={\operatorname{cap}}_1\bigl[\operatorname{rot}_{\alpha-\pi/2}(u)
,\bigl(1-(1-h)\sin(\alpha)\bigr)\bigr]\setminus{\cal
S}_{\alpha,h},
\]
where $\operatorname{rot}_{\theta}$ is the classical rotation of angle
$\theta\in[0,2\pi)$ defined on $\S^1$.

We remark that $x$ is outside the convex hull if and only if either
${\mathcal A}_{\lambda,h}$ is bigger than $\pi/2$, or ${\cal
F}_{h,\alpha}$ is empty. Consequently, we have for $u\in\S^1$
\[
P[r(u,{\mathcal P}_{\lambda})\ge
h] =P[{\mathcal A}_{\lambda,h}\ge
\pi/2]+\int_0^{\pi/2}\exp\{-\lambda\ell({\mathcal
F}_{h,\alpha})\}\,dP_{{\mathcal A}_{\lambda,h}}(\alpha),
\]
where
$dP_{X}$ denotes the distribution of $X$. Applying Lemma
\ref{theta} yields the result.
\end{pf}%

 The next proposition
provides the asymptotic behavior of the distribution of
$\widehat{r}_{\lambda}(\cdot)$:
%
\begin{prop}\label{limitr}
We have for all $h\ge0$ and $u \in\S^1$,
\begin{eqnarray*}
\lim_{\lambda\to\infty}P[\lambda^{2/3}r(u,{\mathcal P}_{\lambda
})\ge
h] &=&\exp\biggl \{-\frac{4\sqrt{2}h^{3/2}}{3} \biggr\}\\
&&{}+2\int
_0^{\infty}
\exp \biggl\{-\frac{4\sqrt{2}}{3}\biggl(h+\frac{t^2}{2}\biggr)^{3/2} \biggr\}t^2\,dt
-1.
\end{eqnarray*}
\end{prop}
\begin{pf} We focus on the asymptotic behavior of the integral in
the relation (\ref{eq:3}) where $h$ is replaced with
$\lambda^{-2/3}h$. We proceed with the change of variable
$\alpha=\frac{\pi}{2}-\lambda^{-1/3}t$, which gives
%
\begin{eqnarray}\label{eq:int}
&&\lambda\int_0^{\pi/2}\frac{\partial\ell
({\cal S}_{\alpha,h})}{\partial\alpha}(\alpha,\lambda^{-2/3}h)
\exp\bigl\{-\lambda\ell\bigl({\operatorname{cap}}_1\bigl[u,\bigl(1-(1-\lambda^{-2/3}h)\sin
(\alpha)\bigr)\bigr]\bigr)\bigr\}\,d\alpha\nonumber\hspace*{-15pt}\\
&& \qquad =\lambda^{2/3}\int_0^{\fracd{\pi}{2}\lambda^{1/3}}\frac{\partial
\ell
({\cal
S}_{\alpha,h})}{\partial
\alpha}\biggl(\frac{\pi}{2}-\lambda^{-1/3}t,\lambda^{-2/3}\biggr)\hspace*{-15pt}\\
&&\hphantom{=\lambda^{2/3}\int_0^{\fracd{\pi}{2}\lambda^{1/3}}}
 \qquad{}\times \exp\bigl\{-\lambda\ell\bigl({\operatorname{cap}}_1\bigl[u,\bigl(1-(1-\lambda
^{-2/3}h)\cos(\lambda^{-1/3}t)\bigr)\bigr]\bigr)\bigr\}\,dt.\nonumber\hspace*{-15pt}
\end{eqnarray}
Using (\ref{eq:1}), we find the exponential part of the integrand,
which yields
%
\begin{eqnarray}\label{eq:exp}\quad
&&\lim_{\lambda\to
\infty}\exp\biggl\{-\lambda\ell\biggl({\operatorname{cap}}_1\biggl[u,\biggl(1-(1-\lambda^{-2/3}h)
\sin\biggl(\frac{\pi}{2}-\lambda^{-1/3}t\biggr)\biggr)\biggr]\biggr)\biggr\}\nonumber
\\[-8pt]
\\[-8pt]
&& \qquad =\exp \biggl\{-\frac
{4\sqrt{2}}{3} \biggl(h+\frac{t^2}{2} \biggr)^{3/2} \biggr\}.
\nonumber
\end{eqnarray}
Moreover, the derivative of the area of ${\mathcal S}_{\alpha,h}$ is
\[
\frac{\partial\ell
({\cal
S}_{\alpha,h})}{\partial\alpha}=1+(1-h)^2\cos(2\alpha)-2(1-h)\cos
(\alpha)\sqrt{1-(1-h)^2\sin^2(\alpha)}.
\]
In particular, we have
%
\begin{equation}\label{eq:partial}\quad
\frac{\partial\ell
({\cal
S}_{\alpha,h})}{\partial\alpha}\biggl(\frac{\pi}{2}-\lambda
^{-1/3}t,\lambda^{-2/3}h\biggr)\mathop{\sim}\limits_{\lambda\to\infty}2\lambda^{-2/3}
 \bigl[h+t^2-t\sqrt{2h+t^2} \bigr].
\end{equation}
Inserting (\ref{eq:exp}) and (\ref{eq:partial}) into (\ref{eq:int})
and using (\ref{eq:3}), we obtain the required result.
\end{pf}%

\begin{figure}

\includegraphics{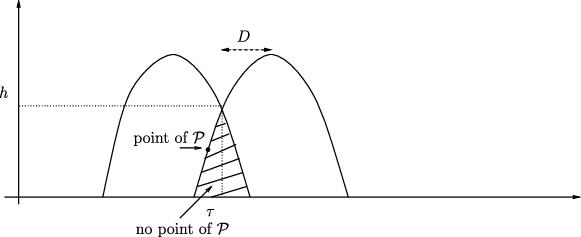}

\caption{Definition of the r.v. D.}
\label{picture2}
\end{figure}

\begin{rema}\label{rem2} In connection with Section
\ref{g+hullprocess}, the above calculation could have been
alternatively based on the limiting hull process related to
$\widehat{r}$. Indeed, for fixed $v\in{\mathbb R}, h\in{\mathbb
R}_+$, saying
that $\partial\Psi(v)$ is greater than $h$ means that there
is no translate of the standard downward parabola
$\Pi^{\downarrow}$ containing two extreme points on its boundary
and
lying underneath the point $(v,h)$. 
We
define a random variable $D$ related to the point $(v,h)$; see
Figure~\ref{picture2}. If ${\mathcal P}\cap((v,h) \oplus
\Pi^{\downarrow})$ is empty, then we take $D=0$. Otherwise, we
consider all the translates of $\Pi^{\downarrow}$ containing on the
boundary at least one point from ${\cal P}\cap((v,h)\oplus
\Pi^{\downarrow})$ and the point $(v,h)$. There is almost
surely precisely one among them which has the farthest peak (with
respect to the first coordinate) from $(v,h)$. The random
variable $D$ is then defined as the difference between the
$v\mbox{-}\mathit{coordinate}$ of the farthest peak and $v$. The
distribution of $|D|$ can be made explicit:
\[
P[|D|\le
t]=\exp\bigl \{-\tfrac{2}{3}(2h+t^2)^{3/2}+t\bigl(2h+\tfrac{2}{3}t^2\bigr)
\bigr\},\qquad
t\ge0.
\]
Conditionally on $|D|$, $\partial\Psi(v)$ is
greater than $h$ if and only if the region between the $v$-axis and the
parabola with the farthest peak does not contain any
point of ${\cal P}$ in its interior. Consequently, we have
\begin{eqnarray*}
 P[\partial\Psi(v)\ge h]&=&P[D=0]
\\&&{} +\int_0^{\infty}\exp \biggl\{ \biggl(-\frac{4\sqrt{2}}{3}
\biggl(h+\frac{t^2}{2}\biggr)^{3/2}\\
&&\hphantom{{} +\int_0^{\infty}\exp \biggl\{ \biggl(}{}-\frac{2}{3}(2h+t^2)^{3/2}-t\biggl(2h+\frac
{2}{3}t^2\biggr) \biggr) \biggr\}\,dP_{|D|}(t),
\end{eqnarray*}
%
which provides the result of Proposition~\ref{limitr}.
\end{rema}

%

The final proposition is the analog of Proposition~\ref{fidis} where the
radius-vector function of the flower is replaced by the one of the
convex hull itself.

\begin{figure}

\includegraphics{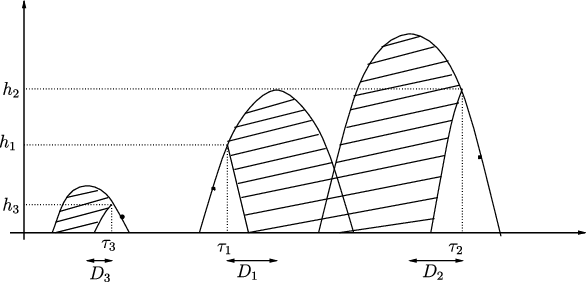}

\caption{Definition of the area $F$ (hatched region). The black points belong
to ${\mathcal P}$.}
\label{picture3}
\end{figure}

\begin{prop}
Let $n\in\N$, $0\le\theta_1<\theta_2\cdots<\theta_n<2\pi$
and $h_i \in(0,\infty)$ for all $i =1,\ldots,n$. 
Then
\begin{eqnarray*}
&&P\Biggl[\lambda^{2/3}r(\exp_1(\theta_1),{\mathcal P}_{\lambda})\ge
h_1;\ldots; \lambda^{2/3}
r(\exp_1(\theta_n),{\mathcal P}_{\lambda})\ge h_n \Biggr]\\
&& \qquad  \mathop{\sim}\limits_{\lambda\to\infty}\prod_{i=1}^nP[\lambda^{2/3}r(\exp_1(\theta
_i),{\mathcal P}_{\lambda})\ge h_i].
\end{eqnarray*}
Moreover, for every $v_1<v_2<\cdots<v_n\in\R$, we have
\begin{eqnarray*}
&& \lim_{\la\to\infty}
P[\lambda^{2/3}r(\exp_1(\la^{-1/3}v_1),{\mathcal P}_{\lambda})\ge
h_1;\ldots; \lambda^{2/3}
r(\exp_1(\la^{-1/3}v_n),{\mathcal P}_{\lambda})\ge h_n]\\
&& \qquad  =
\int_{\R^n}\exp \{-F((t_i,h_i,v_i)_{1\le i\le n})
\}\,dP_{{(D_1,\ldots,D_n)}}(t_1,\ldots,t_n),
\end{eqnarray*}
where $D_1,\ldots,D_n$ are symmetric variables such that
%
\begin{eqnarray}\label{distrD}
&&P[|D_1|\le
t_1;\ldots;|D_n|\le t_n]\nonumber
\\[-8pt]
\\[-8pt]
&& \qquad =\exp \biggl(-\int\sup_{1\le i\le
n}\biggl[\biggl(h_i+\frac{t_i^2}{2}-\frac{(|v-v_i|+t_i)^2}{2}\biggr)\vee0\biggr]\,dv \biggr),
\nonumber
\end{eqnarray}
and $F$ is the area
%
\begin{eqnarray}
\label{funcF}
&&F((t_i,h_i,v_i)_{1\le i\le n})=\int_{\R} \biggl\{\sup_{1\le i\le
n} \biggl[ \biggl(h_i+\frac{t_i^2}{2}-\frac{(v-v_i-t_i)^2}{2}
\biggr)\vee0 \biggr]
\nonumber
\\[-8pt]
\\[-8pt]&&\hphantom{F((t_i,h_i,v_i)_{1\le i\le n})=\int_{\R} \biggl\{}  {}-\sup_{1\le i\le
n} \biggl[ \biggl(h_i+\frac{t_i^2}{2}-\frac{(|v-v_i|+t_i)^2}{2}
\biggr)\vee
0 \biggr] \biggr\}\,dv.
\nonumber
\end{eqnarray}
%
\end{prop}
\begin{pf} We prove the first assertion and denote by
${\mathcal A}_1,\ldots,{\mathcal A}_n$ the angles (as defined by Lemma
\ref{theta}) corresponding to the couples
$(\theta_1,\lambda^{-2/3}h_1),\ldots,(\theta_n,\allowbreak \lambda^{-2/3}h_n)$.
Conditionally on $\{{\mathcal A}_i=\alpha_i\}$, the event
$\{\lambda^{2/3}r(\exp_1(\theta_i),{\mathcal
P}_{\lambda})\ge h_i\}$ only involves the points of the point process
${\mathcal P}_{\lambda}$ included in the circular cap
$\operatorname{cap}_{1}[\theta_i-\frac{\pi}{2}+\alpha_i,(1-(1-\lambda
^{-2/3}h_i)\sin(\alpha_i))]$; see the proof of Lemma \ref
{explicitr}. Moreover there exists
$\delta\in(0,\pi/2)$ such that for $\lambda$ large enough and
$\alpha_i\in(\delta,\frac{\pi}{2})$ for every $i$, these circular
caps are all disjoint. Consequently, we obtain that, conditionally on
$\{{\mathcal A}_i>\delta\ \forall i\}$, the events
$\{\lambda^{2/3}r(\exp_1(\theta_i),{\mathcal
P}_{\lambda})\ge h_i\}$ are independent.
It remains to remark that Lemma
\ref{theta} implies
\[
\lim_{\lambda\to\infty}P[\exists1\le i\le n;  {\mathcal A}_i\le
\delta]=0.
\]

Let us consider now the second assertion, which could be obtained by
a direct estimation of the joint distribution of the angles
${\mathcal A}_i$ [corresponding to the points
$(\lambda^{-1/3}v_i,\lambda^{-2/3}h_i)$]. But it is easier to
prove it with\vadjust{\goodbreak} the use of the boundary $\partial\Psi$ of the hull
process. As in Remark~\ref{rem2}, we define for each point $(v_i,h_i)$,
the random variable $D_i$ as the difference between the
$v$-coordinate of the farthest peak of a downward parabola
arising as a translate of $\Pi^{\downarrow}$ 
(denoted by $\operatorname{Par}_i$) containing on its boundary $(v_i,h_i)$
and a point of ${\cal P}$. Then
$|D_i|$ is less than $t_i$ for every $1\le i\le n$ if and only if there
is no point of ${\cal P}$ inside a region
delimited by the $v$-axis and the supremum of $n$ functions
$g_1,\ldots,g_n$ defined in the following way: $g_i(v_i+\cdot)$ is an even
function with a support equal to $[t_i-\sqrt{2h_i+t_i^2},t_i+\sqrt
{2h_i+t_i^2}]$ and identified with the parabola
$\operatorname{Par}_i(\cdot-v_i)$ on the segment $[t_i-\sqrt
{2h_i+t_i^2},0]$; see
Figure~\ref{picture2}. We
deduce from this assertion the result (\ref{distrD}). Conditionally on
$\{D_1=t_1,\ldots,D_n=t_n\}$, $\partial\Psi(v_i)$ is greater than $h_i$
for every $i$ if and only if the region between the functions $g_i$ and
the parabolae
$\operatorname{Par}_i$ does not contain any point of ${\cal P}$; see Figure
\ref{picture3}. This implies
result (\ref{funcF}) and completes the proof.\looseness=-1
\end{pf}

\begin{rema}\label{rem3} Convergence of the fidis
of the radius-vector function of
the convex hull of $n$ uniform points in the disk has already been
derived in Theorem 2.3 of \cite{hsing}. Still we feel that the
results presented in this section are obtained in a more direct and
explicit way. Moreover, they are characterized by the parabolic
growth and hull processes, which provides an elementary
representation of the asymptotic distribution.
The explicit fidis and the convergence of these fidis to those of
$\partial\Psi$ and $\partial\Phi$ can be used to obtain explicit
formulae for second-order characteristics of the point process of
extremal
points.
\end{rema}

\section{Stabilizing functional representation for convex hull characteristics}\label{sec6}\label{StabRepr}

The purpose of this section is to link the convex hull characteristics
considered in Section~\ref{INTRO} with the theory of stabilizing
functionals, a
tool for proving limit theorems in geometric contexts; see \cite{BY2,PeEJP,PeBer,PY1,PY4,PY5} and \cite{Sch}.


\textit{The collection $\Xi$ of basic geometric functionals.} We let
$\Xi$ be the collection of four basic functionals $\{\xi_s, \xi_r,
\xi_{\vartheta_k}, \xi_{f_k}\}$, where each $\xi_{\cdot}$ is defined
on pairs $(x, \X)$, with $x \in\X\subset\B^d$, according to the
following definitions. When $x \notin\X,$ we write
$\xi(x, \X)$ instead of $\xi(x, \X\cup\{x\})$.

The point-configuration functional $\xi_s(x,\X),  x \in\X
\subset\B^d,$
for finite $\X\subset\B^d$ is set to be zero if $x$ is not a vertex of
$\conv(\X)$, and otherwise it is defined as follows. Let ${\cal F}(x,
\X)$ be
the (possibly empty) collection of faces $f$ in $\Faces_{d-1}(\conv
(\X))$ such that
$x = \Top(f),$ where we recall from (\ref{FaceEmpM}) that $\Top(f)$ is
the point in $f$ which is closest to $\S^{d-1}$.
Let $\cone({\cal F}(x, \X)) := \{ r y,  r > 0,  y \in{\cal F}(x,
\X) \}$
be the corresponding cone.
Recalling that $F(\cdot)$ is the Voronoi flower defined at
(\ref{flower}), for $x \notin\ext(\X)$, we put $\xi_s(x,\X) = 0$,
and for $x \in\ext(\X)$, we put
\[
\xi_s(x,\X):=\Vol\bigl([\B^d \setminus F(\X)] \cap\cone({\cal F}(x,
\X))\bigr).
\]
Then the volume of $\B^d\setminus F({\mathcal P}_{\la})$ equals
$\sum_{x\in\P_{\la}}\xi_s(x,\P_{\la})$. Also, we know from
\eqref{ajout} that $W_{\la}$ is asymptotically equivalent to
the volume of $\B^d\setminus F({\mathcal P}_{\la})$. 

Likewise, given $x \in\X\subset\B^d,$ for
$x \notin\ext(\X)$, we put $\xi_r(x,\X) = 0$, and otherwise we put
\[
\xi_r(x,\X)= \Vol\bigl([\B^d \setminus\conv(\X)] \cap\cone({\cal F}(x,
\X))\bigr).
\]
Thus the volume of $\B^d\setminus K_{\la}$ equals $\sum_{x \in
\P_{\la}} \xi_r(x,\P_{\la})$, and we note that $V_{\la}$ is
asymptotically equivalent to the volume of $\B^d\setminus K_{\la}$.
%

 The $k$th order projection avoidance functional $\xi_{\vartheta
_k}(x,\X),  x \in\X\subset\B^d,$
with $k \in\{1,\ldots,d\}$ is zero if $x \notin\ext(\X),$ and
otherwise equal to
\[
\xi_{\vartheta_k}(x,\X) := \int_{[\B^d \setminus\conv(\X)] \cap
\cone({\cal F}(x, \X))}
\frac{1}{|x|^{d-k}} \vartheta^{\X}_k(x)\,dx;
\]
see (\ref{PAVOF}). In particular, (\ref{KUBOTAREPR}) yields %
%
\begin{equation}\label{INTVOLSTABREPR}
V_k(\B^d) - V_k(K_{\la}) = \frac{ (^{d-1}_{k-1} )}{\kappa_{d-k}}
\biggl [ \sum_{x \in\P_{\la}} \xi_{\vartheta_k}(x,\P_{\la
}) \biggr].
\end{equation}

The $k$-face functional $\xi_{f_k}(x; \X)$, defined for finite
$\X$ in $\B^d$, $x \in\X$ and $k \in\{0,1,\ldots,d - 1\}$, is
the number of $k$-dimensional faces $f$ of $\conv(\X)$
with $x = \Top(f),$ if $x$ belongs to $\Vertices[\conv(\X)],$
and zero otherwise. Thus $\sum_{x \in\X} \xi_{f_k}(x, \X)$
is the total number of $k$-faces in $\conv(\X)$.
In particular, setting $\X:= \P_{\la},$ 
the total mass of $\mu^{f_k}_{\la}$ is
%
\begin{equation}\label{FKREPR}
f_{k}(K_\la) = \sum_{x \in\P_{\la}} \xi_{f_k}(x, \P_{\la}).
\end{equation}
It is readily seen by the definition
of $\mu^{f_k}_{\la}$ at (\ref{FaceEmpM}) that 
%
\begin{equation}\label{FKEMAPPROX}
\mu^{f_k}_{\la} := \mu^{\xi_{f_k}}_{\la} := \sum_{x \in\P_{\la
}} \xi_{f_k}(x,\P_{\la}) \delta_x.
\end{equation}

\textit{The collection $\Xi^{(\infty)}$ of scaling
counterparts to elements of $\Xi$}. In the spirit of the local
scaling Section
\ref{LocScal}, we shall construct \textit{scaling counterparts} to each
functional $\xi\in\Xi$; we shall define these counterparts
in terms of the paraboloid growth and hull processes. To reflect this
correspondence we write $\xi^{(\infty)}_{\cdot}$ to denote the local
scaling limit analog
of $\xi_{\cdot}$ with the $(\infty)$ superscript.

The functional $\xi^{(\infty)}_s(x,\P)$ is defined to be zero
if $x\notin\ext(\Psi)$
and otherwise is defined as follows.
Let $\F^{\infty}(x, \P)$ stand for the collection
of paraboloid faces $f$ of $\Phi$ for which $x=\Top(f)$ [recall (\ref
{FaceEmpM})] and let $\vcone(\F^{\infty}(x, \P))$
be the cylinder (vertical cone) in $\R^{d-1} \times\R_+$ generated
by $\F^{\infty}(x, \P),$
that is to
say, $\vcone(\F^{\infty}(x, \P)) := \{ (v,h),  \exists {h'}:
(v,h') \in\F^{\infty}(x, \P) \}.$
Then, if $x \in\ext(\Psi),$ we set
\[
\xi_s^{(\infty)}(x,\P) := \Vol\bigl(\vcone(\F^{\infty}(x, \P))
\setminus
\Psi\bigr).
\]
Formally we should define $\xi^{(\infty)}_s(x;\X)$ for
general $\X$ rather than just for $\P,$ but we bypass this formality
so as to avoid extra
notation. We will mainly consider $\X= \P$ anyway and the general
definition can
be readily recovered by formally conditioning on $\P= \X.$ This
simplifying convention
will also be applied for the remaining local scaling functionals below.

Likewise, $\xi^{(\infty)}_r(x,\P)$ is zero if
$x\notin\ext(\Psi)$, and otherwise
\[
\xi^{(\infty)}_r(x,\P) :=
\Vol\bigl(\vcone(\F^{\infty}(x, \P)) \setminus\Phi\bigr).
\]

The $k$th order projection avoidance functional $\xi^{(\infty
)}_{\vartheta_k}(x,\P)$ is
zero if $x\notin\ext(\Psi)$, and otherwise
%
\begin{equation}\label{XIVARTH}
\xi^{(\infty)}_{\vartheta_k}(x, \P) := \int_{\vcone(\F^{\infty
}(x,\P)) \setminus\Phi} \vartheta^{\infty}_k(u)\,du
\end{equation}
with $\vartheta^{\infty}_k(\cdot)$ defined in (\ref{PARAVOF}).
Note that the extra factor $\frac{1}{|x|^{d-k}}$ in (\ref
{KUBOTAREPR}), where $x \in\B^d \setminus K_\la$,
converges to one under
the scaling $T^{\la}$ defined at (\ref{SCTRANSF}), and thus is not
present in the asymptotic functional.

The $k$-face functional $\xi^{(\infty)}_{f_k}(x, \P)$, defined for
$x \in\P$, and $k \in\{0, 1,\ldots, \allowbreak d - 1\}$, is the number of
$k$-dimensional paraboloid
faces $f$ of the hull process $\Phi$ for which $x = \Top(f),$ if $x$
belongs to $\ext(\Psi) = \Vertices(\Phi)$, and zero otherwise.

\textit{The collection $\Xi^{(\la)}$ of finite-size scaling
counterparts to elements of $\Xi$.} For each of the four basic
functionals $\xi\in\Xi$, and each $\la\geq1$, consider the
collection $\Xi^{(\la)}$ of the finite-size scaling counterparts
$\xi^{(\la)}$ given by
%
\begin{equation}\label{FUNCTRANSF}
 \qquad \xi^{(\la)}(x,\X) := \xi([T^{\la}]^{-1} x,[T^{\la}]^{-1} \X), \qquad
x \in\X\subset{\cal R}_{\la}
\subset\R^{d-1} \times\R_+,
\end{equation}
where $T^{\la}$ is the scaling transform (\ref{SCTRANSF}) and ${\cal
R}_{\la}$ its image (\ref{RLA}).
Again resorting to the theory developed in Section~\ref{LocScal}, we
see for $\xi\in\Xi$ that $\xi^{(\la)}$ can be regarded as \textit{interpolating} between
$\xi$ and $\xi^{(\infty)}$; as such it is the analog of $\xi^{(\la
)}$ defined at
(1.17) of \cite{SY}. However, due to the differing natures of the functionals
considered here, different scaling pre-factors are needed to ensure nontrivial
scaling behaviors. More precisely, for each $\xi^{(\la)} \in\Xi
^{(\la)}$ we
define its \textit{proper scaling prefactor} $\la^{\eta[\xi]}$
where:
\begin{itemize}
\item$\eta[\xi_s] = \eta[\xi_r] = \eta[\xi_{\vartheta_k}] =
\beta(d-1) + \gamma,  k \in\{0, 1,\ldots, d -
1\},$ since for each of these three functionals, the spatial scaling
involves dilation by $\la^{\beta}$,
whereas the time scaling involves $\la^{\gamma}$.
[Note that the
re-scaled projection avoidance function (\ref{PAVOFRESC}) involves no
scaling prefactor.]

\item$\eta[\xi_{f_k}] = 0$ because the number of $k$-faces does not
undergo any scaling.
\end{itemize}

To proceed, for any measurable $D \subseteq\R^{d-1} \times\R_+$
and generic scaling limit functional $\xi^{(\infty)} \in
\Xi^{(\infty)}$, by its \textit{restricted version} we mean by
$\xi^{(\infty)}_D(x, \P):=
\xi^{(\infty)}(x, \P\cap D)$, $x \in\R^{d-1} \times\R_+$. %
Note that the so-defined restricted functionals in case of $D$
bounded, or of bounded spatial extent, clearly involve growth and
hull processes built on input of bounded spatial extent, in which
case the definition (\ref{PHIDEF}) for $\P$ replaced with $\P\cap
D$ yields infinite vertical faces at the boundary of $D's$ spatial
extent. This makes some of the functionals considered in this paper
infinite or even undefined for points close to these infinite faces.
For such points $x, x := (v_x, h_x)$, and such sets $D$, we may
formally put $\xi_D^{(\infty)} = \infty$. Fortunately, such
pathologies do not arise in the sequel. Indeed, we will restrict
to cylinder sets $D$ centered around the vertical axis $\{(v_x, h),
 h \geq0 \}$ whose radius (termed stabilization radius below) is
sufficiently large so that with probability one the faces meeting
$x$, as defined by the input $\P\cap D$, coincide with the faces
meeting $x$ when the input is $\P$. We now make these ideas
precise.

Having now defined the class $\Xi$ of four basic
functionals, together with the finite-size scaling version
$\Xi^{(\la)}, \la\geq1,$ and the infinite scaling version
$\Xi^{(\infty)}$, we are ready to establish some crucial
localization properties of the functionals in $\Xi^{(\la)}$ and
$\Xi^{(\infty)}$. Recalling that $B_{d-1}(v,r)$ is the $(d-1)$
dimensional ball
centered at $v \in\R^{d-1}$ with radius $r$, let
$C_{d-1}(v,r)$ be the cylinder $B_{d-1}(v,r) \times\R_+$.
Given a generic scaling limit functional $\xi^{(\infty)} \in
\Xi^{(\infty)}$, we shall
write $\xi^{(\infty)}_{[r]} := \xi^{(\infty)}_{C_{d-1}(v,r)}.$
Likewise, for the
finite scaling functionals $\xi^{(\la)} \in\Xi^{(\la)},$ we shall
use the notation
$\xi^{(\la)}_{[r]}$ with a fully analogous meaning.

Given $\xi^{(\infty)} \in
\Xi^{(\infty)}$, a random variable $R := R^{\xi^{(\infty)}}[x]$ is
called a
\textit{localization radius} for $\xi^{(\infty)}$ if and only if a.s.
\[
\xi^{(\infty)}(x,\P) = \xi^{(\infty)}_{[r]}(x,\P)  \qquad \mbox{for all
} r \geq R.
\]
%
Given $\xi^{(\la)} \in\Xi^{(\la)},$
we analogously define $R^{\xi^{(\la)}}[x]$ to be a localization
radius for
$\xi^{(\la)}$ if and only if a.s.
\[
\xi^{(\la)}\bigl (x, \P^{(\la)} \bigr) = \xi^{(\la)}_{[r]} \bigl(x, \P^{(\la)}
\bigr)  \qquad \mbox{for all } r \geq R.
\]

The notion of localization, developed in \cite{SY},
is a variant of a general concept of stabilization \cite{BY2,PeEJP,PY4,PY5}.
A crucial property of the functionals $\xi^{(\la)} \in\Xi^{(\la)},
\la\geq1$ and $\xi^{(\infty)} \in
\Xi^{(\infty)}$
is that they admit localization radii with tails decaying
super-exponentially fast.
%
\begin{lemm}\label{LocalLem} For each $\xi\in\Xi$, the functionals
$\xi^{(\infty)}$ and $\xi^{(\la)}, \la\geq1,$ admit
localization radii with the property that
%
\begin{eqnarray}\label{LocalExpr}
P\bigl[R^{\xi^{(\infty)}}[x] > L\bigr] &\leq& C \exp\biggl(- \frac{L^{d+1}}{ C}\biggr)
 \quad \mbox{and}  \nonumber
 \\[-8pt]
 \\[-8pt]   P\bigl[R^{\xi^{(\la)}}[x] > L\bigr] &\leq& C \exp\biggl(- \frac
{L^{d+1}} { C}\biggr)
\nonumber
\end{eqnarray}
for some finite positive constant $C,$ uniformly in $\la$ large
enough and uniformly in~$x$.
\end{lemm}

\begin{pf}
The proof is given only for the scaling limit functionals $\xi
^{(\infty)} \in
\Xi^{(\infty)}$;
the argument for the finite scaling functionals $\xi^{(\la)} \in\Xi
^{(\la)}, \la\geq1,$ is
fully analogous and is omitted.

For a point $x := (v,h) \in\P$, denote by $\P[[x]]$ the collection of
all vertices of $(d-1)$-dimensional faces of $\Phi$ meeting at $x$
if $x \in\Vertices(\Phi)$ and $\P[[x]] := \{ x \}$ otherwise.
If $x \in\Vertices(\Phi),$ the collection $\P[[x]]$ uniquely
determines the \textit{local facial structure} of $\Phi$ at $x,$
understood as the collection of all $(d-1)$-dimensional faces
$f_1[x],\ldots,f_m[x],  m = m[x] < \infty$ meeting at $x.$
We shall show that there exists a random variable $R':= R'[x]$ with
these two properties:

$\bullet$ With probability one the facial structure $\P_{[r]}[[x]]$ at
$x$ determined upon restricting to $\P:= \P\cap C_{\R^{d-1}}(v,r)$
coincides with $\P[[x]]$ for all $r \geq R'$; in the sequel
we say that $\P[[x]]$ is fully determined within radius $R'$
in such a case.

$\bullet$  We have
%
\begin{equation}\label{RDECAY}
P[R' > L] \leq C \exp\biggl(- \frac{L^{d+1}} { C}\biggr).
\end{equation}
Before proceeding, we note that to conclude the statement of Lemma~\ref
{LocalLem}, it is enough
to establish (\ref{RDECAY}). Indeed, this is because of these
three observations:

$\bullet$ The values of functionals $\xi^{(\infty)}_s,  \xi_r^*$
and $\xi^{(\infty)}_{f_k},  k \in\{0,\ldots,d-1\},$ at \mbox{$x \in\P
$}, are uniquely
determined given $\P[[x]]$, and thus $R'$ can be taken as the localization
radius.

$\bullet$ The values of functionals $\xi^{(\infty)}_{\vartheta_k}(x,\P
),  k \in\{1,\ldots,d-1\},   x:= (v,h),$
are determined, given the intersection of the hull process $\Phi$ with
$\Theta[x] :=\break [\vcone(\F^{(\infty)}(x, \allowbreak \P)) \setminus\Phi]
\oplus\Pi^{\downarrow};$ see
(\ref{PARAVOF}) and the definition of $\xi^{(\infty)}_{\vartheta
_k}$ at \eqref{XIVARTH}. It is readily seen
that this intersection $\Theta[x] \cap\Phi$ is in its turn uniquely
determined
by $\Theta[x] \cap\Vertices(\Phi)$.\vadjust{\goodbreak} Thus, to know it, it is enough
to know
the facial structure at $x$ and at all vertices of $\Phi$ falling into
$\Theta[x]$. To proceed, note that the spatial diameter of $\Theta
[x]$ is
certainly bounded by $R'[x]$ plus $2\sqrt{2}$ times the square root of
the highest
height coordinate of $\partial\Phi$ within spatial distance $R'[x]$
from $v.$
Use (\ref{SuperexpBd}) to bound this height coordinate and thus to establish
a superexponential bound $\exp(-\Omega(L^{d+1}))$ for tail
probabilities of
the spatial diameter $R''[x]$ of $\Theta[x].$ Finally, we set the localization
radius to be $\max_{y \in\Vertices(\Phi), y \in C_{\R
^{d-1}}(v,R''[x])} R'[y]$ which is
again easily verified to exhibit the desired tail behavior as the number
of vertices within $C_{\R^{d-1}}(v,R''[x])$ grows polynomially in
$R''[x]$ with
overwhelming probability; see Lemma 3.2 in \cite{SY}.

%
To proceed with the proof, suppose first that $x$ is not extreme in
$\Phi.$
Then, by Lemma 3.1 in \cite{SY} and its proof, there exists $R' =
R'[x]$ satisfying (\ref{RDECAY})
and such that the extremality status of $x$ localizes within radius $R'.$
In this particular case of $x$ not extreme in $\Phi$, this also
implies localization
for $\P[[x]] = \{ x \}.$ Assume now that $x$ is an extreme point in
$\P$.
Enumerate the $(d-1)$-dimensional faces meeting $x$ by $f_1,\ldots,f_m.$ The
local facial structure $\P[[x]]$ is determined by the parabolic faces
of the space--time region $\bigcup_{i \leq m} \Pi^{\downarrow}[f_i]$,
which by (\ref{PHIDEF2}) is devoid of points from $\P$. Note that this
region contains all vertices of $f_1,\ldots,f_m$ on its upper boundary.
Moreover, Poisson points outside this region do not change the status
of the faces $f_1,\ldots,f_m$ as these faces will not be subsumed by
larger faces meeting $x$ unless Poisson points lie on the boundary
of the hull process, an event of probability zero. It follows that
$\P[[x]]$ is fully determined by the point configuration $\P\cap
C_{d-1}(v, R')$ where $R'$ is the smallest integer $r$ such that
%
\begin{equation}\label{stabregion} \bigcup_{i \leq m} [\Pi
^{\downarrow}[f_i]
\cap(\R^{d-1} \times\R_+)] \subset C_{d-1}(v, r).
\end{equation}
To establish (\ref{RDECAY}) for $R'$, we note that if $R'$ exceeds $L$,
then, by standard geometry, within distance $O(L^2)$ from $x$, we can
find a point $x'$ in $\Z^d$ with the properties that:
%
\begin{itemize}
\item the downwards parabolic solid $x' \oplus\Pi^{\downarrow}$ is contained
in $\bigcup_{i \leq m} \Pi^{\downarrow}[f_i]$ and thus in particular
devoid of points of $\P;$
\item the spatial diameter (the diameter of spatial projection on $\R^{d-1}$)
of $[x' \oplus\Pi^{\downarrow}] \cap(\R^{d-1} \times\R_+)$ does
exceed $L/2.$
\end{itemize}
Since the intensity measure of $\P$ assigns to such $[x' \oplus\Pi
^{\downarrow}]
\cap(\R^{d-1} \times\R_+)$ mass of order at least $\Omega
(L^{d+1})$ [in fact
even $\Omega(L^{d+1+2\delta})$; see the proof of Lemma~3.1 in \cite
{SY} for
details in a much more general set-up], the probability of having
$x' \oplus\Pi^{\downarrow}$ devoid of points of $\P$ is $\exp
(-\Omega(L^{d+1})).$
Since the cardinality of $B_{d}(x, L^2) \cap\Z^{d-1}$ is bounded by $CL^{2d},$
Boole's inequality gives
\[
P[R' > L] \leq CL^{2d} \exp\biggl(- \frac{L^{d+1}} {C}\biggr),
\]
which yields the required inequality (\ref{RDECAY}) and thus
completes the proof of Lemma~\ref{LocalLem}.
\end{pf}

\section{Variance asymptotics and Gaussian limits for empirical measures}\label{sec7}\label{GAUSS} 

Sections~\ref{sec1}--\ref{sec5}
establish the asymptotic embedding in $\R^{d-1} \times\R_+$ of
convex polytope characteristics, whereas Section~\ref{sec6} establishes their
localization proper\-ties. The present section establishes variance
asymptotics and Gaussian limits of these characteristics by
exploiting this embedding within the \mbox{framework} of general methods of
stabilization theory for point processes on \mbox{$\R^{d-1} \times\R_+$}.

Given a generic functional $\xi\in\Xi$, recall from
(\ref{FUNCTRANSF}) its finite size scaling counterpart $\xi^{(\la)}
\in\Xi^{(\la)}$, namely
\[
\xi^{(\la)}(x,\X) := \xi([T^{\la}]^{-1} x,[T^{\la}]^{-1} \X), \qquad
x \in\X\subset{\cal R}_{\la}
\subset\R^{d-1} \times\R_+.
\]
Put
%
\begin{equation}\label{GENEMPM}
\mu^{\xi}_{\la} := \sum_{x \in\P^{(\la)}} \xi^{(\la)}\bigl(x,\P
^{(\la)}\bigr) \delta_x
\end{equation}
and $\bar{\mu}^{\xi}_{\la} := \mu^{\xi}_{\la} - \E\mu^{\xi
}_{\la}.$

As in Section~\ref{sec6}, we write $\xi^{(\infty)} \in\Xi^{(\infty)}$ to
denote the local scaling limit analog
of $\xi$; $\xi^{(\infty)}$ is defined on pairs $(x,
\X)$, with $x \in\X\subset\R^{d-1} \times\R_+$.
Recall that when $x \notin\X$, we write
$\xi(x, \X)$ instead of $\xi(x, \X\cup\{x\})$, with a similar
convention for $\xi^{(\la)}$.
Recall from Definition~\ref{fund-def} that $\P$ is a Poisson point
process in the upper half-space $\R^{d-1} \times\R_+$ with intensity
density $h^{\delta}\,dh\,dv$. Following \cite{SY}, we define
the second order correlation functions for $\xi^{(\infty)}$ given
by
%
\begin{equation}\label{SOScorr1}
\varsigma_{\xi^{(\infty)}}(x) := \E\bigl[\xi^{(\infty)}(x,\P)\bigr]^2,  \qquad  x
\in\R^{d-1} \times\R_+,
\end{equation}
whereas, for all $x, y \in\R^{d-1} \times\R_+$, we put
%
\begin{eqnarray}\label{SOcorr1}
\varsigma_{\xi^{(\infty)}}(x,y) &:=& \E\bigl[\xi^{(\infty)}(x, \P\cup
\{ y \})
\xi^{(\infty)}(y,\P\cup\{ x \})\bigr]\nonumber
\\[-8pt]
\\[-8pt] &&{}- \E\bigl[\xi^{(\infty)}(x,\P)\bigr]
\E\bigl[\xi^{(\infty)}(y,\P)\bigr].
\nonumber
\end{eqnarray}
Define also the asymptotic variance expression
%
\begin{eqnarray}\label{VARasympt1}
\sigma^2\bigl(\xi^{(\infty)}\bigr) &:=& \int_0^{\infty} \varsigma_{\xi
^{(\infty)}}((\0,h))h^{\delta}\,dh \nonumber
\\[-8pt]
\\[-8pt]&&{}+
\int_0^{\infty} \int_0^{\infty} \int_{\R^{d-1}}
\varsigma_{\xi^{(\infty)}}(({\bf0},h),(v',h')) h^{\delta}
h'^{\delta}\,dh\,dh'\,dv'.
\nonumber
\end{eqnarray}
These expressions are the counterparts to (1.7) and (1.8) in
\cite{SY}; recall that here we are working in the isotropic regime,
corresponding to $\rho_0 \equiv1$ in \cite{SY}.

Given $\xi\in\Xi$, consider the sum $\sum_{x \in\P^{(\la)}}
\la^{\eta[\xi]} \xi^{(\la)}(x, \P^{(\la)}).$ There are roughly
$\la^{\beta(d-1)}$ terms which do not vanish, and thus one expects
growth of order $\la^{\tau}$ with 
%
\begin{equation}\label{TAUDEF}
\tau= \beta(d-1)= \frac{d-1}{d+1+2\delta}.
\end{equation}
Upon centering and scaling by $\la^{-\tau/2}$, one may also expect
asymptotic normality as $\la\to\infty$. The following theorem, one
of the main results of this paper, makes this intuition precise. It
establishes a weak law of large numbers, variance asymptotics and a
central limit theorem for the afore-mentioned sums as well as for
$\la^{-\tau/2} \langle g, \la^{\eta[\xi]} \mu^{\xi}_{\la}
\rangle=
\la^{\zeta/2} \langle g, \mu^{\xi}_{\la} \rangle,  g \in{\cal
C}(\B^d)$, where we have $\zeta= -\tau+ 2 \eta$ from \eqref{ZETA}.

%
\begin{theo}\label{MainCLT} For all $\xi\in\Xi$ and all 
$g \in{\cal C}(\B^d)$, we have
%
\begin{equation}\label{EXPconv}
 \qquad \lim_{\la\to\infty} \la^{-\tau} \E\bigl[\bigl\langle g, \la^{\eta[\xi]}
\mu^{\xi}_{\la} \bigr\rangle\bigr]
= \int_0^{\infty} \E\bigl[ \xi^{(\infty)}(\0,h)\bigr] h^\delta \,dh \int_{\S
^{d-1}} g(u)\,d\sigma_{d-1}(u).
\end{equation}
The integral in (\ref{VARasympt1}) converges, and for all $g \in{\cal
C}(\B^d)$, we have
%
\begin{equation}\label{VARconv}
  \lim_{\la\to\infty} \la^{-\tau} \Var\bigl[\bigl\langle g, \la^{\eta[\xi
]} \overline{\mu}^{\,\xi}_{\la} \bigr\rangle\bigr]
= V^{\xi^{(\infty)}}[g] := \sigma^2\bigl(\xi^{(\infty)}\bigr) \int_{\S
^{d-1}} g^2(u)\,d\sigma_{d-1}(u).\hspace*{-35pt}
\end{equation}
%
Furthermore, the random variables $\la^{-\tau/2} \langle g, \la
^{\eta[\xi]}
\bar\mu^{\xi}_{\la} \rangle$ converge in law to $N(0,\allowbreak V^{\xi
^{(\infty)}}[g])$
as $\la\to\infty$.
Finally, if $\sigma^2(\xi^{(\infty)}) > 0,$ then for all $g \in
{\cal C}(\B^d)$
not identically zero, we have
%
\begin{eqnarray}\label{ConvRate}
&&\sup_t  \biggl| P \biggl[ \frac{\langle g, \la^{\eta[\xi]} \bar\mu
^{\xi}_{\la} \rangle}{\sqrt{\Var
[\langle g, \la^{\eta[\xi]} \bar\mu^{\xi}_{\la} \rangle]}} \leq
t  \biggr] - P[ N(0,1) \leq t]  \biggr|
\nonumber
\\[-8pt]
\\[-8pt]
&& \qquad = O (\la^{-\tau/2 }(\log\la)^{3d + 4 \delta+ 1} ).
\nonumber
\end{eqnarray}
\end{theo}

\begin{remarks*} (i) The expectation limit (\ref{EXPconv}) and
variance limit \eqref{VARconv} generalize the analogous limits
appearing at (2.2) and (2.3) in Theorem 2.1 of \cite{SY}, which is
restricted to the case that $\xi$ is the $k$-face functional with $k
= 0$. Likewise, convergence in law of $\la^{-\tau/2} \langle g,
\la^{\eta[\xi]} \bar\mu^{\xi}_{\la} \rangle$ and the rate result
(\ref{ConvRate})
extend the distributional results of Theorem 2.1 of \cite{SY}.


(ii) We refer to the statements \eqref{VARconv} and
\eqref{ConvRate} as \textit{measure-level variance asymptotics and
measure level central limit theorems} for $\la^{\eta[\xi]}
\mu^{\xi}_{\la},$
with scaling exponent $-\tau/2$ and with variance density $\sigma
^2(\xi^{(\infty)}).$
When $g \equiv1$, we obtain the limit theory for the total mass of
$\mu^{\xi}_{\la}$, giving \textit{scalar variance asymptotics} and
central limit theorems. For all $\xi\in\Xi$, Theorem
\ref{MainCLT} admits a multivariate version giving a central limit
theorem for the random vector $(\la^{-\tau/2} \langle g_1,
\la^{\eta[\xi]} \bar\mu^{\xi}_{\la} \rangle,\ldots,\la^{-\tau/2}
\langle g_m, \la^{\eta[\xi]} \bar\mu^{\xi}_{\la} \rangle)$, with
$g_i \in\C(\B^d)$ for all $i = 1,\ldots,m$, which follows from the
Cram\'er--Wold
device.\vadjust{\goodbreak} 

(iii) Given $\xi^{(\infty)}\in\Xi^{(\infty)}$, the question whether
$\sigma^2(\xi^{(\infty)})$ is strictly positive is
nontrivial, and the
application of general techniques of stabilization theory designed to check
this condition may be far from straightforward. These issues are
discussed at the end of this section.

(iv) We have not tried for optimal rates in (\ref{ConvRate}) and
expect that the exponents on the logarithm can be improved.
\end{remarks*}

  The proof of Theorem~\ref{MainCLT} depends on the
following three lemmas, which establish further properties of \textit{the scaling limit functionals $\xi^{(\infty)}\in\Xi^{(\infty)}$ and
local scaling functionals  $\xi^{(\la)} \in\Xi^{(\la)},   \la
\geq1.$}


\begin{lemm}\label{MOMBD}
For all $p > 0$ and all $\xi\in\Xi$, we have
%
\begin{eqnarray}\label{LIMITBD}
\sup_{x \in\R^{d-1}} \E\bigl[\bigl|\xi^{(\infty)}(x,\P)\bigr|^p\bigr] &<& \infty
 \quad \mbox{and} \nonumber
 \\[-8pt]
 \\[-8pt]
 \sup_{\la\geq1} \sup_{x \in{\cal R}_{\la}} \E\bigl[\bigl|\la^{ \eta[\xi
]} \xi^{(\la)}\bigl(x,\P^{(\la)}\bigr)\bigr|^p\bigr] &<&  \infty.
\nonumber
\end{eqnarray}
\end{lemm}

\begin{pf} We only give the proof for $\xi^{(\infty)},$ the
finite scaling case $\xi^{(\la)}$ being fully analogous. This is
done separately for all functionals considered.

For $\xi_s^{(\infty)}(x,\P)$ and $\xi_r^{(\infty)}(x,\P)$
we only consider the case
of $x$ extreme, for otherwise both functionals are zero. With
$x \in\Vertices(\Phi)$ we make use of (\ref{SuperexpBd}) to bound
the height and of (\ref{RDECAY}) and (\ref{stabregion}) to bound the spatial
size of the regions whose volumes define $\xi_s^{(\infty)}$ and $\xi
_r^{(\infty)}.$ Since these
bounds yield superexponential decay rates on each dimension separately, the
volume admits uniformly controllable moments of all orders. Finally, by
(\ref{XIVARTH}),
$0 \leq\xi^{(\infty)}_{\vartheta_k} \leq\xi^{(\infty)}_r$ whence
(\ref{LIMITBD}) follows
for $\xi_{\vartheta_k}$ as well.

For $\xi_{f_k}^{(\infty)}(x,\P)$, we only consider the case $x
\in\Vertices(\Phi)$,
and we let $N:= N[x]$ be the number of extreme points in $\P\cap
C_{d-1}(v, R'[x])$ with
$R'$ as in (\ref{stabregion}). Then $\xi^{(\infty)}_{f_k}(x, \P)$
is upper bounded by
${N \choose k - 1}$. By Lemma 3.2 of \cite{SY}, the probability
that a point $(v_1,h_1)$ is extreme in $\Phi$ falls off
superexponentially fast
in $h_1$; see again (\ref{SuperexpBd}). Consequently, in view of
(\ref{RDECAY}), the random variables
${N \choose k-1}$ and $\xi_{f_k}^{(\infty)}(x,\P)$ admit finite
moments of all
orders.
%

The proof of Lemma~\ref{MOMBD} is now complete.
\end{pf}

For all $h \in\R_+$, $(v',h') \in{\cal R}_\la,$ and $\xi\in\Xi$,
we put
\begin{eqnarray*}
c^{(\la)}((\0,h), (v',h')) &:=& \E\bigl[
\la^{\eta[\xi]}\xi^{(\la)}\bigl((\0,h),\P^{(\la)}
\cup(v',h')\bigr)
\\
&&\hphantom{\E\bigl[}{}\times
\la^{\eta[\xi]}\xi^{(\la)}\bigl((v',h'),\P^{(\la)} \cup(\0, h)\bigr) \bigr]
\\
&& {}-\E\bigl[\la^{\eta[\xi]}\xi^{(\la)}\bigl((\0,h),\P^{(\la)} \bigr)\bigr] \E
\bigl[\la^{\eta[\xi]}\xi^{(\la)}\bigl((v',h'),\P^{(\la)} \bigr) \bigr].
\end{eqnarray*}

The next lemma makes use of the moment bounds of Lemma~\ref{MOMBD}
and is proved through straightforward modifications of the proofs of
Lemmas 3.3 and 3.4 in \cite{SY}.

\begin{lemm} \label{lem3.4} For all $h \in\R_+$, $(v',h') \in{\cal
R}_\la$ and $\xi\in\Xi$ we
have as $\la\to\infty$,
\[
\E\bigl[\la^{\eta[\xi]}
\xi^{(\la)}\bigl((\0,h),\P^{(\la)}\bigr)\bigr] \to
\E\bigl[\xi^{(\infty)}((\0,h),\P)\bigr]
\]
and
\[  c^{(\la)}((\0,h), (v',h')) \to\varsigma_{\xi^{(\infty)}}((\0,h),
(v',h')).
\]
\end{lemm}

The next lemma is the analog of Lemma 3.5 in \cite{SY} and is proved
similarly.

\begin{lemm} \label{lem3.5} There is a constant $C < \infty$ such that
for all $h \in\R_+$, $(v',h') \in{\cal R}_\la$, and all $\xi\in
\Xi$, we have
\[
\bigl|c^{(\la)}((\0,h), (v',h'))\bigr| \leq C \exp \biggl( {-1 \over C} \max
(|v'|, h,
h')  \biggr)
\]
and
\[
\bigl| \varsigma_{\xi^{(\infty)}}((\0,h), (v',h'))\bigr| \leq C \exp\biggl ( {-1
\over C} \max(|v'|, h, h')  \biggr).
\]
\end{lemm}

Equipped with these lemmas, we now prove
Theorem~\ref{MainCLT}. We shall give separate proofs for (\ref{EXPconv}),
(\ref{VARconv}) and (\ref{ConvRate}), following closely the
methods of~\cite{SY}.\vspace*{12pt}

\textit{Proof of the expectation formula (\ref{EXPconv}).} For $g \in
{\cal
C}(\B^d),$ we have for all $\xi\in\Xi$,
%
\begin{equation}\label{EQEXPL1}
\E[\langle g, \mu^{\xi}_{\la} \rangle] = \lambda\int_{\B^d}
g(x) {\Bbb E}[\xi(x,\P_{\la})] (1-|x|)^{\delta}\,dx.
\end{equation}
%

By rotation invariance, we have that $\xi(x, \P_\la) \eqd
\xi(x^{\theta}, \P_\la^{\theta})$, where $x^{\theta}$ is $x$ rotated
by the angle $\theta$, and similarly for $\P_\la^{\theta}$. Letting
$\theta:=\theta_x$ be the rotation sending $x/|x|$ to $u_0$, gives
$\E\xi(x, \P_\la)= \E\xi^{(\la)}((\0,h), \P^{(\la)})$, where $h:=
\la^{\gamma}(1 - |x|)$.
Thus we rewrite (\ref{EQEXPL1}) as
%
\begin{eqnarray}\label{new1}
\E[\langle g, \mu^{\xi}_{\la} \rangle] &=& \lambda\int_{\B^d}
g(x) {\Bbb E} \bigl[\xi^{(\la)} \bigl(
(\0,h),\P^{(\la)}  \bigr) \bigr] \la^{-\gamma\delta} h^{\delta}\,dx
\nonumber
\\
&=& \lambda^{1-\gamma\delta} \int_{\S^{d-1}} \int_0^{\la^{\gamma
}} g\bigl(u (1-\la^{-\gamma}h)\bigr)
{\Bbb E} \bigl[\xi^{(\la)}  \bigl( (\0,h),\P^{(\la)}  \bigr)
\bigr]\\
&&\hphantom{\lambda^{1-\gamma\delta} \int_{\S^{d-1}} \int_0^{\la^{\gamma
}}}{}\times h^{\delta}
(1-\la^{-\gamma} h)^{d-1} \la^{-\gamma}\,dh\,d\sigma_{d-1}(u).
\nonumber
\end{eqnarray}
Noting that $\tau= 1 -\delta\gamma- \gamma$ and multiplying through
by $\la^{-\tau+ \eta[\xi]}$, we obtain
%
\begin{eqnarray}\label{EQEXPL3}
\la^{-\tau+ \eta[\xi]} \E[\langle g, \mu^{\xi}_{\la} \rangle]
&=& \int_{\S^{d-1}}
\int_0^{\la^{\gamma}} g\bigl(u (1-\la^{-\gamma}h)\bigr)
{\Bbb E} \bigl[\la^{\eta[\xi]} \xi^{(\la)}  \bigl( (\0,h),\P
^{(\la)}  \bigr) \bigr]\nonumber\hspace*{-35pt}
\\[-8pt]
\\[-8pt]
&&\hphantom{\int_{\S^{d-1}}
\int_0^{\la^{\gamma}}}{}\times
(1-\la^{-\gamma} h)^{d-1} h^{\delta}\,dh\,d\sigma_{d-1}(u).
\nonumber\hspace*{-35pt}
\end{eqnarray}

Notice that $\E[ \la^{\eta[\xi]} \xi^{(\la)}((\0,h),\P^{(\la)}]
h^\delta$ is dominated
by an integrable function of $h$, as the contribution coming from
large $h$ is well controllable as in Lemma~3.2 in \cite{SY}---in
particular we exploit that
$\xi(x,\X) = 0$ whenever $x$ is nonextreme in $\X$ and, roughly
speaking, only points close enough to
the boundary $\S^{d-1}$ have a nonnegligible chance of being extreme
in $\P_{\la}.$
Thus letting $\la\to\infty$ in (\ref{EQEXPL3}), applying the first
part of Lemma~\ref{lem3.4}, using
$\lim_{\la\to\infty} (1-\la^{-\gamma}h)^{d-1} =
1$ and $\lim_{\la\to\infty} g(u(1-\la^{-\gamma}h)) = g(u)$ for
all $u \in\S^{d-1}$
and applying the dominated convergence theorem as in, for example,
Section~3.2 in \cite{SY}, we finally
get from (\ref{EQEXPL3}) the required relation (\ref{EXPconv}).\vspace*{12pt}

\textit{Proof of variance convergence (\ref{VARconv}).}
We have for $g
\in{\cal C}(\B^d)$ and $\xi\in\Xi$, that
\begin{eqnarray*}
&& \la^{-\tau+ 2 \eta[\xi] } \Var[\langle g, \overline{\mu}_\la
^{\,\xi}
\rangle] \\[-1pt]
&&  \qquad = \la^{-\tau+ 2 \eta[\xi] + 1} \int_{\B^d} g^2(x) \E[
\xi(x, \P_\la)^2] (1 - |x|)^{\delta}\,dx
\\[-1pt]
&&  \qquad  \quad {}+ \la^{-\tau+ 2 \eta[\xi] + 2}  \int_{\B^d} \int_{\B^d} g(x) g(y)
 \bigl( \E[\xi(x, \P_\la\cup y) \xi(y, \P_\la\cup x)]\\[-2.5pt]
&& \hspace*{63.3pt} \hphantom{{}\times \int_{\B^d} \int_{\B^d} g(x) g(y)
 \bigl(}\qquad  \quad {}  - \E[\xi(x,
\P_\la)] \E[ \xi(y, \P_\la)] \bigr)\\[-1pt]
&&\qquad  \quad\hphantom{ {}+ \la^{-\tau+ 2 \eta[\xi] + 2}  \int_{\B^d} \int_{\B^d}}{}\times (1 - |x|)^{\delta} (1 -
|y|)^{\delta}\,dx\,dy
\\[-2.5pt]
&&  \qquad := I + \mathit{II}.
\end{eqnarray*}

As in \eqref{new1}, we write term I as\vspace*{-1pt}
\[
I = \la^{-\tau+ 2 \eta[\xi] + 1} \int_{\B^d} g^2(x) \E\bigl[
\xi^{(\la)}\bigl ( \bigl(\0, \la^{\gamma}(1 - |x|)\bigr), \P^{(\la)}\bigr)^2\bigr] (1 -
|x|)^{\delta}\,dx.\vspace*{-1pt}
\]
Now put $h:= \la^\gamma(1 - |x|)$, and write $dx = (1 -
\la^{-\gamma}h)^{d-1}\,d\sigma_{d-1}(u) \la^{-\gamma}\,dh.$ This
transforms I as follows:
\begin{eqnarray*}
I &=& \la^{-\tau+ 1 - \gamma- \delta\gamma} \int_{\S^{d-1}}
\int_0^{\la^\gamma} g^2\bigl(u (1 - \la^{-\gamma} h)\bigr) \E\bigl[\bigl(
\la^{\eta[\xi]} \xi^{(\la)}\bigl((\0, h), \P^{(\la)}\bigr)\bigr)^2\bigr]\\[-1pt]
&&\hphantom{\la^{-\tau+ 1 - \gamma- \delta\gamma} \int_{\S^{d-1}}
\int_0^{\la^\gamma}}{}\times h^{\delta}
(1 -
\la^{-\gamma}h)^{d-1}\,dh\,d\sigma_{d-1}(u).
\end{eqnarray*}
Lemma~\ref{lem3.4} and the moment bounds of Lemma~\ref{MOMBD} give
\[
\lim_{\la\to\infty} \E\bigl[\bigl( \la^{\eta[\xi]} \xi^{(\la)}\bigl((\0, h),
\P^{(\la)}\bigr)\bigr)^2\bigr] = \E\bigl[ \xi^{(\infty)}((\0,h), \P))^2\bigr]
:= \varsigma_{\xi^{(\infty)}}((\0,h)).
\]
Since $\tau:= 1 - \delta\gamma- \gamma$, by the dominated
convergence theorem, we obtain, as $\la\to\infty$, that
%
\begin{equation}\label{termI} I \to\int_{\S^{d-1}} \int_0^{\infty}
g^2(u)\varsigma_{\xi^{(\infty)}}((\0,h)) h^{\delta}\,dh\,d\sigma_{d-1}(u).
\end{equation}
%

We now consider term II. Recall $x:= (u_x,h_x) \in\B^d$ and $y:=
(u_y,h_y) \in\B^d$. We rotate all points in $\P_\la\cup\{x,y \}$
in such a way that $x/|x|$ gets sent to $u_0$.
Denote the rotated point set by $\P'_\la\cup\{x',y' \}$, where $x':=
(\0,h_{x'}), y':= (v_{y'}, h_{y'}),$ with
$h_{x'}= 1 - |x'|, h_{y'}= 1 - |y'|.$

We write term II as
\[
\mathit{II} = \la^{-\tau+ 2 - 2 \gamma\delta} \int_{\B^d} \int_{\B^d}
g(x') g(y') [\cdotsm] h_{x'}^{\delta} h_{y'}^{\delta}\,dx'\,dy',
\]
where
%
\begin{equation}\label{[E]}
 \qquad  [\cdotsm] := \E[ \xi(x', \P_{\la}
\cup
\{y'\} ) \xi(y',
\P_{\la} \cup\{x'\} )] - \E\xi(x', \P_{\la} )
\E\xi( y', \P_{\la}).
\end{equation}

Write
\begin{eqnarray*}
T^{(\la)}(x')&:=& (\0, \la^{\gamma}h_{x'}):= (\0,h);  \qquad   T^{(\la)}(y'):= (\la^{\beta} v_{y'}, \la^{\gamma}h_{y'}):= (v',h');
 \\
  T^{(\la)}(\P'_\la)   &:=& \P'^{(\la)}.
\end{eqnarray*}
Under these
transformations, the expression $[\cdotsm]$ in \eqref{[E]}
transforms to
%
\begin{eqnarray}\label{[EE]}
\hspace*{28pt} [\cdotsm]' &=& \E\bigl[ \la^{\eta[\xi]}
\xi^{(\la)}
\bigl((\0, h), \P'^{(\la)} \cup(v',h') \bigr) \la^{\eta[\xi]} \xi^{(\la)}
\bigl((v',h'),
\P'^{(\la)} \cup(\0, h) \bigr)\bigr]
\nonumber
\\[-6pt]
\\[-10pt]
&&   {}-\E\la^{\eta[\xi]} \xi^{(\la)} \bigl((\0, h), \P'^{(\la)} \bigr)
\la^{\eta[\xi]} \E\xi^{(\la)}\bigl ((v',h'),
\P^{(\la)}\bigr).
\nonumber
\end{eqnarray}

Recalling the definitions of $x'$ and $y'$, we obtain $h_{x'}^{\delta}
= \la^{-\gamma\delta} h^{\delta}$,
$h_{y'}^{\delta} = \la^{-\gamma\delta} (h')^{\delta}$, with
\[
dx' =
(1 - \la^{-\gamma}h)^{d-1}\,d\sigma_{d-1}(u) \la^{-\gamma}\,dh,
\]
%
and
\[
dy' = \la^{-\beta(d-1)}\,dv' \lambda^{-\gamma}\,dh'.
\]

Thus the polynomial $\la$ multiplier in term II gets replaced\vspace*{1pt} by
$\la^{-\tau+ 2 - 2 \gamma\delta}\*\la^{-2 \gamma} \la^{-\beta
(d-1)}$, and so the differential $\la^{-\tau+ 2 - 2 \gamma
\delta}\,dx'\,dy'$ on $\B^d \times\B^d$ in term II transforms to the
differential
%
\begin{equation}\label{differential}\la^{-\tau+ 2 - 2 \gamma
\delta}\la^{-\beta(d-1)} \lambda^{-2\gamma} (1 -
\la^{-\gamma}h)^{d-1}\,d\sigma_{d-1}(u)\,dv'\,dh'\,dh
\end{equation}
on $\S^{d-1}
\times T^{(\la)}(\S^{d-1}) \times[0, \la^\gamma] \times[0,
\la^\gamma].$ The pre-factor in (\ref{differential}) involving
powers of $\la$ reduces to unity in view of the identity $\tau= 2 -
2 \gamma- 2\gamma\delta- \beta(d-1)$. Thus (\ref{differential})
transforms to
%
\begin{equation}\label{trans1} (1 - \la^{-\gamma}h)^{d-1}\,d\sigma_{d-1}(u) \la^{-\gamma}\,dv'\,dh'\,dh.
\end{equation}

For all triples $(v', h', h) \in T^{(\la)}(\S^{d-1}) \times[0,
\la^\gamma] \times[0, \la^\gamma]$, the covariance term $[\cdots ]'$ at
(\ref{[EE]}) may be expressed as
%
\begin{equation}\label{trans2} [\cdotsm]' =
c^{(\la)}((\0,h), (v',h')).
\end{equation}
By Lemma~\ref{lem3.4} we have for
all triples $(v', h', h) \in T^{(\la)}(\S^{d-1}) \times[0,
\la^\gamma] \times[0, \la^\gamma]$, that as $\la\to\infty$
%
\begin{equation}
\label{conv2} [\cdotsm]' = c^{(\la)}((\0,h), (v',h')) \to
\varsigma_{\xi^{(\infty)}}((\0,h), (v',h')).
\end{equation}

Finally, for all $y' \in\B^d$, consider the factor $g(y')$ in the
integrand of term~II.
The factor $g(y')$ transforms to $g((T^{(\la)})^{-1}(v',h'))$.
For all pairs $(v',h') \in T^{(\la)}(\S^{d-1}) \times[0, \la^\gamma
]$, we have
$(T^{(\la)})^{-1}(v',h') \to(u_{x'}, 0)$ as $\la\to\infty$. By
continuity of
$g$ we obtain
%
\begin{equation}\label{conv3} g \bigl(\bigl(T^{(\la)}\bigr)^{-1}(v',h')\bigr) \to
g(u_{x'}, 0)
\end{equation}
as $\la\to\infty$.


Therefore, combining (\ref{trans1}), (\ref{trans2}), 
we may rewrite term II as
\begin{eqnarray*}
\mathit{II}& =& \int_{\S^{d-1}} \int_{T^{(\la)}(\S^{d-1})} \int_0^{\la
^\gamma}
\int_0^{\la^\gamma} g\bigl(\bigl(u (1 - \la^{-\gamma}h)\bigr)
g \bigl(\bigl(T^{(\la)}\bigr)^{-1}(v',h')\bigr)
\\
&&\hphantom{\int_{\S^{d-1}} \int_{T^{(\la)}(\S^{d-1})} \int_0^{\la
^\gamma}
\int_0^{\la^\gamma}}{}\times c^{(\la)} ((\0,h), (v',h'))\\
&&\hphantom{\int_{\S^{d-1}} \int_{T^{(\la)}(\S^{d-1})} \int_0^{\la
^\gamma}
\int_0^{\la^\gamma}}{}\times h^{\delta} h'^{\delta}(1 -
\la^{-\gamma}h)^{d-1}\,dh'\,dh\,dv'\,d\sigma_{d-1}(u).
\end{eqnarray*}

By Lemma~\ref{lem3.5}, the integrand is dominated by the function
\[
(u,v',h',h) \mapsto C h^{\delta} h'^{\delta} \exp\biggl ( {-1 \over
C} \max(|v'|, h, h') \biggr),
\]
which is integrable on $\S^{d-1} \times\R^{d-1}\times(0,\infty
)^2.$ The dominated
convergence theorem, combined with the limits (\ref{conv2}) and
(\ref{conv3}), together with $T^{(\la)}(\S^{d-1}) \uparrow
\R^{d-1}$, show that as $\la\to\infty$, we have
%
\begin{eqnarray}\label{termII}
&& \mathit{II} \to\int_{\S^{d-1}} \int_{\R^{d-1} } \int_0^{\infty}
\int_0^{\infty} g(u)^2 \varsigma_{\xi^{(\infty)}}((\0,h),
(v',h'))\nonumber
\\[-8pt]
\\[-8pt]
&&\hphantom{\mathit{II} \to\int_{\S^{d-1}} \int_{\R^{d-1} } \int_0^{\infty}
\int_0^{\infty}}{}\times h^{\delta} h'^{\delta}\,dh'\,dh\,dv'\,d\sigma_{d-1}(u).
\nonumber
\end{eqnarray}
The second
part of Lemma~\ref{lem3.5} implies that the integral in
(\ref{termII}) is finite. Combining (\ref{termI}) with
(\ref{termII}) gives the desired limit (\ref{VARconv}).\vspace*{12pt}

\textit{Proof of Gaussian convergence (\ref{ConvRate}).}
 The proof uses
the Stein method for dependency graphs and is inspired by the proof
of Theorem 2.1 of \cite{PY5}, which involves a dependency graph
structure on nonscaled sample points in a rectangular solid. Since
the sample points of this paper belong to the unit ball, we find it
more convenient to put a dependency graph on the re-scaled points
$\P^{(\la)}$ in ${\cal R}_\la$. Additionally, we do not use all of
the re-scaled points $\P^{(\la)}$, but only those with\vadjust{\goodbreak} a small
height coordinate. These differences complicate the approach and, in
an effort to make this paper reader friendly, we include the
details.
As in \cite{PY5}, the argument makes use of the
following lemma of Chen and Shao \cite{CSh}. For any random
variable $X$ and any $p > 0$, let
$\|X\|_p := (\E|X|^p)^{1/p}.$ Let $\Phi$ denote the cumulative
distribution function of the standard normal.

%
\begin{lemm}[(See Theorem 2.7 of \cite{CSh})] \label{ChenShao}
Let $2 < q \leq3$. Let $W_i,  i \in{\cal V},$ be random
variables indexed by the vertices of a dependency graph. Let $W =
\sum_{i \in{\cal V} } W_{i}.$ Assume that $\E W^2 = 1, \E W_i =
0,$ and $\|W_i \|_q \leq\theta$ for all $i \in{\cal V}$ and for
some $\theta> 0.$ Then
%
\begin{equation}
\label{CS1} \sup_t | P [ W \leq t ] - \Phi(t) | \leq75 D^{5(q
- 1) } | {\cal V} | \theta^q.
\end{equation}
\end{lemm}

Fix $\xi\in\Xi$ to be one of the basic functionals discussed in
Section~\ref{sec6}. For all $g \in{\cal C}(\B^d)$, we have
\[
\bigl\langle g, \la^{\eta[\xi]} \mu_\la^\xi\bigr\rangle= \sum_{x \in\P
_\la}
\la^{\eta[\xi]} \xi(x, \P_\la) g(x) = \sum_{x' \in\P^{(\la)} }
\la^{\eta[\xi]} \xi^{(\la)}\bigl(x', \P^{(\la)} \bigr) g([T^\la]^{-1} x').
\]
Recalling that $x':= (\la^\beta v_x, h)$, we define for all $L > 0$
and $g \in{\cal C}(\B^d)$,
\[
T_\la^\xi(L, g):= \sum_{x' \in\P^{(\la)}, h \leq L \log\la}
\la^{\eta[\xi]} \xi^{(\la)}\bigl(x', \P^{(\la)} \bigr) g([T^\la]^{-1} x').
\]

By the analog of Lemma 3.2 of \cite{SY}, given $K > 0$ and large, we
may choose $L:=L(K)$ large so that $\langle g, \la^{\eta[\xi]}
\mu_\la^\xi\rangle$ and $T_\la^\xi(L, g)$ coincide everywhere
except on a set with probability $O(\la^{-K})$. It follows that
$\langle g, \la^{\eta[\xi]} \mu_\la^\xi\rangle$ and $T_\la^\xi(L,
g)$ have the same asymptotic distribution as $\la\to\infty$, and it
may be shown that they share the same variance asymptotics. It
suffices to find a rate of convergence to $N(0,1)$ for
$(T_\la^\xi(L, g) - \E T_\la^\xi(L, g))/ \sqrt{ \Var[T_\la^\xi(L,
g)] }$.

To prepare for dependency graph arguments, we put $\rho_\la:= L \log
\la$, $L$ a constant, and we subdivide $\la^\beta\B_{d-1}(\pi)$
into $V(\la):= (2\pi\la^{\beta})^{d-1}(\rho_\la)^{-(d-1)}$
sub-cubes $Q_i,
i = 1,\ldots,V(\la)$,
of edge length $\rho_\la$ and of volume
$(\rho_\la)^{d-1}$. Enumerate the points $\P^{(\la)} \cap[Q_i
\times L \log\la]$ by $\{X'_{ij} \}_{j=1}^{N_i}$ so that
\[
T_\la^\xi(L, g)= \sum_{i=1}^{V(\la)} \sum_{j=1}^{N_i}
\la^{\eta[\xi]} \xi^{(\la)}\bigl(X'_{ij}, \P^{(\la)} \bigr) g([T^\la]^{-1}
X'_{ij}).
\]

The random variable $N_i$ is Poisson whose mean $\nu_i$ equals the
$\P^{(\la)}$ intensity measure of the rectangular solid $Q_i \times
L \log\la$, and thus $\nu_i$ is bounded by the product of $\Vol(Q_i
\times L \log\la)$ and the maximum of the intensity of $\P^{(\la)}$
on this solid. Recalling the intensity of
$\P^{(\la)}$ at (\ref{DENS1}), we obtain 
%
\[
\nu_i := \E N_i \leq C\bigl( \Vol(Q_i \times L \log\la)\bigr) (L \log
\la)^{\delta} = C (\rho_\la)^{d + \delta}.
\]

The following result is the analog of Lemma 4.3 of \cite{PY5} and is
proved similarly. For $1 \leq i \leq V(\la)$, and
for $j \in\{
1,2,\ldots\}$, we define for the fixed functional $\xi$,
\[
\xi_{i,j} :=
\cases{\displaystyle \la^{\eta[\xi]} \xi^{(\la)}
\bigl(X'_{i,j}, \P^{(\la)} \bigr) ,&\quad  if    $N_i \geq j,   X'_{i,j} \in\la^\beta\B_{d-1}(\pi)
\times[0, L \log\la] $,\cr\displaystyle
0 ,&\quad  otherwise.
}
\]
With $\xi\in\Xi$ still fixed, note that $\xi$ satisfies the moment
condition (\ref{LIMITBD}) for all $p > q \geq1$.

\begin{lemm} \label{lem4.3} With $p > q \geq1$,
there exists $C:= C(p,q)$, such that for $1 \leq i \leq V(\la)$, we
have
\[
 \Biggl\| \sum_{j = 1}^{\infty} |\xi_{i,j} |\Biggr\|_q
\leq
C \rho_\la^{(d + \delta)(p + 1)/p}.
\]
\end{lemm}
%

Continuing with the dependency graph arguments, we let $p > q$ and $q
\in(2,3]$. Recall from Section~\ref{sec6} that $R^{\xi^{(\la)}}(x)$ is
localization radius for the functional $\xi^{(\la)}$ if and only if a.s.
\[
\xi^{(\la)}\bigl (x, \P^{(\la)} \bigr) =
\xi^{(\la)}_{[r]}\bigl (x, \P^{(\la)} \bigr)
\]
for all $r \geq
R^{\xi^{(\la)}}(x)$. Put $U(t):= \sup_{\la\geq1, x \in{\cal
R}_\la} P[R^{\xi^{(\la)}}(x) > t],$ which is the analog of the
$\tau$ function defined at (2.2) of \cite{PY5}. With the choice
$\rho_\la= L \log\la$, Lemma~\ref{LocalLem} implies that for $L$
large, we have that $U(\rho_\la)$ has polynomial decay of high order,
and so
we have
%
\begin{eqnarray}\label{dis4.8}
  \quad \qquad V(\la) \rho_\la^{(d + \delta)(p +
1)/p}  \bigl(\la^{\beta(d-1)}(\log\la)^{1 + \delta}
U(\rho_\la) \bigr)^{(q-2)/2q} &<& \la^{-3 - (\beta d/2)}  \quad \mbox{and}   \nonumber
\\[-8pt]
\\[-8pt]
  \quad \qquad U(\rho_\la) &<& \la^{-\beta(d-1) - 3},
\nonumber
\end{eqnarray}
which is
the analog of display (4.8) in \cite{PY5}. 
We also have $\rho_\la^{d + \delta} < C \la^{p/(p + 2)}$, the analog
of display (4.9) in \cite{PY5}.

For all $1 \leq i \leq V(\la)$ and all $j =1,2,\ldots,$ let
$R^{(\la)}_{i,j}$ denote the radius of stabilization of
$\xi^{(\la)}$ at $X'_{i,j}$
if $1 \leq j \leq N_i$ and $X'_{i,j} \in\la^\beta\B_{d-1}(\pi)
\times[0,
\la^\gamma]$; let $R_{i,j}$ be zero otherwise.
Let $E_{i,j} := \{ R^{(\la)}_{i,j} \leq\rho_\la\}$. Let $E_\la
:= \bigcap_{i=1}^{V(\la)} \bigcap_{j=1}^{\infty} E_{i,j}.$ Then
%
\begin{eqnarray}
\label{dis4.11}
 P[E_\la^c] &\leq&\E \Biggl[ \sum_{i=1}^{V(\la)}
\sum_{j=1}^{N_i} \1({E_{i,j}^c})  \Biggr] \nonumber\\
&=& \int_{\la^\beta
\B_{d-1}(\pi) \times[0, L \log\la]} P\bigl[R^{\xi^{(\la)}}(x) \geq
\rho_\la\bigr] h^\delta \,dv \,dh
\\
&\leq& C \la^{\beta(d-1)}(\log\la)^{1 + \delta} U(\rho_\la).
\nonumber
\end{eqnarray}

We have thus
\[
T_\la^\xi(L, g):= \sum_{i=1}^{V(\la)} \sum_{j=1}^{N_i} \la^{\eta
[\xi]}
\xi^{(\la)}\bigl(X'_{i,j}, \P^{(\la)} \bigr) g([T^\la]^{-1} (X'_{i,j})).
\]
To obtain rates of normal approximation, we consider a modified
version of $T_\la^\xi(L, g)$ having more independence between
terms, namely
\[
T'_\la(L, g):= \sum_{i=1}^{V(\la)} \sum_{j=1}^{N_i} \la^{\eta[\xi]}
\xi^{(\la)}\bigl(X'_{i,j}, \P^{(\la)} \bigr) {\bf1}(E_{i,j}) g([T^\la]^{-1}
(X'_{i,j})).
\]
For all $1 \leq i\leq V(\la)$, define
\[
S_i := S_{Q_i}:= (\Var T'_{\la}(L,g))^{-1/2} \sum_{j =1}^{N_i}
\la^{\eta[\xi]} \xi^{(\la)}\bigl(X'_{i,j}, \P^{(\la)} \bigr) \mathbf{1}(E_{i,j}) g([T^\la]^{-1} (X'_{i,j}))
\]
and put
%
\begin{equation}\label{Sdef} \quad  S := (\Var T'_{\la}(L,g))^{-1/2}
\bigl(T'_{\la}(L,g) - \E T'_{\la}(L,g)\bigr) = \sum_{i=1}^{V(\la)} (S_i - \E
S_i).
\end{equation}
Note that $\Var S = \E S^2 = 1$.

We define a graph $G_\la:= ({\cal V}_\la, {\cal E }_\la)$ as
follows. The set ${\cal V}_\la$ consists of the sub-cubes
$Q_1,\ldots,Q_{V(\la) }$ and edges $(Q_i, Q_j)$ belong to ${\cal E
}_\la$ if $d(Q_i, Q_j) \leq2 \rho_\la$, where $d(Q_i,Q_j) := \inf
\{ |x - y|, x \in Q_i, y \in Q_j \}$. To prepare for dependency
graph arguments, we make the following five observations, paralleling
those in
\cite{PY5}:
\begin{longlist}[(iii)]
\item[(i)] $V(\la) := | {\cal V}_\la|$.

\item[(ii)] Since the number of cubes in $Q_1,\ldots,Q_V$ distant at most $2
\rho_\la$ from a given cube is bounded by $ 5^d$,
it follows that the maximal degree $D$ of $G_\la$ satisfies $D:=
D_\la\leq5^d$.

\item[(iii)] For all $1 \leq i \leq V(\la)$ and all $q \geq1$, we have, by
Lemma~\ref{lem4.3},
%
\begin{eqnarray}\| S_i\|_q &\leq& C (\Var
T'_{\la}(l,g))^{-1/2}
 \Biggl\| \sum_{j=1}^{\infty} |\xi_{i,j} |  \Biggr\|_q
\nonumber
\\[-8pt]
\\[-8pt]&\leq& C (\Var T'_{\la}(L,g))^{-1/2} \rho_\la^{(d + \delta)(p+1)/p}.
\label{0511}
\nonumber
\end{eqnarray}

\item[(iv)] $T'_\la(L,g)$ is the sum of $V(\la)$ random variables,
which, by the case $q = 2$ of Lemma~\ref{lem4.3}, each have a
variance bounded by a constant multiple of $\rho_\la^{2(d +
\delta)(p + 1)/p}.$ The covariance of any pair of the $V(\la)$
random variables is zero when the indices of the random variables
correspond to nonadjacent cubes. For adjacent cubes, the
Cauchy--Schwarz inequality implies that the covariance is also
bounded by a constant multiple of $\rho_\la^{2(d + \delta)(p +
1)/p}.$ This gives the analog of (4.13) of \cite{PY5}, %
namely
%
\begin{eqnarray}\Var[T'_\la(L,g)] = O\bigl(
\rho_\la^{2(d + \delta)(p+1)/p} V(\la)\bigr). 
\label{dis4.13}
\end{eqnarray}
%

\item[(v)] $\Var[T'_{\la}(L,g) ]$ is close to $\Var[T_\la^\xi(L,g)]$ for
$\la$ large. We require more estimates to show this. Note that
$| T'_{\la}(L,g) - T_\la^\xi(L,g)| = 0$ except possibly on the set
$E_\la^c $. Lemma~\ref{lem4.3}, along with Minkowski's inequality,
yields the upper bound
%
\begin{eqnarray}
 \quad \Biggl\| \sum_{i=1}^{V(\la)} \sum_{j = 1 }^{N_i}
\bigl\vert
\la^{\eta[\xi]} \xi^{(\la)}\bigl(X'_{i,j}, \P^{(\la)} \bigr)\bigr \vert
 \Biggr\|_q &\leq& C V(\la) \rho_\la^{(d + \delta)(p+1)/p }\nonumber
 \\[-8pt]
 \\[-8pt]
 \quad &=& C \la^{\beta(d-1)} \rho_\la^{-(d-1 )} \rho_\la^{(d + \delta
)(p+1)/p }. \label{dis4.14}
\nonumber
\end{eqnarray}
Since $T_\la^\xi(L,g) =
T'_\la(L,g)$ on the event $E_\la$, as in \cite{PY5}, the H\"older
and Minkowski inequalities yield
\begin{eqnarray*}
\|T_{\la}^\xi(L,g) - T'_{\la}(L,g)\|_2 &\leq&
\|T^\xi_{\la}(L,g) - T'_{\la}(L,g)\|_q
P[E_\la^c]^{(1/2)-(1/q)}
 \\
&\leq&
\bigl(\|T^\xi_\la(L,g)\|_q+ \|T'_\la(L,g)\|_q\bigr) P[E_\la^c]^{(q-2)/(2q)}.
\end{eqnarray*}
\end{longlist}

Hence, by (\ref{dis4.11}) and the first inequality in
(\ref{dis4.14}),
\begin{eqnarray*}
&&\|T^\xi_{\la}(L,g) - T'_{\la}(L,g)\|_2\\
&& \qquad \leq C V(\la) \rho_\la^{(d + \delta)(p + 1)/p}
 \bigl(\la^{\beta(d-1)} (\log\la)^{1 + \delta}
U(\rho_\la) \bigr)^{(q-2)/2q}.
\end{eqnarray*}
By (\ref{dis4.8}) this yields
%
\begin{eqnarray}\|T^\xi_{\la}(L,g) - T'_{\la}(L,g)\|_2 \leq C \la
^{-3 - \beta
d/2}
\label{dis4.15}
\end{eqnarray}
which clearly implies
%
\begin{equation}
\label{L1bound}
\E[|T'_{\la}(L,g) -T_{\la}(L,g) |] \leq C \la^{-3},
\end{equation}
which we use later. As in \cite{PY5}, we obtain the analog of
(4.17) of \cite{PY5}, that is,
%
\begin{equation}
\label{dis4.17}
| \Var[T^\xi_{\la}(L,g) ] -
\Var[T'_\la(L,g)] | \leq C\la^{-2},
\end{equation}
concluding observation (v).

  We may now use Lemma~\ref{ChenShao} and dependency graph
arguments to
establish the error bound (\ref{ConvRate}). We apply the bound
(\ref{CS1}) of Lemma~\ref{ChenShao} to $W_i := S_i - \E S_i,  1
\leq i \leq V(\la),$
with
\[
\theta: = C ( \Var T'_{\la}(L,g) )^{-1/2}
\rho_\la^{(d + \delta)(p+1)/p}.
\]

Observe that $\E W_i = 0$, $\E(\sum_{i=1}^{V(\la)} W_i)^2 = 1$,
$\|W_i\|_q \leq\theta$ and recall from \eqref{Sdef} that $S =
\sum_{i=1}^{V(\la)} W_i$. Lemma~\ref{ChenShao} along with
observation (i)\vadjust{\goodbreak} above yields the counterpart of (4.18) of
\cite{PY5}, namely
%
\begin{eqnarray}\label{diffbound}
\sup_t  | P[S \leq t] - \Phi(t)  |
&\leq& C V(\la) (\Var T'_\la(L,g))^{-q/2} \rho_\la^{q (d + \delta)(p +
1)/p}
\nonumber
\\[-8pt]
\\[-8pt]
&\leq& C V(\la) (\Var T^\xi_\la(L,g))^{-q/2} \rho_\la^{q (d + \delta
)(p +
1)/p}.
\nonumber
\end{eqnarray}
The last line follows since by \eqref{VARconv} we have
$\Var[T^{\xi}_{\la}(L,g)] = \Theta(\la^{\tau}), \tau\in(0,1),$ and
thus by (\ref{dis4.17}) we get for $\la$ large that
$\Var[T'_{\la}(L,g)] \geq
\Var[T^\xi_\la(L,g) ]/2$.
Now put $q
= 3$ in (\ref{diffbound}). Since $T^\xi_\la(L,g)$ and $\langle g,
\la^{\eta[\xi]} \mu_\la^\xi\rangle$ have the same variance
asymptotics, it follows from the assumption $\sigma^2(\xi^{(\infty)}
) > 0$ that when $q = 3$, we have $(\Var T^\xi_\la(L,g))^{-q/2} =
\Theta(\la^{-3 \tau/2}).$ Since $V(\la) = \la^\tau
\rho_\la^{-(d-1)}$ and since $q/p < 1$, display
(\ref{diffbound}) becomes
\[
\sup_t  | P[S \leq t] - \Phi(t)  | \leq C \la^{-\tau/2}
(\log\la)^{3d + 4 \delta+ 1}.
\]



This gives a rate of convergence for $S$, as defined at
\eqref{Sdef}. Following verbatim the last part of the proof of
Theorem 2.1 of \cite{PY5} [starting three lines after (4.18) of that
paper], we deduce a rate of convergence for $T^\xi_\la(L,g)$, namely
\[
\sup_t  \biggl| P \biggl[ {T_\la^\xi(L,g) - \E T_\la^\xi(L,g) \over
\sqrt{\Var T_\la^\xi(L,g) } } \leq t  \biggr] - P[N(0,1) \leq t]
 \biggr| \leq C \la^{-\tau/2} (\log\la)^{3d + 4 \delta+ 1}.
\]
This yields (\ref{ConvRate}), concluding the proof of Theorem
\ref{MainCLT}.





\textit{Positivity of asymptotic
variances.} For $\xi^{(\infty)}\in\Xi^{(\infty)}$, we now consider
the question whether $\sigma^2(\xi^{\infty})$ is strictly positive.
Fortunately, the variances $\sigma^2(\xi_r^{(\infty)})$, $\sigma
^2(\xi_s^{(\infty)})$
and $\sigma^2(\xi^{(\infty)}_{\vartheta_k}),  k \in\{1,\ldots
,d-1\}$, admit alternative expressions enjoying
monotonicity properties in the underlying Poisson input process $\P$, enabling
us to use suitable positive correlation inequalities and to conclude
the required
positivity for variance densities. The underlying Poisson input process
$\P$
depends on the parameter $\delta$ [recall (\ref{basicPPP})], and the
following lemma
holds for all $\delta> 0$. More precisely, we have:
%
\begin{lemm}\label{AVARREPR} We have
\begin{eqnarray*}
\sigma^2_s&:=& \sigma^2\bigl(\xi_s^{(\infty)}\bigr) = \int_{\R^{d-1}} \Cov
(\partial\Psi(0),\partial\Psi(v))\,dv,
\\
\sigma^2_r&:=&\sigma^2\bigl(\xi_r^{(\infty)}\bigr) = \int_{\R^{d-1}} \Cov
(\partial\Phi(0),\partial\Phi(v))\,dv
\end{eqnarray*}
and
\[
\sigma^2_k:=\sigma^2\bigl(\xi^{(\infty)}_{\vartheta_k}\bigr) = \int_{\R^{d-1}}
\Cov \biggl( \int_0^{\Phi(0)} \vartheta^{\infty}_k((\0,h))\,dh,
\int_0^{\Phi(v)} \vartheta^{\infty}_k((v,h))\,dh  \biggr)\,dv.
\]
\end{lemm}
\begin{pf} We only consider the functional $\xi^{(\infty)}_s,$ as
the remaining cases are analogous.
Recalling \eqref{SOScorr1}--\eqref{VARasympt1}, the general theory
of stabilizing functionals
(see, e.g., \cite{BY2}, \cite{PeEJP}) shows that if $\xi^{(\infty)} $
is a generic exponentially stabilizing functional on the Poisson
input $\P$ on $\R^{d-1} \times\R_+$, then
\[
\lim_{T\to\infty} \frac{1}{T^{d-1}}
\Var \biggl( \sum_{x = (v,h) \in\P,  v \in[0,T]^{d-1}}
\xi^{(\infty)}(x,\P)  \biggr) = \sigma^2\bigl(\xi^{(\infty)}\bigr),
\]
that is to say, the scaled variance limit of $\sum_{x = (v,h) \in
\P,  v \in[0,T]^{d-1}} \xi^{(\infty)}(x,\P)$ coincides with
\[
\lim_{\la\to\infty} \la^{-\tau} \Var \biggl( \sum_{x \in
\P^{(\la)}} \xi^{(\la)}\bigl(x, \P^{(\la)}\bigr)  \biggr).
\]

Since $\xi_s^{(\infty)}$ is an exponentially stabilizing functional
on the Poisson input $\P$ (recall Lemma
\ref{LocalLem}), it follows 
that
$\sigma^2(\xi_s^{(\infty)})$ is the asymptotic variance density for
$\xi^{(\infty)}_s,$
that is to say,
\[
\sigma^2\bigl(\xi_s^{(\infty)}\bigr) = \lim_{T\to\infty} \frac{1}{T^{d-1}}
\Var\biggl ( \sum_{x = (v,h) \in\P,  v \in[0,T]^{d-1}}
\xi^{(\infty)}_s(x,\P)  \biggr).
\]
For $x := (v,h) \in\ext(\Psi) = \Vertices(\Phi)$ denote by $V[x] :=
V[x;\P]$
the set of all $v' \in\R^{d-1}$ for which there exists $h'$ with
$(v',h') \in\F^{\infty}(x,\P)$---in other words, $V[x]$ is the
spatial projection of all faces $f$ of $\Phi$ with
$x = \Top(f).$ Clearly, $\{ V[x],  x \in\ext(\Psi) \}$, forms a
tessellation of $\R^{d-1}.$
Thus, by definition of $\xi_s^{(\infty)},$
\[
\sigma^2\bigl(\xi_s^{(\infty)}\bigr) = \lim_{T\to\infty} \frac{1}{T^{d-1}}
\Var
 \biggl(\sum_{x = (v,h) \in\Vertices(\Psi),  v \in[0,T]^{d-1}}
\int_{V[x]} \partial\Psi(u)\,du  \biggr).
\]
Consequently,
\begin{eqnarray*}
\sigma^2\bigl(\xi_s^{(\infty)}\bigr) &=& \lim_{T\to\infty} \frac{1}{T^{d-1}}
\Var \biggl(\int_{[0,T]^{d-1}} \partial\Psi(u)\,du  \biggr)\\
&=& \lim_{T\to\infty} \frac{1}{T^{d-1}} \int_{([0,T]^{d-1})^2} %
\Cov(\partial\Psi(u),\partial\Psi(u'))\,du'\,du,
\end{eqnarray*}
%
where the existence of the integrals 
follows from the exponential localization
of $\xi^{(\infty)}_s,$ as stated in Lemma~\ref{LocalLem}, implying
the exponential decay of correlations.
Further, by stationarity of the process $\partial\Psi(\cdot),$ we obtain
\begin{eqnarray*}
\sigma^2(\xi_s^{(\infty)}) &=& \lim_{T\to\infty} \frac{1}{T^{d-1}}
\int_{[0,T]^{d-1}} \int_{[0,T]^{d-1}}
\Cov\bigl(\partial\Psi({\bf0}),\partial\Psi(u'-u)\bigr)\,du'\,du
\\
&=& \lim_{T\to\infty} \int_{[-T,T]^{d-1}} \frac{\Vol([0,T]^{d-1}
\cap([0,T]^{d-1} + w))}{T^{d-1}}
\\
&&\hphantom{\lim_{T\to\infty} \int_{[-T,T]^{d-1}}}{}\times\Cov(\partial\Psi({\bf0}),\partial\Psi(w))\,dw
\\
&=& \int_{\R^{d-1}} \Cov(\partial\Psi({\bf0}), \partial\Psi(w))\,dw
\end{eqnarray*}
as required, with the penultimate equality following again by
exponential localization
of $\xi^{(\infty)}_s$, implying the exponential decay of correlations
and thus
allowing us to apply dominated convergence theorem to determine the limit
of integrals. This completes the proof of Lemma~\ref{AVARREPR}.
\end{pf}

Observe that for each $v,$ $\partial\Psi(v),  \partial\Phi(v)$ as
well as
$\int_0^{\Phi(v)} \vartheta^{\infty}_k((v,h))\,dh$ are all nonincreasing
functionals of $\P$ and therefore
\[
\Cov(\partial\Psi(\0),\partial\Psi(v)) \geq0, \qquad
  \Cov(\partial\Phi(\0),\partial\Phi(v)) \geq0
\]
 and
\[
\Cov \biggl( \int_0^{\Phi(\0)} \vartheta^{\infty}_k((\0,h))\,dh,
\int_0^{\Phi(v)} \vartheta^{\infty}_k((v,h))\,dh  \biggr) \geq0
\]
for all $v \in\R^{d-1}$ in view of the positive correlations property
of Poisson point
processes; see Proposition 5.31 in \cite{RES}. It is also readily seen
that these
covariances are not identically zero, because for $v=0$ they are just
variances of
nonconstant random variables and, depending continuously on~$v$,
they are strictly positive on a nonzero measure set of
$v'$s. 
Thus, the integrals in the variance expressions
given in Lemma~\ref{AVARREPR} are all strictly positive. Consequently,
we have
\begin{coro}\label{VARPOS} For all $\delta> 0$, the variance densities
$\sigma^2(\xi_r^{(\infty)}), \sigma^2(\xi_s^{(\infty)})$
and $\sigma^2(\xi^{(\infty)}_{\vartheta_k}),  k \in\{1,\ldots
,d-1\}$ are all strictly positive.
\end{coro}


\begin{remark*} When $\delta= 0$ the variance positivity for
$\sigma^2(\xi^{(\infty)}_{\vartheta_k})$ has been established in a
slightly different, but presumably equivalent, context (binomial
input) in \cite{BFV}, Theorem 1.

We also believe that for all $\delta> 0$, the variance density
$\sigma^2(\xi_{f_k}^{(\infty)})$ is strictly positive as well---this is because of the asymptotic nondegeneracy of the
corresponding so-called add-one cost functional
\cite{BY2,PeEJP,PY1,PY5}. 
However, making this intuition precise requires additional
technical considerations, as does extending the important work of
Reitzner \cite{Re} to the case $\delta> 0$, which shows strict
variance positivity for $\delta= 0$.
\end{remark*}



\textit{Variance asymptotics
and central limit theorems for mean widths, volumes, intrinsic
volumes and $k$-face functionals.} We now deduce from Theorem
\ref{MainCLT} and Corollary~\ref{VARPOS} the limit theory for the
convex hull functionals\vadjust{\goodbreak} described at the outset of this paper. We
require some preliminary observations which will also be needed in
Section~\ref{InvPrinc}. Define for $v \in\R^{d-1}$ the defect width
functional
%
\begin{equation}\label{HSV}
H^{\xi_s}_{\la}(v) := \sum_{x \in\P_{\la},  x/|x| \in\exp([\0,v])} \xi_s(x;\P_{\la})
\end{equation}
and the defect volume functional
%
\begin{equation}\label{HRV}
H^{\xi_r}_{\la}(v) := \sum_{x \in\P_{\la},  x/|x| \in\exp([\0,v])}
\xi_r(x,\P_{\la}).
\end{equation}
%

The next lemma shows that the centered defect width functional
approximates its asymptotic counterpart $W_\la$ and likewise for the
centered defect volume functional.

\begin{lemm} \label{PClem} We have, uniformly in $v \in\R^{d-1}$,
%
\begin{equation}
\label{Diff1} \lim_{\la\to\infty} \la^{\zeta/2}
\bigl|\bigl(H^{\xi_s}_{\la}(v) - \E H^{\xi_s}_{\la}(v)\bigr) - \bigl(W_\la(v) - \E
W_\la(v)\bigr)\bigr| \stackrel{{P}}{=}0
\end{equation}
and
%
\begin{equation}\label{Diff2} \lim_{\la
\to\infty}
\la^{\zeta/2} \bigl|\bigl(H^{\xi_r}_{\la}(v) - \E H^{\xi_r}_{\la}(v)\bigr) -
\bigl(V_\la(v) - \E V_\la(v)\bigr)\bigr| \stackrel{{P}}{=}0.
\end{equation}
\end{lemm}

\begin{pf} We first prove \eqref{Diff1}.
It is enough to show the
two following limits:
%
\begin{equation}\label{Diff02} \lim_{\la\to\infty}
\la^{\zeta/2} \bigl|\bigl(H^{\xi_s}_{\la}(v) - \E H^{\xi_s}_{\la}(v)\bigr) -
\bigl(\Vol({\mathcal C}(v))-\E\Vol({\mathcal C}(v))\bigr)\bigr| \stackrel{{P}}{=}0
\end{equation}
and
%
\begin{equation}
\label{Diff03} \lim_{\la\to\infty} \la^{\zeta/2} \bigl|\bigl(\Vol
({\mathcal
C}(v))-\E\Vol({\mathcal C}(v))\bigr)
- \bigl(W_\la(v) - \E
W_\la(v)\bigr)\bigr| \stackrel{{P}}{=}0,
\end{equation}
where ${\mathcal C}(v):=[\B
^d\setminus
F({\mathcal P}_{\lambda})]\cap\cone(\exp([\0,v]))$.

We start by proving \eqref{Diff02}. For all $v \in\R^{d-1}$,
$|H^{\xi_s}_{\la}(v) - \Vol({\mathcal C}(v))|$
is bounded by the volume of the set
\begin{eqnarray*}
  \Delta_\la(v)  &:=&[\B^d\setminus F({\mathcal P}_{\la})]\\
  &&{}\cap
\biggl (\cone(\exp([\0,v]))  \Delta\bigcup_{x\in{\mathcal
P}_{\la}\cap\cone(\exp([\0,v]))}\cone({\mathcal F}(x,{\mathcal
P}_{\la})) \biggr).
\end{eqnarray*}
Let
\[
{\mathcal F}_\la(v):= \bigcup\{f\in{\mathcal F}_{d-1}(K_{\la
})\dvtx f\cap\partial
\cone(\exp([\0,v])) \ne\varnothing\}.
\]
By the usual scaling via the
transformation $T^\la$, we get that $T^\la(\Delta_\la(v))$ is a
solid $D_\la(v)$, with $D_\la(v) \subset T^{\la}({\mathcal
F}_\la(v))$. Consider the $(d-2)$-dimensional surface given by
\[
S_\la(v):= \partial \bigl( T^\la[ \cone(\exp([\0,v])) \cap
\S^{d-1}] \bigr).\vadjust{\goodbreak}
\]
Then the maximal height coordinate of $D_\la(v)$,
with respect to the surface $S_\la(v)$, satisfies the exponential
decay \eqref{SuperexpBd}. Also, the maximal spatial distance between
$T^{\la}({\mathcal F}_\la(v))$ and $S_\la(v)$ has exponentially
decaying tails, as in Lemma~\ref{LocalLem}.
By mimicking the proof of Theorem
\ref{MainCLT}, but with now $\tau$ taken to be $\beta(d-2)$ instead
of $\beta(d-1)$ , it follows that $\la^{-\beta(d-2)/2 +
\eta}(\Delta_\la(v) - \E\Delta_\la(v))$ converges to
a normal random variable. Since 
$ \la^{\zeta/2}:= \la^{\beta(d-1)/2 + \gamma} =
o(\la^{-\beta(d-2)/2 + \eta})$, this gives the convergence
\eqref{Diff02}.

We prove now \eqref{Diff03}. We deduce from \eqref{ajout} that
\begin{eqnarray*}
W_{\la}(v)-\Vol({\mathcal C}(v))&=& \int_{\exp([\0,v])}
 \biggl(s_{\la}(u)-\frac{1-(1-s_{\la}(u))^d}{d} \biggr)\,d\sigma_{d-1}(u)
\\
&=&O \biggl(\int_{\exp([\0,v])}s_{\la}^2(u)\,d\sigma_{d-1}(u) \biggr).
\end{eqnarray*}
For every $x\in{\mathcal P}_{\la}$, let
\[
\widetilde{\xi}_{s}(x,{\mathcal P}_{\la}):=
\int_{\cone({\mathcal F}(x,{\mathcal P}_{\la})\cap
\exp([\0,v]))}s_{\la}^2(u,{\mathcal P}_{\la})\,d\sigma_{d-1}(u).
\]
In
particular, we have
\[
\int_{\exp([\0,v])}s_{\la}^2(u)\,d\sigma_{d-1}(u)=\sum_{x\in
{\mathcal P}_{\la}}\widetilde{\xi}_{s}(x,{\mathcal P}_{\la}).
\]
Exactly as for $\xi\in\Xi$, where $\Xi$ is the class of
functionals defined in Section~\ref{StabRepr}, $\widetilde{\xi}_s$
has an associated scaling prefactor $\la^{\eta[\widetilde{\xi}_s]}$
with $\eta[\widetilde{\xi}_s]=\beta(d-1)+2\gamma$ (recall that
$s_{\la}$ is of order $\la^{\gamma}$). Moreover $\widetilde{\xi}_s$
is seen to satisfy the exponential decay \eqref{LocalExpr}.
Following verbatim the proof of Theorem~\ref{MainCLT}, we obtain
that
$\la^{-\tau/2+\eta[\widetilde{\xi}_s]}\int_{\exp([\0
,v])}(s^2(u,{\mathcal
P}_{\la}) - \E s^2(u,{\mathcal P}_{\la}))\,d\sigma_{d-1}(u)$
converges to a normal random variable. Since
$\tau/2-\eta[\widetilde{\xi}_s]<-\zeta/2$, the asserted limits
\eqref{Diff03} and \eqref{Diff1} follow.

To prove \eqref{Diff2}, it suffices to recall from (\ref{RGR}) that
\[
 V_{\la}(v) := \int_{\exp([\0,v])} r_{\la}(u)\,d\sigma_{d-1}(u)
 \]
  and
to follow arguments similar to those given above for
$H^{\xi_s}_\la.$ This completes the proof of Lemma~\ref{PClem}.
\end{pf}

Letting $d \kappa_d$ be the total surface measure of $\S^{d-1}$ and
recalling $\zeta:= (d+3)/ (d+1+2\delta)$ from (\ref{ZETA}), the
following theorem gives scalar variance asymptotics and scalar
central limit theorems for the basic functionals discussed in the
\hyperref[sec1]{Introduction}.


\begin{theo}\label{VWSCALAR}
\textup{(i)} The volume functional $V(K_\la)$ satisfies
%
\begin{equation}\label{scalarv}
\lim_{\la\to\infty} \la^{\zeta} \Var[V(K_\la)] = \sigma^2_V :=
\sigma^2\bigl(\xi_r^{(\infty)}\bigr)\,d\kappa_d\vadjust{\goodbreak}
\end{equation}
and
%
\begin{equation}\label{scalarclt}
\la^{\zeta/2} \bigl(V(K_\la) - \E V(K_\la)\bigr) \tod N(0,\sigma^2_V),
\end{equation}
where $\sigma^2_V$ is strictly positive.

\textup{(ii)} The volume functional $V_{\la}(\infty)$ satisfies the identical
asymptotics whereas the mean width functional $W_{\la}(\infty)$ and
the intrinsic volume functionals $V_k(K_{\la}),\break  k \in
\{1,\ldots,d-1\}$ satisfy \eqref{scalarv} and \eqref{scalarclt} with
strictly positive variances $\sigma^2_W :=
\sigma^2(\xi_s^{(\infty)})\,d\kappa_d$ and $\sigma_{V_k}^2 :=
\sigma^2(\xi^{(\infty)}_{\vartheta_k})\,d\kappa_d$, respectively.
\end{theo}

\begin{remark*} Recalling (\ref{SGRESC}) and setting
$\delta= 0$, Theorem~\ref{VWSCALAR} yields the
asserted variance limits (\ref{widthvar}), (\ref{volvar}) and (\ref
{intrinsicvar}).
In Section~\ref{InvPrinc} we shall show convergence of the $\R
^{d-1}$-indexed processes
$W_{\la}(\cdot)$ and $V_{\la}(\cdot)$.
\end{remark*}

\begin{pf*}{Proof of Theorem~\ref{VWSCALAR}} To prove the assertion for $V(K_\la)$, it suffices to
put $g \equiv1$ and $\xi\equiv\xi_r$ in Theorem~\ref{MainCLT}, to
recall that $\Vol(\B^d\setminus K_{\la})= \sum_{x \in\P_{\la}}
\xi_r(x,\P_{\la})$, and to use $\la^{-\tau} \la^{2 \eta[\xi_r]} =
\la^{(d+3)/(d + 1 + 2 \delta)}$. Corollary~\ref{VARPOS} yields
positivity of the limiting variance $\sigma^2_V$. The limit theory
for $V_{\la}(\infty)$ holds since we may follow verbatim the proof of Lemma
\ref{PClem} to show that $\sum_{x\in\P_{\la}}\xi_r(x,\P_\la)$
approximates $V_\la(\infty)$.

Similarly, to prove the asserted limit theory for $W_{\la}(\infty)$,
we put $g \equiv1$ and $\xi\equiv\xi_s$ in Theorem~\ref{MainCLT},
we use $\la^{-\tau} \la^{2 \eta[\xi_s]} = \la^{(d+3)/(d + 1 + 2
\delta)}$, and we follow verbatim the proof of Lemma~\ref{PClem} to
show that $\sum_{x\in\P_{\la}}\xi_s(x,\P_\la)$ approximates
$W_\la(\infty)$.
Corollary~\ref{VARPOS} yields positivity of the limiting variance
$\sigma^2_W$. Finally, the asserted limit theory for $V_k(K_{\la})$
follows by putting $g \equiv1$ and $\xi\equiv\xi_{\vartheta_k}$
in Theorem~\ref{MainCLT} and using Corollary~\ref{VARPOS} to deduce
the positivity of the limiting variance.
\end{pf*}

Next, using (\ref{FKEMAPPROX}) and Theorem~\ref{MainCLT} we
obtain the limit theory for the $k$-face empirical measures $\mu
^{f_k}_{\la}$ defined at
$\eqref{FaceEmpM}$.
%
\begin{theo} \label{Vkface}
For each $k \in\{0,\ldots,d-1\},$ the $k$-face empirical measures
$\mu^{f_k}_{\la}$
satisfy the measure-level variance asymptotics and central limit
theorem with scaling exponent
$\tau/2$ and with variance density
$\sigma^2(\xi_{f_k}^{(\infty)})$ where $\tau:= (d-1)/(d+1+2\delta)$.
In particular, the total number $f_k(K_{\la})$ of
$k$-faces for $K_{\la}$ satisfies the scalar variance asymptotics and
central limit theorem with
scaling exponent $\tau/2$ and variance $\sigma_{f_k}^2 :=
\sigma^2(\xi_{f_k}^{(\infty)})\,d\kappa_d.$
\end{theo}

\begin{remarks*} (i) Setting $\delta= 0$ in Theorem~\ref{Vkface}
gives the asserted variance limit (\ref{introvar}).

(ii) We expect that the variance asymptotics of Theorems
\ref{VWSCALAR} and~\ref{Vkface} can be de-Poissonized, that is to
say, that there are analogous variance limits when the polytope
$K_\la$ is replaced by the polytope $K_n$ generated by $n$ i.i.d.
uniformly distributed points in $\B^d$. We leave these issues for
further study.
\end{remarks*}

\section{Global regime and Brownian limits}\label{sec8}\label{InvPrinc}
In this section we establish a functional central limit theorem for
the integrated convex hull processes $\hat{W}_{\la}$
and $\hat{V}_{\la}$, defined at (\ref{SGRESC}). The methods extend
to yield functional
central limit theorems for stabilizing functionals in general,
thus extending \cite{Sh}.

For any $\sigma^2 > 0$ let $B^{\sigma^2}$ be the Brownian sheet of variance
coefficient $\sigma^2$ on the injectivity region $\B_{d-1}(\pi)$ of
$\exp:=\exp_{\S^{d-1}};$ that
is to say, $B^{\sigma^2}$ is the mean zero continuous path
Gaussian process indexed by $\R^{d-1}$ with
\[
\operatorname{Cov}(B^{\sigma^2}(v),B^{\sigma^2}(w)) = \sigma^2
\cdot\sigma
_{d-1}\bigl(\exp([\0,v] \cap[\0,w])\bigr),
\]
%
where, recall, $\sigma_{d-1}$ is the $(d-1)$-dimensional surface
measure on $\S^{d-1}.$
Recalling from Lemma~\ref{AVARREPR} the shorthand notation
$\sigma^2_s:= \sigma^2(\xi_s^{(\infty)})$ and
$\sigma^2_r:=\sigma^2(\xi_r^{(\infty)})$,
we have the following limit result, the main result of this section.
Via Lemma
\ref{PClem} and Theorem~\ref{VWSCALAR}, this theorem also yields
Brownian sheet limits for the defect width and volume functionals
given at \eqref{HSV} and \eqref{HRV}, respectively. We remark that
the same holds for the two processes $\Vol([\B^d\setminus K_{\la
}]\cap\cone(\exp([\0,v])))$ and $\Vol([\B^d\setminus F({\mathcal
P}_{\la})]\cap\cone(\exp([\0,v])))$, $v\in\R^{d-1}$.

\begin{theo}\label{GLOBALLT}
As $\la\to\infty$, the random functions $\hat{W}_{\la} \dvtx  \R^{d-1}
\to\R$
converge in law to $B^{\sigma^2_s}$ in the space ${\cal C}(\R^{d-1}).$
Likewise, the random functions $\hat{V}_{\la} \dvtx  \R^{d-1} \to\R$
converge in law
to $B^{\sigma^2_r}$ in ${\cal C}(\R^{d-1}).$
\end{theo}

\begin{pf*}{Proof of Theorem~\ref{GLOBALLT}} Our argument relies heavily on
the theory developed in \cite{SY}
and is further extended in Section~\ref{StabRepr}. For $v \in\R
^{d-1}$ and $x \in\B^d$, define
%
\begin{equation}\label{INDIC}
{\bf1}^{[\0,v]}_{\B^d}(x) :=
\cases{\displaystyle 1  ,&\quad   if
$x/|x| \in\exp([\0,v])$,\cr\displaystyle
0  ,&\quad   otherwise.
}
\end{equation}

We thus have 
the identities
\[
\la^{\zeta/2}\bigl( H^{\xi_s}_{\la}(v) - \E H^{\xi_s}_{\la}(v)\bigr) =
\la^{\zeta/2} \bigl\langle{\bf1}^{[\0,v]}_{\B^d}, \bar{\mu}^{\xi
_s}_{\la} \bigr\rangle
\]
and
\[
\la^{\zeta/2}\bigl( H^{\xi_r}_{\la}(v) - \E H^{\xi_r}_{\la}(v)\bigr) =
\la^{\zeta/2}\bigl\langle{\bf1}^{[\0,v]}_{\B^d}, \bar{\mu}^{\xi
_r}_{\la} \bigr\rangle.
\]
Recalling from (\ref{SGRESC}) that $ \hat{W}_{\la}(v) :=
\la^{\zeta/2} (W_{\la}(v) - \E W_{\la}(v))$
and $ \hat{V}_{\la}(v) := \la^{\zeta/2} (V_{\la}(v) - \E
V_{\la}(v))$, and using \eqref{Diff1} and \eqref{Diff2} from Lemma
\ref{PClem}, we obtain,
uniformly in $v$,
%
\begin{eqnarray}\label{NASYMPTEQ}
  \lim_{\la\to\infty} \bigl| \hat{W}_{\la}(v) - \la^{\zeta/2}
\bigl\langle{\bf1}^{[\0,v]}_{\B^d}, \bar\mu^{\xi_s}_{\la}\bigr \rangle\bigr|
&\stackrel{{P}}{=}&0,\nonumber
\\[-8pt]
\\[-8pt]
\lim_{\la\to\infty} \bigl| \hat{V}_{\la}(v) - \la^{\zeta/2}
\bigl\langle{\bf1}^{[\0,v]}_{\B^d}, \bar\mu^{\xi_r}_{\la}\bigr \rangle\bigr|
&\stackrel{{P}}{=}&
0.
\nonumber
\end{eqnarray}
%


Even though ${\bf1}^{[\0,v]}_{\B^d}$ is not a continuous function, it
is easily seen that the proofs in \cite{SY}
hold for functions which are almost everywhere continuous with respect
to the uniform
measure on $\B^d$, and, in fact, the central limit
theorems and variance asymptotics of \cite{SY} hold for all bounded
functions on $\B^d$. Thus
Theorem~\ref{MainCLT}
for $\xi_s$ and $\xi_r$ remain valid upon setting the test function
$g$ to ${\bf1}^{[\0,v]}_{\B^d}.$ This application of Theorem
\ref{MainCLT},
combined with (\ref{NASYMPTEQ}), yields that the fidis of
$(\hat{W}_{\la}(v))_{v \in\R^{d-1}}$ converge
to those of $(B^{\sigma^2(\xi_s^{(\infty)})}(v))_{v \in\R^{d-1}}$
and, likewise, the fidis of $(\hat{V}_{\la}(v))_{v \in\R^{d-1}}$
converge to those of $(B^{\sigma^2(\xi_r^{(\infty)})}(v))_{v \in\R^{d-1}}.$
Additionally, for all $v \in\R^{d-1}$, we have
\[
\lim_{\la\to\infty} \Var[\hat{W}_{\la}(v)] = \sigma^2\bigl(\xi
_s^{(\infty)}\bigr)(v),
\]
with similar variance asymptotics for $\hat{V}_{\la}(v)$; see also Theorems
1.2 and 1.3 in \cite{SY}.

We claim that the fidis convergence of $\hat{W}_{\la}$ and $\hat
{V}_{\la}$
can be strengthened to convergence in law in ${\cal C}(\R^{d-1}).$
It suffices to establish the
tightness of the processes $(\hat{W}_{\la}(v))_{v \in\R^{d-1}}$
and $(\hat{V}_{\la}(v))_{v \in\R^{d-1}}.$ We shall
focus on $\hat{W}_{\la},$ the argument for $\hat{V}_{\la}$ being
analogous, and we shall proceed to some extent along the lines of the
proof of Theorem 8.2 in \cite{DD},
which is based on \cite{BW}.
We extend the definition of $W_{\la}$ to subsets of $\R^{d-1}$ putting
for measurable $B \subseteq\R^{d-1}$
\[
W_{\la}(B) := \int_{\exp_{d-1}(B)} s_{\la}(u)\,d\sigma_{d-1}(u)
\]
%
and letting
%
\begin{equation}\label{WHAT}
\hat{W}_{\la}(B) := \la^{\zeta/2}\bigl(W_{\la}(B) - \E W_{\la}(B)\bigr).
\end{equation}
%
It is enough to show
%
\begin{equation}\label{BillMomBds}
\E (\hat{W}_{\la}([v,v']) )^4
= O(\Vol([v,v'])^2), \qquad   v, v' \in\R^{d-1},
\end{equation}
for then $\hat{W}_{\la}$ satisfies condition (2) on page 1658 of
\cite{BW}, thus
belongs to the class ${\cal C}(2,4)$ of \cite{BW} and is tight in view
of Theorem 3
on page 1665 of \cite{BW}.

To this end, we put
%
\begin{equation}\label{WTILDE}
{W}^{\#}_{\la}(B) := \la^{\eta[\xi_s]} W_{\la}(B) = \la^{\beta
(d-1)+\gamma}
W_{\la}(B),
\end{equation}
where we recall from the definition of $\Xi^{(\la)}$ in Section \ref
{StabRepr}
that $\eta[\xi_s] = \beta(d-1)+\gamma$ is the proper scaling exponent
for $\xi_s.$ The crucial point now is that in analogy to the proof of Lemma
5.3 in \cite{BY2}, and similar to (3.24) in the proof of Theorem~1.3
in \cite{SY}, by a stabilization-based argument all cumulants of
${W}^{\#}_{\la}([v,w])$ over rectangles $[v,w]$
are at most linear in $\la^{\tau} \Vol([v,w])$ with $\tau:= \beta
(d-1)$ as in
(\ref{TAUDEF}). In other words, for all $k \geq1,$ we have
%
\begin{equation}\label{CumBds}
 |c^k({W}^{\#}_{\la}([v,w])) |
\leq C_k \la^{\tau} \Vol([v,w]), \qquad   v,w \in\R^{d-1},
\end{equation}
where $c^k(Y)$ stands for the $k$th order cumulant of the random
variable $Y$ and where $C_k$
is a constant.
Thus, putting (\ref{WHAT}) and (\ref{WTILDE}) together, we get from
(\ref{CumBds})
%
\begin{eqnarray}\label{CumBds2}
 |c^k(\hat{W}_{\la}([v,w])) | &\leq& C_k
\la^{k[\zeta/2-\eta[\xi_s]]} \la^{\tau} \Vol([v,w])
\nonumber
\\[-8pt]
\\[-8pt]
&=&C_k \la^{k[\zeta/2-\beta(d-1)-\gamma]} \la^{\beta(d-1)} \Vol
([v,w]).
\nonumber
\end{eqnarray}
To proceed, we use the identity
$\E(Y-\E Y)^4 = c^4(Y) + 3 (c^2(Y))^2$ valid for any random variable $Y.$
Recalling that $\gamma= 2 \beta$ and $\zeta= \beta(d-1) + 2
\gamma$, as in \eqref{BETAGAMMA} and \eqref{ZETA}, respectively, we
obtain from \eqref{BillMomBds}--\eqref{CumBds2} that for $v,w \in
\R^{d-1}$,
%
\begin{eqnarray}\label{MomBds}
\E (\hat{W}_{\la}([v,w]) )^4
&=& O\bigl(\la^{4[\zeta/2-\beta(d-1)-\gamma]} \la^{\beta(d-1)} \Vol([v,w])\bigr)
\nonumber\\
&&{}+ O\bigl(\bigl[\la^{2[\zeta/2-\beta(d-1)-\gamma]} \la^{\beta(d-1)} \Vol([v,w])\bigr]^2\bigr)
\\
&=& O\bigl(\la^{-\beta(d-1)} \Vol([v,w])\bigr) + O(\Vol([v,w])^2),
\nonumber
\end{eqnarray}
which is of the required order $O(\Vol([v,w])^2)$ as soon as
$\Vol([v,w]) =  \break \Omega(\la^{-\beta(d-1)}).$ Thus we have
shown (\ref{BillMomBds}) for $\Vol([v,w]) = \Omega(\la^{-\beta(d-1)})$,
and we have to show it holds for $\Vol([v,w]) = O(\la^{-\beta(d-1)})$
as well.
To this end, we use that $W_{\la}([v,w]) = \la^{-\gamma} O_P(\Vol([v,w]))$
with $\gamma$ being the height coordinate re-scaling exponent, and
$\E[W_{\la}([v,w])-\E W_{\la}([v,w])]^4 = \la^{-4\gamma}
O(\Vol([v,w])^4)$.
Thus by (\ref{WHAT})
\[
\E (\hat{W}_{\la}([v,w]) )^4 = \la^{2\zeta} \la
^{-4\gamma} O(\Vol([v,w])^4).
\]

Recalling $\zeta= \beta(d-1) + 2\gamma$
and using that $\Vol([v,w]) = O(\la^{-\beta(d-1)})$,
we conclude that
\[
\E (\hat{W}_{\la}([v,w]) )^4 = O\bigl(\la^{2\beta(d-1)} \Vol
([v,w])^4\bigr) = O(\Vol([v,w]^2))
\]
%
as required, which completes the proof of the required relation (\ref
{BillMomBds}).
Having obtained the required tightness, we get the convergence in law
of $(\hat{W}_{\la}(v))_{v \in\R^{d-1}}$
to $(B^{\sigma^2_s}(v))_{v \in\R^{d-1}}$ and,
likewise, of $(\hat{V}_{\la}(v))_{v \in\R^{d-1}}$
to $(B^{\sigma^2_r}(v))_{v \in\R^{d-1}}$in
${\cal C}(\R^{d-1})$. This completes the proof of Theorem
\ref{GLOBALLT}.
\end{pf*}

\section*{Acknowledgments}
It is a pleasure to thank an anonymous referee for a careful reading of the
original manuscript, resulting in an improved and more accurate exposition.


%

\printaddresses

\end{document}